\documentclass[11pt,letterpaper,twoside,reqno,nosumlimits]{amsart}

\usepackage[usenames,dvipsnames]{xcolor}
\usepackage{fancyhdr}
\usepackage{amsmath,amsfonts,amsbsy,amsgen,amscd,amssymb,amsthm}
\usepackage{url}

\usepackage{mathtools}

\mathtoolsset{showonlyrefs}

\usepackage[font=small,margin=0.25in,labelfont={bf},labelsep={colon}]{caption}

\usepackage{subcaption}

\usepackage{tikz}
\usepackage{microtype}
\usepackage{enumitem}

\definecolor{dark-gray}{gray}{0.3}
\definecolor{dkgray}{rgb}{.4,.4,.4}
\definecolor{dkblue}{rgb}{0,0,.5}
\definecolor{medblue}{rgb}{0,0,.75}
\definecolor{rust}{rgb}{0.5,0.1,0.1}

\usepackage[colorlinks=true]{hyperref}

\hypersetup{urlcolor=rust}
\hypersetup{citecolor=dkblue}
\hypersetup{linkcolor=dkblue}

\usepackage{setspace}

\usepackage{graphicx}
\usepackage{booktabs,longtable,tabu} 
\usepackage{multirow} 

\usepackage{float}

\usepackage[full]{textcomp}

\usepackage[scaled=.98,sups,osf]{XCharter}
\usepackage[scaled=1.04,varqu,varl]{inconsolata}
\usepackage[type1]{cabin}
\usepackage[charter,vvarbb,scaled=1.07]{newtxmath}
\usepackage[cal=boondoxo]{mathalfa}
\usepackage[T1]{fontenc}

\usepackage{bm} 

\graphicspath{{figures/}}

\newtheorem{theorem}{Theorem}[section]
\newtheorem{lemma}[theorem]{Lemma}

\newtheorem{proposition}[theorem]{Proposition}

\newtheorem{corollary}[theorem]{Corollary}

\theoremstyle{definition}

\newtheorem{remark}[theorem]{Remark}

\numberwithin{equation}{section} 
\numberwithin{figure}{section}
\numberwithin{table}{section}

\floatstyle{plaintop}
\newfloat{recipe}{thp}{lor}
\floatname{recipe}{Recipe}
\numberwithin{recipe}{section}

\renewcommand{\phi}{\varphi}

\newcommand{\econst}{\mathrm{e}}

\providecommand{\mathbbm}{\mathbb} 

\newcommand{\R}{\mathbbm{R}}

\newcommand{\Z}{\mathbbm{Z}}

\newcommand{\sgn}[1]{\operatorname{sgn}{#1}}

\newcommand{\diff}{\mathrm{d}}
\newcommand{\idiff}{\, \diff}

\newcommand{\triplenorm}[1]{{\left\vert\kern-0.25ex\left\vert\kern-0.25ex\left\vert #1
    \right\vert\kern-0.25ex\right\vert\kern-0.25ex\right\vert}}

\newcommand{\calG}{\mathcal{G}}

\newcommand{\calQ}{\mathcal{Q}}
\newcommand{\calR}{\mathcal{R}}

\newcommand{\calT}{\mathcal{T}}

\newcommand{\calX}{\mathcal{X}}

\def\XXint#1#2#3{{%
\setbox0=\hbox{$#1{#2#3}{\int}$}
\vcenter{\hbox{$#2#3$}}\kern-.5\wd0}}

\renewcommand{\leq}{\leqslant}
\renewcommand{\geq}{\geqslant}
\renewcommand{\subset}{\subseteq}

\newcommand{\om}{\omega}

\usepackage{mathrsfs}
\providecommand{\mathbbm}[1]{\mathbb{#1}}
\providecommand{\sgn}{\operatorname{sgn}}
\newcommand{\dr}[1]{\left(#1\right)}
\allowdisplaybreaks
\mathtoolsset{showonlyrefs}
\hypersetup{bookmarksdepth=2,pdftitle={Steady states and rescaled profiles for periodic gCLM}}
\begin{document}
\title[Steady states and profiles for periodic gCLM]{The Periodic Generalized Constantin--Lax--Majda Equation I:\\ Steady States and Rescaled Profiles}
\author{Thomas Y. Hou}
\address{Applied and Computational Mathematics,
California Institute of Technology, Pasadena, CA, USA}
\email{hou@cms.caltech.edu}

\author{Xiang Qin}
\address{Applied and Computational Mathematics,
California Institute of Technology, Pasadena, CA, USA}
\email{xqin2@caltech.edu}

\author{Peicong Song}
\address{Applied and Computational Mathematics,
California Institute of Technology, Pasadena, CA, USA}
\email{psong2@caltech.edu}
\date{}
\begin{abstract}
We construct periodic steady states for every \(a>0\) and rescaled
profiles for every \(0<a<1\) in the generalized Constantin--Lax--Majda
equation. A single compact fixed-point map selects a finite-period or
whole-line profile with nonnegative spatial rate. Our estimates guarantee
focusing profiles for \(0<a\leq(17+16\log2)/(53-16\log2)\) and
finite-period profiles for \(64/(176-9\pi^2)\leq a<1\), while connected
continuation yields a critical whole-line profile at a parameter strictly
between these two bounds. We determine sharp vorticity regularity at the
steady-state zeros and periodic endpoints, together with algebraic or
Gaussian behavior at infinity. All constructions are analytic, using
variational and compact fixed-point methods and a sharp first-moment
inequality.
\end{abstract}
\maketitle
\tableofcontents

\section{Introduction}
\label{sec:introduction}

We consider the generalized Constantin--Lax--Majda equation on the circle,
\begin{equation}\label{eq:periodic_gclm}
 \partial_t\omega+a u\partial_x\omega=(\partial_xu)\omega,
 \qquad \partial_xu=H_{\mathbb T}\omega,
 \qquad u(0,t)=0,
 \qquad \mathbb T:=\mathbb R/(2\pi\mathbb Z).
\end{equation}
The parameter \(a\) measures the strength of transport relative to
stretching. We work with odd vorticity and the corresponding odd velocity.
Our convention for the periodic Hilbert transform is
\begin{equation}\label{eq:hilbert_convention}
 H_{\mathbb T}v(x)
 =\frac1{2\pi}\operatorname{p.v.}\int_{-\pi}^{\pi}
 v(y)\cot\left(\frac{x-y}{2}\right)\idiff y.
\end{equation}
The original model of Constantin, Lax, and Majda corresponds to \(a=0\),
and the De Gregorio equation corresponds to \(a=1\)
\cite{CLM1985,DG1990,DG1996,OSW2008}.

Our main motivation is to understand the spatial shapes of self-similar
and nearly self-similar blowup in the periodic problem. We therefore seek
stationary solutions after rescaling the growing vorticity. In this
question the spatial domain is part of the problem: concentration
can send the rescaled period to infinity, whereas a bounded rescaled
period retains the periodic Biot--Savart law. A fixed-point construction
must therefore select both a profile and its spatial geometry. As a
complementary problem, we also construct steady states of the original
equation \eqref{eq:periodic_gclm}; these are time-independent physical
solutions, without rescaling.

This paper develops that construction and the regularity theory of its
solutions. We obtain periodic steady states for every \(a>0\), and, for
every \(0<a<1\), a normalized rescaled profile on either a finite circle
or the real line. The selected profiles have a precise leading behavior
at the identified endpoints of a period or at spatial infinity. We also
prove the existence of a critical whole-line profile with zero spatial
scaling rate, and sharpen the analytic lower bound for its parameter.
The proofs use variational compactness, Schauder's fixed-point theorem,
and continuation of compact fixed-point families.

\subsection{Stationary equations and main results}

The steady-state and blowup-profile problems lead to different stationary
equations. A physical steady state satisfies
\(a u\omega_x=u_x\omega\), interpreted in distributions when necessary.
For blowup, the amplitude and spatial scales vary with time, and it is
the rescaled shape that is stationary. We first derive the latter equation
and then state the two existence results separately.

For \(L\in(0,\infty)\), write
\(\mathbb T_L=\mathbb R/(2L\mathbb Z)\), and define
\begin{equation}\label{eq:H_L_intro}
 H_Lv(x)=\frac1{2L}\operatorname{p.v.}\int_{-L}^{L}
 v(y)\cot\left(\frac{\pi(x-y)}{2L}\right)\idiff y.
\end{equation}
The notation \(H_\infty=H_{\mathbb R}\) refers to the whole-line transform
with kernel \(1/(\pi(x-y))\). To seek the limiting shape of a self-similar
or nearly self-similar singularity, introduce a rescaled time \(\tau\) and positive
scale factors \(C_l(\tau),C_\omega(\tau)\), with
\(\idiff t/\idiff\tau=C_\omega^{-1}\). The change of variables
\begin{equation}\label{eq:intro_dynamic_rescaling}
\begin{gathered}
 \Omega(x,\tau)=C_\omega^{-1}\omega(C_l^{-1}x,t(\tau)),
 \qquad U(x,\tau)=C_lC_\omega^{-1}u(C_l^{-1}x,t(\tau)),\\
 \partial_\tau\Omega+(c_lx+aU)\partial_x\Omega
   =(c_\omega+\partial_xU)\Omega,
 \qquad \partial_xU=H_{\pi C_l}\Omega
\end{gathered}
\end{equation}
gives the normalized evolution, where
\(c_l=\idiff(\log C_l)/\idiff\tau\) and
\(c_\omega=-\idiff(\log C_\omega)/\idiff\tau\).
Thus the vorticity amplitude is \(C_\omega\), the physical length scale
is \(C_l^{-1}\), and the rescaled half-period is \(\pi C_l\).
The stationary equation with limiting rates
\(\bar c_l,\bar c_\omega\) and half-period \(L\) leads us to define
a rescaled profile as a tuple
\((\bar\Omega,\bar U,\bar c_l,\bar c_\omega,L)\) satisfying
\begin{equation}\label{eqt:general_rescaled_profile_equation_periodic}
 (\bar c_lx+a\bar U)\partial_x\bar\Omega
   =(\bar c_\omega+\partial_x\bar U)\bar\Omega,
 \qquad \partial_x\bar U=H_L\bar\Omega,
 \qquad \bar U(0)=0.
\end{equation}
The constants \(\bar c_l\) and \(\bar c_\omega\) are therefore the
limiting spatial and amplitude scaling rates. A negative
\(\bar c_\omega\) gives increasing physical amplitude. If
\(\bar c_l>0\), the physical length scale shrinks exponentially in
\(\tau\), and the rescaled period tends to infinity; this is the
\emph{focusing} case. When \(\bar c_l=0\), the rescaled period may
still diverge, giving the \emph{critical} whole-line case, or approach
a finite limit, giving the \emph{finite-period} case. We impose
\(\bar c_l\geq0\), \(\bar c_\omega<0\), and
\begin{equation}\label{eq:profile_geometry_normalization}
 L<\infty\ \Longrightarrow\ \bar c_l=0,
 \qquad
 \bar\Omega(x)=-x+x^3+O(x^5)\quad(x\to0).
\end{equation}
The word critical refers throughout to the zero spatial rate in this
whole-line stationary equation. The two origin coefficients in
\eqref{eq:profile_geometry_normalization} fix the amplitude and length
normalizations, as explained in Section~\ref{sec:period_selection}.

We begin with the complementary equilibrium result for the original
periodic equation. The functions in the next theorem are independent
of physical time.

\begin{theorem}[Periodic steady states]\label{thm:intro_steady}
For every \(a>0\), there are a nonzero odd periodic function \(u\) and a
constant \(\lambda>0\) such that
\begin{equation}\label{eq:intro_steady_pair}
 \omega=-\lambda |u|^{1/a-1}u,
 \qquad u_x=H_{\mathbb T}\omega,
 \qquad a(u\omega)_x=(1+a)u_x\omega.
\end{equation}
The last identity holds in distributions. The velocity is positive on
\((0,\pi)\), its only zeros are \(0\) and \(\pi\) modulo \(2\pi\),
and both zeros are simple. The pair is smooth away from these zeros.
Writing \(p=1/a\), its global regularity is as follows:
\begin{enumerate}[label=\emph{(\roman*)}]
\item If \(p\) is an odd positive integer, then \(u,\omega\in C^\infty\).
\item If \(p=m+\alpha\notin\mathbb N\), with \(m\in\mathbb N_0\)
and \(0<\alpha<1\), then
\(u\in C^{m+1,\alpha}\) and \(\omega\in C^{m,\alpha}\).
\item If \(p=2N\), with \(N\geq1\), then
\(u\in C^{2N,\alpha}\) for every \(0<\alpha<1\) and
\(\omega\in C^{2N-1,1}\).
\end{enumerate}
For \(0<a\leq1\), the steady-state equation holds pointwise on the circle.
\end{theorem}

The simple zeros determine the exact leading vanishing order
\(|\omega(x)|\asymp\operatorname{dist}(x,\{0,\pi\})^{1/a}\).
Theorem~\ref{thm: steady state} and
Proposition~\ref{prop:steady_endpoint_behavior} prove the assertions,
including the sharp noninteger and even-integer regularity thresholds
for the vorticity.

We next construct stationary shapes of the rescaled evolution
\eqref{eq:intro_dynamic_rescaling}.

\begin{theorem}[Rescaled profiles]\label{thm:intro_profiles}
For every \(0<a<1\), there is a nontrivial odd profile satisfying
\eqref{eqt:general_rescaled_profile_equation_periodic}--
\eqref{eq:profile_geometry_normalization} for some \(L\in(1,\infty]\).
On the positive half of its domain, the function
\begin{equation}\label{eq:intro_f_definition}
 f(x):=-\frac{\bar\Omega(x)}x
\end{equation}
is positive and nonincreasing, and \(s\mapsto f(\sqrt{s})\) is convex.
The profile is smooth at every interior point, including the origin.
Its endpoint behavior is given by the following alternatives.
\begin{enumerate}[label=\emph{(\roman*)}]
\item If \(L<\infty\), then \(\bar\Omega\in C^1(\mathbb T_L)\),
and there are \(A>0\) and
\(p_L\geq(3-a)/(2a)>1\) such that
\begin{equation}\label{eq:intro_finite_asymptotic}
 f(L-t)=At^{p_L}(1+o(1)),
 \qquad f'(L-t)=-p_LAt^{p_L-1}(1+o(1))
 \quad(t\downarrow0).
\end{equation}
\item If \(L=\infty\) and \(\bar c_l>0\), then, for some \(A>0\),
\begin{equation}\label{eq:intro_focusing_asymptotic}
 f(x)\sim Ax^{-\rho},
 \qquad \rho=1-\frac{\bar c_\omega}{\bar c_l}>2,
 \qquad x\to\infty.
\end{equation}
Every ordinary derivative of \(f\) is \(O((1+|x|)^{-\rho})\).
\item If \(L=\infty\) and \(\bar c_l=0\), then
\(f,\bar\Omega\in\mathcal S(\mathbb R)\), and
\begin{equation}\label{eq:intro_critical_asymptotic}
 \lim_{x\to\infty}\frac{\log f(x)}{x^2}
 =-\frac{(1-a)b}{2a\mathsf b_2},
 \qquad \bar U(x)\sim\frac{\mathsf b_2}{x},
 \qquad
 b=\frac2\pi\int_0^\infty f\idiff x,
 \quad \mathsf b_2=\frac2\pi\int_0^\infty x^2f\idiff x.
\end{equation}
\end{enumerate}
\end{theorem}

Here \(\mathcal S(\mathbb R)\) denotes the Schwartz class. In particular,
rapid decay of the critical vorticity coexists with an algebraic velocity
tail. The existence assertion follows from
Theorems~\ref{thm:fixed_point_existence} and
\ref{thm:fixed_point_profile_reconstruction};
Section~\ref{sec:regularity} proves the stated regularity and asymptotics.

The endpoint power in the finite-period case determines a sharp
regularity classification, proved in
Theorem~\ref{thm:periodic-endpoint-regularity}. A noninteger power
\(p_L=m+\alpha\) gives precisely \(C^{m,\alpha}\) vorticity at the
joining point; an even integer \(p_L=2N\) gives
\(C^{2N-1,1}\setminus C^{2N}\); and an odd integer gives a smooth
physical profile. These assertions concern \(\bar\Omega\) on
\(\mathbb T_L\); the scalar zero extension used in the fixed-point
map is a different function. In the critical case,
Proposition~\ref{prop:critical-full-tail} sharpens
\eqref{eq:intro_critical_asymptotic} to a Gaussian asymptotic with an
explicit algebraic prefactor and a relative \(O(x^{-2})\) error.

The parameter bounds are expressed by the two explicit constants
\begin{equation}\label{eq:intro_parameter_constants}
 \begin{aligned}
 a_{\mathrm{foc}}&:=\frac{17+16\log2}{53-16\log2}
       =0.670259908298743\ldots,\\
 a_{\mathrm{per}}&:=\frac{64}{176-9\pi^2}
       =0.7341675585\ldots.
 \end{aligned}
\end{equation}

\begin{theorem}[Focusing range and a critical profile]
\label{thm:intro_critical}
The following assertions hold for the normalized fixed-point construction.
\begin{enumerate}[label=\emph{(\roman*)}]
\item For every \(0<a\leq a_{\mathrm{foc}}\), there is a focusing
whole-line profile. Every HQWW fixed point at such a parameter is focusing.
\item For every \(a_{\mathrm{per}}\leq a<1\), every fixed point of the
period-selection map has finite selected period.
\item There exist \(a_*\in(a_{\mathrm{foc}},a_{\mathrm{per}})\)
and a critical whole-line profile at \(a=a_*\).
\end{enumerate}
\end{theorem}

The whole-line profiles in this theorem are precisely profiles in the
normalization of Huang, Qin, Wang, and Wei \cite{HQWW2024}, abbreviated
HQWW, with nonnegative spatial rate. Proposition~\ref{prop:hqww_correspondence}
makes the correspondence explicit. The critical existence assertion
comes from a zero of the spatial rate along a connected set of fixed points.
The parameter estimates classify whole-line fixed points; finite-period
fixed points are treated by the period-selection map of
Section~\ref{sec:period_selection}.
The constants in \eqref{eq:intro_parameter_constants} are sufficient
parameter bounds. The sharpness of the first-moment estimate refers to
its optimal constant on the admissible shape class.

\subsection{Related work}

The gCLM equation isolates the competition between transport and vortex
stretching in the vorticity formulation of the three-dimensional Euler
equations.  The continuation criterion of Beale, Kato, and Majda
\cite{BKM1984} identifies the time integral of the maximum vorticity as a
quantity controlling classical continuation, while the geometric criteria
of Constantin, Fefferman, and Majda \cite{CFM1996} explain how coherence of
the vorticity direction can weaken stretching.  The role of convection and
dynamic depletion was further investigated by Okamoto and Ohkitani
\cite{OO2005}, Hou and Li \cite{HouLi2006,HouLi2008}, and Hou and Lei
\cite{HouLei2009}.  Although the gCLM model has no incompressibility
constraint or vorticity direction, it retains a nonlocal stretching rate
and permits the relative strength of transport to vary through one
parameter.  In two dimensions, the double-exponential small-scale creation
of Kiselev and \v Sver\'ak \cite{KS2014} provides a complementary example
in which transport produces rapid growth while the vorticity remains
bounded.

Self-similar constructions for Euler provide a second source of motivation.
Elgindi \cite{Elgindi2021} constructed finite-time singularities with
\(C^{1,\alpha}\) velocity, and Elgindi, Ghoul, and Masmoudi
\cite{EGM2021Euler} established stability of the corresponding mechanism.
For smooth data in a domain with boundary, the numerical scenario of Luo
and Hou \cite{LH2014PNAS,LH2014MMS} led to the analytic and computational
work of Chen and Hou \cite{CH2022Analysis,CH2025Numerics};
\cite{CH2025PNAS} presents the resulting boundary singularity construction.
The one-dimensional Hou--Luo model has also been treated through integral
estimates \cite{Choi2017} and, more recently, through a purely analytic
fixed-point construction of smooth profiles \cite{HQWW2025HL}.
Shao, Wei, Zhang, and Zhang \cite{SWZZ2026} construct self-similar profiles
for higher-dimensional axisymmetric Euler without swirl and prove a
finite-codimensional stability theorem.  These works distinguish two
steps that are also essential here: constructing a stationary solution of
the rescaled equation and realizing it as the asymptotic profile of an
evolving solution in the desired spatial domain and regularity class.

Constantin, Lax, and Majda \cite{CLM1985} introduced the model without
transport and obtained finite-time singularity by an explicit solution
formula.  De Gregorio \cite{DG1990,DG1996} added transport with the same
coefficient as stretching.  Okamoto, Sakajo, and Wunsch \cite{OSW2008}
introduced the parameterized family considered here, and Wunsch
\cite{Wunsch2011} developed further aspects of its analysis.
Periodic equilibria and traveling waves were studied by Okamoto, Sakajo,
and Wunsch \cite{OSW2014}.  The related nonlocal scalar transport equations
of C\'ordoba, C\'ordoba, and Fontelos \cite{CCF2005}, and Castro and
C\'ordoba \cite{CastroCordoba2010}, show how transport can also reinforce
singularity formation; Castro's thesis \cite{Castro2010Thesis} develops
these nonlocal fluid models in a broader setting.

The construction of smooth self-similar gCLM profiles began with the
perturbative analysis of Elgindi and Jeong \cite{EJ2020} near the explicit
CLM solution.  Elgindi, Ghoul, and Masmoudi \cite{EGM2021Transport} proved
stability in a corresponding small-parameter regime.  Chen, Hou, and
Huang \cite{CHH2021} developed a dynamic-rescaling and weighted-energy
method that also yields smooth blowup for the De Gregorio equation on
the real line.  The explicit profile at \(a=1/2\), found independently in
\cite{Chen2020,LSS2021}, provides another parameter around which stable
smooth blowup can be analyzed.  Huang, Tong, and Wei \cite{HTW2023}
subsequently constructed infinitely many compactly supported De Gregorio
profiles through a compact self-adjoint operator.  Exactly self-similar
solutions of lower H\"older regularity, for a different range of parameters
and regularity exponents, were constructed in \cite{EJ2020,Zheng2023}.
Recent work of Hou, Qin, and Song \cite{HyperConvective} extends
smooth-data blowup into the range \(a>1\) through a stable local branch
of compactly supported profiles near \(a=3/2\). For parameters sufficiently
close to \(3/2\), this branch gives spatially expanding blowup below
\(3/2\) and global decay above it for nearby smooth compactly supported data.

Huang, Qin, Wang, and Wei \cite{HQWW2024}, abbreviated HQWW, gave a unified
fixed-point construction for every \(a\leq1\).  In our notation, their normalized variable
\(f=-\bar\Omega/x\) is nonnegative, nonincreasing, and convex as a function of
\(x^2\).  These shape properties yield compactness and preserve a convex
set under the nonlinear reconstruction map.  The resulting profiles are
smooth in the interior of their supports, and their spatial scaling rate
distinguishes compact support, Gaussian logarithmic decay, and algebraic
decay.  The whole-line profiles used in the present work belong to this
family.  We use its established existence and regularity theory and make
the correspondence of the two normalizations explicit.  The new
first-moment inequality sharpens the focusing range within this class;
its equality case is determined by a single parabolic cap.
The additional ingredients here are selection of a finite period within
the fixed-point map, reconstruction and sharp regularity at the periodic
joining point, the improved focusing bound, and connected continuation
to a critical profile. We also refine the critical whole-line tail to
include its algebraic prefactor.

On the circle, the domain geometry introduces phenomena that are absent
from the whole-line profile equation.  Jia, Stewart, and \v Sver\'ak
\cite{JSS2019} analyzed the stability of sinusoidal De Gregorio equilibria.
Lei, Liu, and Ren \cite{LLR2020} established global regularity for
sign-definite vorticity and developed weighted estimates near an
equilibrium.  Chen \cite{Chen2021} studied the slightly perturbed
De Gregorio equation on the circle, and \cite{Chen2025Regularity} obtained
further global regularity and blowup criteria in an odd, one-sided-sign
class.  The numerical study \cite{LSS2021} identified different focusing
and finite-period blowup regimes and a transition near
\(a=0.6890665\). The companion work \cite{CompanionDynamics} proves
global regularity under sign, endpoint, and interior-zero hypotheses,
gives a blowup criterion throughout \(0<a<1\), and realizes each focusing
or finite-period profile constructed here by smooth periodic blowup.
Its logarithmic concentration result uses a separately constructed critical
profile with computer-assisted estimates; the qualitative critical existence
and parameter bounds in Theorem~\ref{thm:intro_critical} are proved analytically.

The recent periodic analysis of Guo and Jiu \cite{GuoJiu2025} treats
stability and instability near the first excited sinusoidal state.
Silantyev, Lushnikov, Siegel, and Ambrose \cite{SLAS2025} derive exact
periodic pole solutions for dissipative gCLM equations.
At \(a=0\), Huang, Qin, and Wang \cite{HQW2025Multiscale} constructed
asymptotically self-similar blowups with two collapsing spatial scales,
using the explicit CLM solution formula and its complex pole dynamics.
Huang, Tong, and Wang \cite{HTW2026Singular} study a different normalization
with derivative degeneracy at the origin.  Their work combines explicit
singular profiles and fixed-point constructions of traveling waves with
numerical investigations of additional one-scale and two-scale regimes.
Our normalization instead retains a nonzero vorticity slope at the origin.

Dissipative variants have a related, but distinct, regularity theory.
Schochet \cite{Schochet1986} found explicit viscous-model solutions.
Transport equations with nonlocal drift and dissipation were studied by
Dong \cite{Dong2008}, Li and Rodrigo \cite{LiRodrigo2008}, and Silvestre
and Vicol \cite{SV2016}; Kiselev \cite{Kiselev2010} surveys the surrounding
active-scalar theory.  Chen \cite{Chen2020} and Ambrose, Lushnikov, Siegel,
and Silantyev \cite{ALSS2024} address singularity formation and global
existence for dissipative gCLM equations, with the latter emphasizing the
distinction between the real line and periodic domains.  These works
underscore the need to retain the periodic Biot--Savart law throughout
the reconstruction and dynamical arguments.

\subsection{Strategy of the proof}

The steady-state problem first reduces to the nonlinear eigenvalue equation
\(\Lambda u=\lambda |u|^{1/a-1}u\). Minimization in the odd subspace
removes the constant functions. A comparison of the two reflected
half-period kernels gives positivity, and the Hilbert transform at the
symmetry points gives simple zeros. The scalar power nonlinearity then
determines the regularity at these zeros.

For rescaled profiles, the unknown \(f=-\bar\Omega/x\) records the
shape near the origin. The condition that \(f\) is convex in \(x^2\)
allows a representation as a probability mixture of parabolic caps.
The cap calculation serves two purposes. It proves that the ratio
\(b(f,L)/c(f,L)\) increases strictly with the proposed half-period,
which selects a unique finite period whenever one is needed; it also
proves that the periodic velocity operator preserves the shape
conditions. The periodic input is corrected to
\(x(f(L)-f(x))\), so that its endpoint values agree before a fixed
point has been found. At a fixed point the correction vanishes because
\(f(L)=0\).
The uniqueness in this step concerns the period for a prescribed input
\(f\); the subsequent existence argument does not select a unique
nonlinear fixed point.

The reconstructed shape is obtained by integrating the scalar profile
equation. A common slope bound near the origin gives a decaying upper
barrier. This supplies compactness in the uniform norm and controls
continuity when the selected period tends to infinity. Schauder's
theorem then gives a profile. A commutator identity removes the apparent
division by \(x\) from the differentiated fixed-point equation, making
the interior regularity argument valid through the origin. At a finite
endpoint the velocity has a simple zero, so integrating its logarithmic
singularity gives the exact exponent \(p_L\).

The improved focusing estimate uses a different consequence of the cap
representation. The first-moment interaction has the positive kernel
\(xy/(x+y)\), which factors as an integral of products of Laplace
kernels. Two applications of Cauchy--Schwarz reduce its optimal bound
to one parabolic cap. The profile moment identity turns that bound into
\(a_{\mathrm{foc}}\). Finally, a continuation theorem supplies a connected
set of HQWW fixed points. The mass is continuous on this set by a
positive-kernel identity and Fatou's lemma, even though it is not
continuous on the entire weighted admissible class. The opposite signs
of the spatial rate force a critical profile.

\subsection{Conventions and organization}

The physical variables are \((\omega,u)\), the rescaled variables are
\((\Omega,U)\), and the stationary profiles are \((\bar\Omega,\bar U)\).
The half-period is \(L\), with \(\bar L\) used for a selected stationary
half-period and \(H_\infty=H_{\mathbb R}\).
The stationary rates are \(\bar c_l,\bar c_\omega\), and
\(f=-\bar\Omega/x\) is the normalized scalar profile.
Its focusing decay exponent is \(\rho\), and its finite-endpoint
vanishing exponent is \(p_L\). Auxiliary functions and regularity
exponents introduced within a proof are local to that proof.

The remainder of the paper is organized as follows.
Section~\ref{sec:steady} constructs periodic steady states and determines
their regularity at the symmetry points.
Section~\ref{sec:period_selection} defines the profile map and proves
the estimates that select its period.
Section~\ref{sec:fixed_point} proves invariance, continuity, existence,
and reconstruction of the physical profiles.
Section~\ref{sec:regularity} identifies the whole-line profiles with the
HQWW family and establishes the regularity and endpoint conclusions.
Section~\ref{sec:critical_parameter} proves the sharp first-moment
inequality, the improved focusing range, and existence of a critical
profile. Appendix~\ref{app:hilbert} records the Fourier identities used
in the steady-state reduction.
\section{Periodic steady states}\label{sec:steady}
This section constructs steady states by minimizing a nonlinear Rayleigh quotient in the odd subspace. We first prove existence, then establish positivity on a half-period and determine the regularity at its two endpoints. 
The steady state solution satisfies the nonlinear equation
\begin{equation}\label{eq: steady state}
    a u \om_x = u_x \om,\quad a>0.
\end{equation}
In this section \(H=H_{\mathbb T}\). The Hilbert transform and its
zero-mean velocity potential on the circle of period \(2\pi\) are
\begin{equation}
    u_x(x)=(H \om)(x)= \frac{1}{2\pi}\, \operatorname{p.v.}\int_{-\pi}^{\pi}\om(y)\cot\dr{\frac{x-y}{2}}\idiff y,\quad u(x)=  \frac{1}{\pi}  \int_{-\pi}^{\pi}\om(y)\ln\left|\sin\dr{\frac{x-y}{2}}\right|\idiff y.
\end{equation}
The construction gives the following result.
\begin{theorem}[Existence and regularity of steady states]\label{thm: steady state}
Let \(a>0\). There exist a nonconstant profile
\(u\in H^{1/2}(\mathbb T)\) and a constant \(\lambda>0\) such that
\(\omega:=-\lambda |u|^{1/a-1}u\) satisfies \(u_x=H\omega\)
in distributions. Moreover, \((u,\omega)\) solves \eqref{eq: steady state} globally in the conservative weak sense: for every $\varphi\in C^\infty(\mathbb T)$,
\begin{equation}\label{eq:weak-cons}
    -a\int_{\mathbb T}u\omega\,\varphi_x\,\idiff x
    =(a+1)\int_{\mathbb T}u_x\omega\,\varphi\,\idiff x .
\end{equation}
Equivalently,
\[
    a(u\omega)_x=(a+1)u_x\omega
\]
in the sense of distributions. This formulation is well-defined since \(u\omega\in C^1(\mathbb T)\) and \(u_x\omega\in C(\mathbb T)\). On the open set $\{x\in\mathbb T:u(x)\neq0\}$, the profile is smooth and satisfies \eqref{eq: steady state}
pointwise. The global regularity of $u$ is as follows:
\begin{enumerate}
    \item If $1/a=2N+1$ for some $N\in\mathbb N_0$, then
    \[
        u,\omega\in C^\infty(\mathbb T).
    \]
    \item If $1/a\notin\mathbb N$ and $1/a=m+\alpha$, where $m=\lfloor 1/a\rfloor$ and $0<\alpha<1$, then
    \[
        u\in C^{m+1,\alpha}(\mathbb T),
        \qquad
        \omega\in C^{m,\alpha}(\mathbb T).
    \]
    \item If $1/a=2N$ for some $N\in\mathbb N$, then for every $0<\alpha<1$,
    \[
        u\in C^{2N,\alpha}(\mathbb T),
        \qquad
        \omega\in C^{2N-1,\alpha}(\mathbb T).
    \]
\end{enumerate}
In particular, when $0<a\leq1$, then $(u,\omega)$ is a global classical steady state, namely
\[
    u_x=H\omega,
    \qquad
    a u\omega_x=u_x\omega
\]
hold pointwise on $\mathbb T$.
\end{theorem}

\subsection{Eigenproblem Formulation}
To solve \eqref{eq: steady state}, it suffices to solve for the non-trivial solutions to the following nonlinear eigenproblem on the circle:
\begin{equation}\label{eigen problem_circle}
    \om= -\lambda |u|^{\frac{1}{a}-1}u,\quad H\om = u_x,\quad x\in \mathbb T.
\end{equation}
Indeed, suppose $\om$ and $u$ solve \eqref{eigen problem_circle} (with sufficient regularity). Then, direct computation yields
\[
   \om_x = -\frac{\lambda}{a}u_x|u|^{\frac{1}{a}-1},
\]
and \eqref{eq: steady state} follows.

We parametrize $\mathbb T:= [-\pi,\pi]/\{\pm \pi\}$. The variational construction takes place in the odd subspace
\begin{equation}
    \calX_{\rm odd}:=\{u\in H^{1/2}(\mathbb T):u(-x) = -u(x)\}.
\end{equation}
By Lemma \ref{Appendix: lemma fourier identity for half Laplacian}, since functions in $\calX_{\rm odd}$ have zero means, we have (recall that $\Lambda:= (-\partial^2_x)^{1/2}$)
\[
    \om =-H u_x= -\Lambda u,\quad u\in \calX_{\rm odd}.
\]
Thus, denoting $p:= \frac{1}{a}>0$, the eigenproblem \eqref{eigen problem_circle} reduces to
\begin{equation}\label{eigenproblem_circle_reduced}
    \Lambda u = \lambda |u|^{p-1}u,\quad u\in \calX_{\rm odd}.
\end{equation}
Define
\begin{equation}\label{variational form_epsilon}
    \calQ_{\mathrm{st}}(u):=\frac{\int_{\mathbb T}u\Lambda u\idiff x}{\dr{\int_{\mathbb T}|u|^{p+1}\idiff x}^{2/(p+1)}},\quad u\in \calX_{\rm odd}\backslash\{0\}.
\end{equation}
It is clear that $\calQ_{\mathrm{st}}(u)\geq 0$. Now suppose $\{u_n\}\subset \calX_{\rm odd}$ is a minimizing sequence to $\calQ_{\mathrm{st}}$ with the normalization 
\[
\|u_n\|_{L^{p+1}}=1,\quad \forall\;n.
\]
Since $u_n$ has zero mean, we have
\[
   \|u_n\|_{H^{1/2}(\mathbb T)}\leq C\|\Lambda^{1/2}u_n\|_{L^2(\mathbb T)}
\]
for some universal $C>0$. Thus, $u_n$ is uniformly bounded in $H^{1/2}(\mathbb T)$. Since $H^{1/2}(\mathbb T)\hookrightarrow L^{p+1}(\mathbb T)$ is compact, there exists $u\in H^{1/2}(\mathbb T)$ such that after passing to a subsequence,
\begin{equation}
    u_n\rightharpoonup u \;\text{ in } H^{1/2}(\mathbb T)\quad\text{and}\quad u_n\to u \;\text{ in } L^{p+1}(\mathbb T) . 
\end{equation}
Therefore,
\[
    \|u\|_{L^{p+1}(\mathbb T)} = 1,
\]
so $u\neq 0$, and by weak lower semicontinuity,
\[
    \int_{\mathbb T}u\Lambda u\idiff x
    \leq \liminf_{n\to\infty}\int_{\mathbb T}u_n\Lambda u_n\idiff x.
\]
Hence, $u$ is a minimizer of $\calQ_{\mathrm{st}}$. The Euler–Lagrange equation gives 
\[
    \int_{\mathbb T}u\Lambda\varphi\idiff x
    =\lambda\int_{\mathbb T}|u|^{p-1}u\varphi\idiff x,
    \qquad
    \lambda=\frac{\int_{\mathbb T}u\Lambda u\idiff x}
                   {\int_{\mathbb T}|u|^{p+1}\idiff x}>0,
\]
for all odd test functions $\varphi \in \calX_{\rm odd}$. Since both sides above vanish for even test functions, this actually holds for all $\varphi \in H^{1/2}(\mathbb T)$. Hence
\[
    \Lambda u = \lambda |u|^{p-1}u
\]
in the distributional sense on $\mathbb T$. This proves the following existence statement.
\begin{proposition}\label{prop: existence of weak solution to the nonlinear eigenproblem}There exists a weak solution $u\in \calX_{\rm odd}\backslash\{0\}$ to \eqref{eigenproblem_circle_reduced} with some $\lambda>0$, in the sense that
\[
    \int_{\mathbb T}\varphi\Lambda u\idiff x
    =\lambda\int_{\mathbb T}|u|^{p-1}u\varphi\idiff x,
    \quad \forall\; \varphi\in H^{1/2}(\mathbb T).
\]
    
\end{proposition}
\subsection{Regularity of the nonlinear eigenfunction}
In this subsection, we bootstrap the regularity of the nonlinear eigenfunction via Schauder estimates. The main result is the following:
\begin{proposition}[Regularity classification]\label{prop: regularity classification of steady states}
Let $u\in \calX_{\rm odd}$ be a nontrivial weak solution of \eqref{eigenproblem_circle_reduced}. Then:
\begin{enumerate}[label=\emph{(\roman*)}]
    \item If \(p=2N+1\) is an odd integer,
    then
    \[
        u\in C^\infty(\mathbb T).
    \]

    \item If \(p\notin\mathbb N\), write
    \[
        p=m+\alpha,
        \qquad m=\lfloor p\rfloor,
        \qquad 0<\alpha<1.
    \]
    Then
    \[
        u\in C^{m+1,\alpha}(\mathbb T).
    \]
    Moreover, \(u\) is smooth on every open interval on which \(u\ne0\).

    \item If \(p=2N\) is an even integer,
    then for every \(\alpha\in(0,1)\),
    \[
        u\in C^{2N,\alpha}(\mathbb T).
    \]
    Moreover, \(u\) is smooth on every open interval on which \(u\ne0\).
\end{enumerate}
\end{proposition}

\begin{proof}
Set
\[
    F_p(s):= |s|^{p-1}s,
\]
and \eqref{eigenproblem_circle_reduced} is equivalent to
\begin{equation}\label{eigenproblem_circle_equiv}
    \Lambda u=\lambda F_p(u).
\end{equation}
\medskip
\noindent\textit{Step 1: First bootstrap from \(H^{1/2}\) to almost Lipschitz regularity.}
Since $u\in H^{1/2}(\mathbb T)$,
the one-dimensional Sobolev embedding gives $u\in L^q(\mathbb T)$ for every finite $q<\infty$. Hence $F_p(u)\in L^q(\mathbb T)$
for every finite \(q<\infty\).
Since \(u\) is odd, it has zero mean. Therefore, up to the fixed zero-mean inverse,
\[
    u=\lambda \Lambda^{-1}F_p(u).
\]
For \(1<q<\infty\), the operator $\partial_x\Lambda^{-1}$ has Fourier multiplier
\[
    \frac{ik}{|k|}=i\sgn(k),
\]
which is, up to sign, the periodic Hilbert transform. Hence it is bounded on \(L^q(\mathbb T)\). Therefore $u_x\in L^q(\mathbb T)$ for every finite \(q<\infty\), and thus
\[
    u\in W^{1,q}(\mathbb T)
\]
for every finite \(q<\infty\). By Sobolev embedding in one dimension,
\[
    u\in C^{0,\theta}(\mathbb T)
\]
for every \(0<\theta<1\). In the following we use the standard periodic Schauder estimate for the elliptic Fourier multiplier \(\Lambda\): if \(v\) has zero mean and
\[
    \Lambda v=f,
\]
then for every \(k\in\mathbb N_0\) and every \(0<\alpha<1\),
\[
    f\in C^{k,\alpha}(\mathbb T)
    \quad\Longrightarrow\quad
    v\in C^{k+1,\alpha}(\mathbb T).
\]
Moreover,
\[
    \|v\|_{C^{k+1,\alpha}}
    \lesssim
    \|f\|_{C^{k,\alpha}}.
\]
This is the usual Schauder estimate for an elliptic pseudodifferential operator of order one on the compact manifold \(\mathbb T\).

\medskip
\noindent\textit{Step 2: Odd integer case.}
Assume $p=2N+1$. Then $F_p(u)=u^{2N+1}$ is a smooth function of \(u\). From Step 1,
\[
    u\in C^{0,\theta}(\mathbb T)
\]
for every \(0<\theta<1\). Hence
\[
    F_p(u)=u^{2N+1}\in C^{0,\theta}(\mathbb T)
\]
for every \(0<\theta<1\). By the Schauder estimate,
\[
    u\in C^{1,\theta}(\mathbb T)
\]
for every \(0<\theta<1\). Since the nonlinearity is smooth, the same argument can be iterated. Namely, if $u\in C^{k,\theta}(\mathbb T)$, then
\[
    F_p(u)=u^{2N+1}\in C^{k,\theta}(\mathbb T),
\]
and the equation gives
\[
    u\in C^{k+1,\theta}(\mathbb T).
\]
By induction, we have $u\in C^\infty(\mathbb T)$.

\medskip
\noindent\textit{Step 3: Noninteger case.}
Assume
\[
    p=m+\alpha,
    \qquad m=\lfloor p\rfloor,
    \qquad 0<\alpha<1.
\]
First, we note that
\[
    F_p(s) \in C^{m,\alpha}_{\rm loc}(\R). 
\]
By Step 1, \(u\in C^{0,\theta}\) for every \(0<\theta<1\).

If \(0<p<1\), then \(m=0\) and \(\alpha=p\). Since \(F_p\) is \(p\)-H\"older on bounded sets and \(u\in C^{0,\theta}\), we have
\[
    F_p(u)\in C^{0,p\theta}
\]
for every \(0<\theta<1\). Schauder estimates imply
\[
    u\in C^{1,p\theta}
\]
for every \(0<\theta<1\). In particular, \(u\) is Lipschitz. Since \(F_p\) is \(p\)-H\"older and \(u\) is Lipschitz, it follows that
\[
    F_p(u)\in C^{0,p}.
\]
Applying Schauder once more gives
\[
    u\in C^{1,p}.
\]
This is precisely \(C^{m+1,\alpha}\) in the case \(0<p<1\).

Now assume \(p>1\). Since \(F_p\in C^{m,\alpha}_{\rm loc}\), it is smooth enough to bootstrap up to the threshold \(m+\alpha\). Starting from $u\in C^{0,\theta}$ for all $0<\theta<1$,
and using composition in H\"older spaces, we obtain $F_p(u)\in C^{0,\theta}$
for every \(0<\theta<1\). Schauder estimates give $u\in C^{1,\theta}$
for every \(0<\theta<1\). Iterating this argument until the differentiability threshold of \(F_p\), we obtain
\[
    u\in C^{m+1,\alpha}(\mathbb T).
\]
On any open interval where \(u\ne0\), the scalar map \(F_p\) is smooth as a function of \(u\). Hence the preceding bootstrap can be continued indefinitely on that interval. Therefore \(u\) is smooth away from its zero set.

\medskip
\noindent\textit{Step 4: Even integer case.}
Assume $p=2N$. Then
\[
    F_p(s)=|s|^{2N-1}s\in C^{2N-1,\alpha}_{\rm loc}(\R)
\]
for every \(0<\alpha<1\). The same composition and Schauder argument as in Step 3 applies. This proves the proposition.
\end{proof}
The sharp integer-endpoint conclusion is established in Proposition~\ref{prop:steady_endpoint_behavior}.

\subsection{Proof of Theorem \ref{thm: steady state}}
We now prove Theorem~\ref{thm: steady state}. Take $u\in \calX_{\rm odd}\backslash\{0\}$ as the weak solution of \eqref{eigenproblem_circle_reduced} in Proposition \ref{prop: existence of weak solution to the nonlinear eigenproblem}  and define
\[
    \om:= -\lambda |u|^{1/a-1}u.
\]
Then, $\om = -\Lambda u$ and hence $u_x = H\om$ in the sense of distribution. By Proposition~\ref{prop: regularity classification of steady states}, we know that
\[
   u\om = -\lambda |u|^{1/a+1}\in C^1(\mathbb T). 
\]
Thus, for any $\varphi\in C^\infty(\mathbb T)$, integration by parts yields
\[
\begin{aligned}
 -a\int_{\mathbb T}\om u\varphi_x\idiff x
 &=a\lambda\int_{\mathbb T}|u|^{1/a+1}\varphi_x\idiff x\\
 &=-a\left(\frac1a+1\right)\lambda
      \int_{\mathbb T}u u_x|u|^{1/a-1}\varphi\idiff x\\
 &=(1+a)\int_{\mathbb T}\om u_x\varphi\idiff x.
\end{aligned}
\]
This proves \eqref{eq:weak-cons}. The rest follows from Proposition \ref{prop: regularity classification of steady states}. In particular, when $a\in (0,1]$, both $u$ and $\om$ belong to $C^1(\mathbb T)$, so they are strong solutions to \eqref{eq: steady state}. \qed

\subsection{Positivity of the nonlinear eigenfunction}

In this subsection, we show that the nonlinear eigenfunction obtained from the variational approach has a fixed sign on the positive half period. Since the equation is invariant under the transformation \(u\mapsto -u\), the eigenfunction can be chosen to be strictly positive on \((0,\pi)\).

\begin{proposition}[Positivity of the nonlinear eigenfunction]
\label{prop: positivity of the nonlinear eigenfunction}
Let \(u\in\calX_{\rm odd}\backslash\{0\}\) be a minimizer of \(\calQ_{\mathrm{st}}\) defined in \eqref{variational form_epsilon}. Then \(u\) has a fixed sign on \((0,\pi)\). After replacing \(u\) by \(-u\) if necessary, we have
\[
    u(x)>0,\quad x\in(0,\pi).
\]
\end{proposition}

\begin{proof}
We first prove that \(u\) does not change sign on \((0,\pi)\). Recall the quadratic-form representation of the periodic half-Laplacian:
\begin{equation}\label{eq: quadratic form periodic half laplacian}
    \int_{\mathbb T}u\Lambda u\,\idiff x
    =
    \frac{1}{8\pi}
    \int_{-\pi}^{\pi}\int_{-\pi}^{\pi}
    \frac{|u(x)-u(y)|^2}
    {\sin^2\dr{\frac{x-y}{2}}}
    \idiff y\idiff x.
\end{equation}
Since \(u\) is odd, splitting the integrals into the positive and negative half periods gives
\begin{align}
    \int_{\mathbb T}u\Lambda u\,\idiff x
    =
    \frac{1}{4\pi}
    \int_0^\pi\int_0^\pi
    \bigg(
    &\frac{|u(x)-u(y)|^2}
    {\sin^2\dr{\frac{x-y}{2}}}
    +
    \frac{|u(x)+u(y)|^2}
    {\sin^2\dr{\frac{x+y}{2}}}
    \bigg)
    \idiff y\idiff x.
    \label{eq: odd quadratic form periodic half laplacian}
\end{align}

Define \(v\in\calX_{\rm odd}\) by
\[
    v(x):=
    \begin{cases}
        |u(x)|,&x\in[0,\pi],\\
        -|u(x)|,&x\in[-\pi,0].
    \end{cases}
\]
It is clear that
\begin{equation}\label{eq: same lp norm positivity}
    \int_{\mathbb T}|v|^{p+1}\,\idiff x
    =
    \int_{\mathbb T}|u|^{p+1}\,\idiff x.
\end{equation}
We next compare the quadratic energies of \(u\) and \(v\). For \(x,y\in(0,\pi)\), set
\[
    K_-(x,y):=
    \frac{1}{\sin^2\dr{\frac{x-y}{2}}},
    \qquad
    K_+(x,y):=
    \frac{1}{\sin^2\dr{\frac{x+y}{2}}}.
\]
A direct computation yields
\begin{equation}\label{eq: kernel strict comparison}
\begin{aligned}
    K_-(x,y)-K_+(x,y)
    &=
    \frac{
    \sin^2\dr{\frac{x+y}{2}}
    -
    \sin^2\dr{\frac{x-y}{2}}
    }{
    \sin^2\dr{\frac{x-y}{2}}
    \sin^2\dr{\frac{x+y}{2}}
    }
    =
    \frac{\sin x\sin y}{
    \sin^2\dr{\frac{x-y}{2}}
    \sin^2\dr{\frac{x+y}{2}}
    }>0.
\end{aligned}
\end{equation}
If \(u(x)u(y)\geq0\), then
\[
\begin{aligned}
    &K_-(x,y)|u(x)-u(y)|^2
    +K_+(x,y)|u(x)+u(y)|^2\\
    &\qquad=
    K_-(x,y)|v(x)-v(y)|^2
    +K_+(x,y)|v(x)+v(y)|^2.
\end{aligned}
\]
If \(u(x)u(y)<0\), then
\begin{align}
    &K_-(x,y)|u(x)-u(y)|^2
    +K_+(x,y)|u(x)+u(y)|^2
    \notag\\
    &\quad-
    K_-(x,y)|v(x)-v(y)|^2
    -
    K_+(x,y)|v(x)+v(y)|^2
    \notag\\
    &=
    4|u(x)u(y)|\big(K_-(x,y)-K_+(x,y)\big)>0.
    \label{eq: strict energy improvement}
\end{align}
It follows from \eqref{eq: odd quadratic form periodic half laplacian} that
\[
    \int_{\mathbb T}v\Lambda v\,\idiff x
    \leq
    \int_{\mathbb T}u\Lambda u\,\idiff x.
\]
Together with \eqref{eq: same lp norm positivity}, this gives
\[
    \calQ_{\mathrm{st}}(v)\leq\calQ_{\mathrm{st}}(u).
\]
Since \(u\) is a minimizer, equality must hold. In view of the strict inequality in \eqref{eq: strict energy improvement}, the sets
\[
    \{x\in(0,\pi):u(x)>0\}
    \quad\text{and}\quad
    \{x\in(0,\pi):u(x)<0\}
\]
cannot both have positive measure. Therefore, \(u\) has a fixed sign almost everywhere on \((0,\pi)\). By Proposition \ref{prop: regularity classification of steady states}, \(u\) is continuous. Replacing \(u\) by \(-u\) if necessary, we may assume that
\begin{equation}\label{eq: eigenfunction nonnegative half period}
    u(x)\geq0,\quad x\in(0,\pi).
\end{equation}

It remains to prove that the inequality is strict. Suppose, for the sake of contradiction, that
\[
    u(x_0)=0
\]
for some \(x_0\in(0,\pi)\). Since \(u\geq0\) on \((0,\pi)\), \(x_0\) is a minimum point. The regularity obtained in Proposition \ref{prop: regularity classification of steady states} implies that the pointwise singular-integral formula
\[
    \Lambda u(x)
    =
    \frac{1}{4\pi}\,\operatorname{p.v.}
    \int_{-\pi}^{\pi}
    \frac{u(x)-u(y)}
    {\sin^2\dr{\frac{x-y}{2}}}
    \idiff y
\]
is valid at \(x=x_0\). Using \(u(x_0)=0\) and the oddness of \(u\), we obtain
\begin{align}
    \Lambda u(x_0)
    &=
    -\frac{1}{4\pi}
    \int_0^\pi
    u(y)
    \left(
    \frac{1}{\sin^2\dr{\frac{x_0-y}{2}}}
    -
    \frac{1}{\sin^2\dr{\frac{x_0+y}{2}}}
    \right)
    \idiff y.
    \label{eq: half laplacian at interior zero}
\end{align}
By \eqref{eq: kernel strict comparison}, the expression inside the parentheses is strictly positive for every \(y\in(0,\pi)\). Since \(u\geq0\) on \((0,\pi)\) and \(u\not\equiv0\), it follows that
\[
    \Lambda u(x_0)<0.
\]
On the other hand, the nonlinear eigenvalue equation gives
\[
    \Lambda u(x_0)
    =
    \lambda |u(x_0)|^{p-1}u(x_0)
    =
    0,
\]
which is a contradiction. Therefore,
\[
    u(x)>0,\quad x\in(0,\pi).
\]
This completes the proof.
\end{proof}

\begin{remark}
By the oddness of \(u\), Proposition \ref{prop: positivity of the nonlinear eigenfunction} also implies
\[
    u(x)<0,\quad x\in(-\pi,0).
\]
In particular, the only zeros of \(u\) on one period are the symmetry points
\[
    u(0)=u(\pi)=u(-\pi)=0.
\]
\end{remark}

\paragraph{Integer regularity.}
For an integer \(r\geq1\), the periodic H\"older--Zygmund space
\(C_*^r\) consists of functions whose derivative of order \(r-1\)
belongs to the Zygmund class \(C_*^1\). Here a bounded continuous
function \(v\) belongs to \(C_*^1\) if
\begin{equation}\label{eq:zygmund-definition}
 \sup_{0<|h|\leq1}
 \frac{\|v(\cdot+h)-2v+v(\cdot-h)\|_{L^\infty}}{|h|}<\infty.
\end{equation}
For periodic functions of zero mean, the inverse of \(\Lambda\) maps \(C_*^r\) into
\(C_*^{r+1}\), and \(C^{r-1,1}\subset C_*^r\).
These multiplier bounds, recalled in \cite[Appendix A.1]{Taylor1991}, provide the integer endpoint of the regularity
statements below.

\begin{proposition}[Simple zeros and sharp vorticity cusps]
\label{prop:steady_endpoint_behavior}
Let \(u\) be a minimizer from
Proposition~\ref{prop: positivity of the nonlinear eigenfunction},
chosen positive on \((0,\pi)\), and set
\(p=1/a\) and \(\omega=-\lambda|u|^{p-1}u\).
Its two zeros on the circle are simple. More precisely, the positive
constants
\begin{align}
 s_0:=u'(0)
 &=\frac{\lambda}{\pi}\int_0^\pi
             u(y)^p\cot(y/2)\,\mathrm{d}y>0,
       \label{eq:steady-origin-slope}\\
 s_\pi:=-u'(\pi)
 &=\frac{\lambda}{\pi}\int_0^\pi
             u(y)^p\tan(y/2)\,\mathrm{d}y>0
       \label{eq:steady-antipodal-slope}
\end{align}
satisfy the endpoint expansions
\begin{align}
 \omega(x)
 &=-\lambda s_0^p\operatorname{sgn}(x)|x|^p(1+o(1))
       &&(x\to0),\label{eq:steady-origin-cusp}\\
 \omega(\pi+h)
 &=\lambda s_\pi^p\operatorname{sgn}(h)|h|^p(1+o(1))
       &&(h\to0).\label{eq:steady-antipodal-cusp}
\end{align}
If \(p=m+\alpha\notin\mathbb N\), with
\(m\in\mathbb N_0\) and \(0<\alpha<1\), then the regularity
\(\omega\in C^{m,\alpha}\) is sharp at both zeros:
\(\omega\notin C^{m,\beta}\) for every \(\alpha<\beta\leq1\).
If \(p=2N\) is even, then
\begin{equation}\label{eq:steady-even-endpoint}
 \omega\in C^{2N-1,1}(\mathbb T)\setminus C^{2N}(\mathbb T),
 \qquad u\in C_*^{2N+1}(\mathbb T).
\end{equation}
For odd positive integers \(p\), both functions are smooth.
\end{proposition}

\begin{proof}
The regularity already established in
Proposition~\ref{prop: regularity classification of steady states}
gives \(u\in C^1\), including when \(0<p<1\). Hence
\(u(y)=O(y)\) near zero and \(u(y)=O(\pi-y)\) near \(\pi\).
Since \(p>0\), the integrands in
\eqref{eq:steady-origin-slope} and
\eqref{eq:steady-antipodal-slope} are integrable at their endpoints.
Pairing the positive and negative halves in \(u'=H\omega\) gives
\begin{equation}\label{eq:steady-endpoint-hilbert}
 H\omega(0)=-\frac1\pi\int_0^\pi
                      \omega(y)\cot(y/2)\,\mathrm{d}y,
 \qquad
 H\omega(\pi)=\frac1\pi\int_0^\pi
                      \omega(y)\tan(y/2)\,\mathrm{d}y.
\end{equation}
The strict sign of \(\omega\) on \((0,\pi)\) proves the two slope
formulas. Substitution of
\(u(x)=s_0x+o(x)\) and
\(u(\pi+h)=-s_\pi h+o(h)\) into the definition of \(\omega\)
gives \eqref{eq:steady-origin-cusp} and
\eqref{eq:steady-antipodal-cusp}.

For noninteger \(p=m+\alpha\), either cusp expansion and the known
\(C^{m,\alpha}\) regularity force all derivatives through order
\(m\) to vanish at the corresponding zero. Membership in
\(C^{m,\beta}\), with \(\beta>\alpha\), would then give an
\(O(|x|^{m+\beta})\) Taylor remainder, contradicting the nonzero
coefficient of \(|x|^p\) in \eqref{eq:steady-origin-cusp}.
The same argument applies at \(\pi\).

For \(p=2N\), the scalar function
\(s\mapsto|s|^{2N-1}s\) belongs to \(C^{2N-1,1}\).
The previously proved bound \(u\in C^{2N,\alpha}\) implies that
composition with \(u\) preserves this class. Thus
\(\omega\in C^{2N-1,1}\), and \(-\Lambda u=\omega\), together
with \eqref{eq:zygmund-definition}, gives
\(u\in C_*^{2N+1}\). If \(\omega\) were \(C^{2N}\) at zero,
its Taylor polynomial would give the same limit for
\(\omega(x)/x^{2N}\) on both sides. Equation
\eqref{eq:steady-origin-cusp} gives opposite nonzero limits instead.
The smooth odd-integer case follows from the polynomial nonlinearity.
\end{proof}

\section{The profile map and selection of the period}
\label{sec:period_selection}

The stationary equation leaves both an amplitude and a length normalization
free. After fixing them at the origin, we select the period by a scalar
ratio of two positive functionals. The main issue is to prove that this
ratio is strictly monotone for every admissible shape. A representation
by parabolic caps makes that argument compatible with convexity of the
fixed-point domain.

\subsection{Normalization of the stationary equation}
We first note that the profile equation
\eqref{eqt:general_rescaled_profile_equation_periodic} has a natural scaling
freedom. Suppose that
$(\bar\Omega,\bar U,\bar c_l,\bar c_\omega)$
solves \eqref{eqt:general_rescaled_profile_equation_periodic}, together with
the corresponding Biot--Savart law. Then, for any \(\alpha\neq0\) and
\(\beta>0\), the rescaled tuple $(\bar\Omega_{\alpha,\beta},
\bar U_{\alpha,\beta},
\bar c_{l,\alpha},
\bar c_{\omega,\alpha})$ defined by
\[
\bar\Omega_{\alpha,\beta}(x)
=
\alpha\bar\Omega(\beta x),
\qquad
\bar U_{\alpha,\beta}(x)
=
\frac{\alpha}{\beta}\bar U(\beta x),
\qquad
\bar c_{l,\alpha}
=
\alpha\bar c_l,
\qquad
\bar c_{\omega,\alpha}
=
\alpha\bar c_\omega,
\]
also satisfies the same profile equation, with the spatial period rescaled
accordingly in the periodic case. More precisely, if
$\bar U_x=H_{\bar L}\bar\Omega,$
then
$(\bar U_{\alpha,\beta})_x
=
H_{\bar L/\beta}\bar\Omega_{\alpha,\beta}.$
In the whole-line case, the same transformation preserves the Biot--Savart law
$\bar U_x=H_{\mathbb R}\bar\Omega .$ 
Therefore, we have two degrees of freedom to normalize the profile. In what
follows, we will look for profiles with odd symmetry and fix the normalization
near the origin by requiring
\[
\bar\Omega(x)=-x+x^3+o(x^3),
\qquad
\text{ as }x\to0.
\]
This normalization  rules out the trivial profile and fixes the remaining
scaling ambiguity.
It leaves the amplitude rate \(\bar c_\omega\) determined by the profile:
the compatibility relation \eqref{eqt:compatibility_cl_cw} below gives
\(\bar c_\omega=\bar c_l+(a-1)b\), where \(b=\bar U_x(0)\).
The alternative normalization with amplitude rate \(-1\) is given by
\begin{equation}\label{eq:unit_amplitude_rate_normalization}
 \left(\frac{\bar\Omega}{-\bar c_\omega},
       \frac{\bar U}{-\bar c_\omega},
       \frac{\bar c_l}{-\bar c_\omega},-1\right),
\end{equation}
with the same half-period. Thus the spatial exponent in physical
self-similar time is \(-\bar c_l/\bar c_\omega\).

We now introduce the normalized variables
\begin{equation}\label{eqt:fg_change_variable}
f(x):=-\frac{\bar\Omega(x)}{x},
\qquad
g(x):=
\frac{\bar c_l x+a\bar U(x)}
{x\big(\bar c_l+a\bar U_x(0)\big)} .
\end{equation}
By the oddness of \(\bar\Omega\) and \(\bar U\), both \(f\) and \(g\) are even
functions. Moreover, under the normalization condition,
we have
\[
f(0)=1,
\qquad
f(x)=1-x^2+o(x^2).
\]
We also have \(g(0)=1\), provided
$\bar c_l+a\bar U_x(0)\neq0.$
Taking the limit \(x\to0\) in
\eqref{eqt:general_rescaled_profile_equation_periodic}, using
\(\bar\Omega(x)\sim -x\), gives the compatibility relation
\begin{equation}\label{eqt:compatibility_cl_cw}
\bar c_l+a\bar U_x(0)=\bar c_\omega+\bar U_x(0).
\end{equation}
We can rewrite \eqref{eqt:general_rescaled_profile_equation_periodic} in
terms of \(f\) and \(g\). Since we have
\[
\bar\Omega(x)=-xf(x),
\qquad
\bar c_lx+a\bar U(x)
=
x\big(\bar c_l+a\bar U_x(0)\big)g(x),
\]
the profile equation, together with the compatibility relation
\eqref{eqt:compatibility_cl_cw}, gives
\begin{equation}\label{eqt:f_g_differential_relation}
\frac{f_x}{f}
=
\frac{1}{x}
\left(
\frac{\bar c_\omega+\bar U_x}
{\big(\bar c_l+a\bar U_x(0)\big)g}
-1
\right)
=
\frac{1}{a}\frac{g_x}{g}
+
\frac{1-a}{a}\frac{g-1}{xg}.
\end{equation}
Integrating from \(0\) to \(x\), and using \(f(0)=g(0)=1\), yields
\begin{equation}\label{eqt:f_from_g}
f(x)
=
g(x)^{1/a}
\exp\left(
\frac{1-a}{a}
\int_0^x
\frac{g(y)-1}{y g(y)}
\idiff y
\right).
\end{equation}
We can also express \(g\) in terms of \(f\) through the Biot--Savart law,
\[
\bar U_x=H_{\bar L}\bar\Omega,
\qquad
{\bar L}\in(0,+\infty],
\]
where \(\bar L=+\infty\) corresponds to the whole-line Hilbert transform
\(H_{\mathbb R}\). Since
\(\bar\Omega(x)=-xf(x)\) and \(\bar U(0)=0\), the definition of \(g\) gives
\begin{equation}\label{eqt:g_from_f}
g(x)
=
\frac{
\bar c_lx
+
a\displaystyle\int_0^x
H_{\bar L}(-yf(y))(z)\idiff z
}{
x\left(
\bar c_l
+
aH_{\bar L}(-xf(x))(0)
\right)
}.
\end{equation}
It remains to decide how to choose the parameters \(\bar c_l\) and \(\bar L\). The second normalization condition determines one relation between these parameters. Since
$f(x)=1-x^2+o(x^2)$,
the relation \eqref{eqt:f_from_g} requires
\begin{equation}\label{eqt:g_origin_normalization}
g(x)
=
1-\frac{2a}{3-a}x^2+o(x^2),
\qquad \text{ as } x\to0.
\end{equation}
We now impose this condition on the expression \eqref{eqt:g_from_f}. Since
\[
H_{\bar L}(-xf)(x)
=
H_{\bar L}(-xf)(0)
+
\frac12\partial_x^2\big(H_{\bar L}(-xf)\big)(0)x^2
+
o(x^2),
\]
we obtain from \eqref{eqt:g_from_f} that
\[
g(x)
=
1+
\frac{
a\,\partial_x^2\big(H_{\bar L}(-xf)\big)(0)
}{
6\big(\bar c_l+aH_{\bar L}(-xf)(0)\big)
}
x^2
+
o(x^2).
\]
Comparing this expansion with \eqref{eqt:g_origin_normalization}, we require
the parameters \(\bar c_l\) and \(\bar L\) to satisfy
\begin{equation}\label{eqt:cl_L_constraint}
\frac{
\bar c_l/a+H_{\bar L}(-xf)(0)
}{
-\partial_x^2\big(H_{\bar L}(-xf)\big)(0)/2
}
=
\frac{1-a/3}{2a}.
\end{equation}

The selection first tests the whole-line moments in
\eqref{eqt:cl_L_constraint}. If they are finite and the required rate
\(\bar c_l\) is nonnegative, we keep the whole-line geometry,
\(\bar L=\infty\). Otherwise, we set \(\bar c_l=0\) and choose a finite
period so that the same constraint holds. The transition \(\bar L=\infty,\bar c_l=0\)
is the critical case.

These identities suggest the following reconstruction scheme on the interval
where \(g>0\), with
the velocity in the middle step determined by
\(\bar U_x=H_{\bar L}(-xf)\) and \(\bar U(0)=0\):
\begin{equation}\label{eq:profile_mapping_scheme}
\begin{aligned}
 f&\xrightarrow{\text{period selection}}(\bar L,\bar c_l)\\
  &\xrightarrow{\text{Biot--Savart law}}
 g(x)=\frac{\bar c_lx+a\bar U(x)}
                 {x\bigl(\bar c_l+a\bar U_x(0)\bigr)}\\
  &\xrightarrow{\text{profile equation}}
 \calR(f)(x)=g(x)^{1/a}
 \exp\left(\frac{1-a}{a}\int_0^x
                  \frac{g(y)-1}{yg(y)}\,\idiff y\right).
\end{aligned}
\end{equation}
The three steps use \eqref{eqt:cl_L_constraint},
\eqref{eqt:g_from_f}, and \eqref{eqt:f_from_g}, respectively.
A fixed point makes the input and reconstructed shapes agree.
For continuous inputs, the next subsection defines the moments and
replaces \(-xf(x)\) by the endpoint-corrected function
\(x(f(\bar L)-f(x))\) when \(\bar L<\infty\).
This gives the rigorous map \eqref{eqt:R_map_definition};
the correction vanishes at a fixed point.

\subsection{The admissible shapes and the periodic correction}

We impose a common origin slope and convexity in the squared variable.
The endpoint-corrected periodic input introduced below is continuous
before the fixed-point relation has been imposed.
Consider the Banach space of continuous even functions,
\[
\mathbb V
:=
\left\{
f\in C(\mathbb R):
f(x)=f(-x),\ 
\|f\|_{L^\infty}<+\infty
\right\}.
\]
Fix
\[
\eta:=\frac{1}{3\cdot2^{20}\sqrt{2}}
\]
and define
\begin{equation}\label{eqt:delta_a_definition}
\delta_a
:=
\frac{12(1-a)}{(3-a)(\pi^2+24)}\,\eta.
\end{equation}
In particular, \(0<\delta_a<\eta<1\). Define the decaying upper barrier
\[
f_{a,\delta_a}^{\mathrm{up}}(x):=(1+x^2)^{-\delta_a}.
\]
The choice \eqref{eqt:delta_a_definition} will be recovered from the uniform
tail estimate in the invariant-set proof below. In the
\(L^\infty\)-topology, consider
\[
\begin{aligned}
\mathbb D
:=
\Big\{
f\in\mathbb V:\quad
& f(0)=1,\,
 (1-x^2)_+\leq f(x)\leq f_{a,\delta_a}^{\mathrm{up}}(x),\,
f'_-(1/2)\leq -\eta,\\
& f \text{ is nonincreasing on }[0,+\infty),\,
 f(\sqrt{s}) \text{ is convex in }s\geq0
\Big\}.
\end{aligned}
\]

\begin{lemma}\label{lem:D_properties}
The set \(\mathbb D\) is nonempty, convex,
and compact in \(\mathbb V\) with respect to the \(L^\infty(\mathbb R)\)
norm. Every \(f\in\mathbb D\) satisfies
\[
0\leq f\leq1,
\qquad
\lim_{|x|\to\infty}f(x)=0,
\]
and is globally Lipschitz with Lipschitz constant at most \(2\).
\end{lemma}

\begin{proof}
Let \(f_0(x)=(1-x^2)_+\). Since \(0<\delta_a\leq1\), for
\(0\leq x\leq1\) we have
\[
(1+x^2)^{-\delta_a}
\geq
(1+x^2)^{-1}
\geq
1-x^2.
\]
Thus \(f_0\leq f_{a,\delta_a}^{\mathrm{up}}\). Moreover,
\(s\mapsto(1-s)_+\) is convex and nonincreasing and
\((f_0)'_-(1/2)=-1\leq-\eta\). Hence \(f_0\in\mathbb D\), so
\(\mathbb D\) is nonempty. Convexity follows directly from the defining
conditions.

We next prove closedness. Suppose that \(f_n\in\mathbb D\) and
\(f_n\to f\) uniformly. All the pointwise bounds, evenness, and monotonicity
pass to the limit. The functions \(F_n(s):=f_n(\sqrt{s})\) converge locally
uniformly to \(F(s):=f(\sqrt{s})\), so the secant-slope characterization of
convexity shows that \(F\) is convex. For \(0<h<1/4\), convexity gives
\[
\frac{F_n(1/4)-F_n(1/4-h)}{h}
\leq
(F_n)'_-(1/4)
=
(f_n)'_-(1/2)
\leq-\eta.
\]
Passing to the limit and then taking the supremum over \(h\) yields
\(F'_-(1/4)\leq-\eta\), which is the required slope condition for \(f\).
Thus \(\mathbb D\) is closed.

It remains to prove relative compactness. For \(f\in\mathbb D\), convexity
of \(F(s)=f(\sqrt{s})\) implies, at every differentiability point \(x>0\),
\[
0\leq-\frac{f'(x)}{2x}
\leq
\frac{1-f(x)}{x^2}
\leq
\min\left\{1,\frac{1}{x^2}\right\}.
\]
Consequently \(|f'(x)|\leq2\min\{x,x^{-1}\}\leq2\) almost everywhere,
and the same bound holds for all secant slopes. The family \(\mathbb D\) is
therefore uniformly bounded and equicontinuous. Finally,
\[
\sup_{f\in\mathbb D}\sup_{|x|\geq R}f(x)
\leq
(1+R^2)^{-\delta_a}
\longrightarrow0
\qquad\text{as }R\to\infty.
\]
Arzel\`a--Ascoli on bounded intervals, followed by this uniform tail
estimate, gives relative compactness in the global \(L^\infty\) norm.
Together with closedness, this proves compactness.
\end{proof}

The decay in the upper barrier yields compactness in the uniform norm. Theorem~\ref{thm:R_invariant} proves that the reconstruction preserves this barrier.

For \(f\in\mathbb D\) and \(L\in(0,+\infty)\), set
\[
q_{f,L}(x):=x\big(f(L)-f(x)\big)
\]
on \([-L,L]\), with its odd \(2L\)-periodic extension, and define
\begin{equation}\label{eqt:b_c_definition}
b(f,L):=H_Lq_{f,L}(0).
\end{equation}
For the continuous, possibly nonsmooth profiles in \(\mathbb D\), we take
the following positive integral as the definition of \(c\):
\begin{equation}\label{eqt:c_integral_representation}
c(f,L)
:=
\frac{\pi\big(1-f(L)\big)}{4L}
+
\frac{\pi^2}{4L^3}
\int_0^L
\frac{
y\big(1-f(y)\big)
\cos\left(\frac{\pi y}{2L}\right)
}{
\sin^3\left(\frac{\pi y}{2L}\right)
}
\idiff y.
\end{equation}
For smooth profiles this agrees with
\(-\frac12\partial_x^2H_Lq_{f,L}(0)\); the integral definition avoids
assuming a classical second derivative for every \(f\in\mathbb D\).
The values at \(L=+\infty\) are defined by the limits proved below.
Then given \(L\), the
normalization constraint \eqref{eqt:cl_L_constraint} gives
\begin{equation}\label{eqt:cl_preliminary}
\bar c_l=c_l(f,L):=
\frac{1-a/3}{2}\,c(f,L)-a\,b(f,L).
\end{equation}
We use the abbreviation
\begin{equation}\label{eqt:k_a_definition}
k_a:=\frac{1-a/3}{2a}.
\end{equation}
Since the admissible regimes require
$\bar c_l\geq0,$ and either $\bar c_l=0$ or $\bar L=+\infty$,
we choose \(\bar L\) as a function of \(f\). Define \(\ell(f)\) by
\begin{equation}\label{eqt:ell_a_definition}
\begin{cases}
\ell(f)=+\infty,
& \text{ if } c_l(f,+\infty)\geq0,\\[1mm]
c_l(f,\ell(f))=0,\quad 0<\ell(f)<+\infty,
& \text{ if } c_l(f,+\infty)<0
\end{cases}.
\end{equation}
Existence and uniqueness are
proved in Corollary~\ref{cor:ell_well_defined} below. With this choice,
$c_l(f,\ell(f))=\big(c_l(f,+\infty)\big)_+$.
We now define the fixed-point map. For \(f\in\mathbb D\), define
\begin{equation}\label{eqt:T_map_definition}
\calT(f)(x)
:=
\begin{cases}
0,
&x=0,\\[1mm]
\displaystyle
\frac{1}{x}\int_0^x
H_{\ell(f)}q_{f,\ell(f)}(z)\idiff z
-
b(f,\ell(f)),
&0<|x|\leq \ell(f),\\[3mm]
\calT(f)(\ell(f)),
& |x|>\ell(f).
\end{cases}
\end{equation}
When \(\ell(f)=+\infty\), we use the convention
\(q_{f,+\infty}(y):=-yf(y)\), take the whole-line Hilbert transform, and
omit the last case.
\begin{equation}\label{eqt:G_map_definition}
\begin{aligned}
\calG(f)(x)
&:=
1+
\frac{a\calT(f)(x)}{
c_l(f,\ell(f))+a b(f,\ell(f))}
\\
&=
1+
\frac{2a}{(1-a/3)c(f,\ell(f))}\calT(f)(x).
\end{aligned}
\end{equation}
Let
\[
\tau(f)
:=
\sup\big\{r>0:\ \calG(f)(x)>0\text{ for }0\leq x<r\big\}
\in(0,+\infty].
\]
The nonlinear update is defined only on this positive component and is
continued by zero after its first endpoint:
\begin{equation}\label{eqt:R_map_definition}
\begin{aligned}
\calR(f)(x)
&:=
\calG(f)(x)^{1/a}
\exp\left(
\frac{1-a}{a}
\displaystyle\int_0^{|x|}
\frac{\calG(f)(y)-1}{y\calG(f)(y)}
\idiff y
\right)
\qquad (|x|<\tau(f)).
\end{aligned}
\end{equation}
When \(\tau(f)<+\infty\), we set \(\calR(f)(x)=0\) for
\(|x|\geq\tau(f)\). This convention removes the division by zero in the
finite-period branch.
Corollary~\ref{cor:selected_G_structure} identifies the positive component of \(\calG\), and Theorem~\ref{thm:R_invariant} proves continuity of this extension and invariance of \(\mathbb D\).
We remark that the admissible class \(\mathbb D\), the functionals
\(b,c,c_l,\ell\), and the maps \(\calT,\calG,\calR\) defined above all 
depend on the parameter \(a\in(0,1)\).
For notational simplicity, we suppress this dependence throughout the notation.
For each fixed \(a\in(0,1)\), we seek a fixed point \(f=\calR(f)\).

We first collect some basic properties of the functionals \(b\) and \(c\).
We begin with their integral representations and well-definedness, then derive
uniform bounds for \(c\), and finally establish their continuity with respect
to both the profile \(f\) and the period \(L\).

\begin{lemma}\label{lem:b_c_integral_representation}
Let \(f\in\mathbb D\). For every \(L\in(0,+\infty)\), the quantities
\(b(f,L)\) and \(c(f,L)\) are finite and strictly positive. The functional
\(b\) admits the representation
\begin{equation}\label{eqt:b_integral_representation}
b(f,L)
=
\frac{1}{L}\int_0^L
y\big(f(y)-f(L)\big)
\cot\left(\frac{\pi y}{2L}\right)\idiff y.
\end{equation}
The map \(L\mapsto b(f,L)\) is continuous and strictly increasing on
\((0,+\infty)\), and \(L\mapsto c(f,L)\) is continuous on
\((0,+\infty)\). Moreover, the limits
\[
\begin{aligned}
b(f,+\infty)
&:=
\lim_{L\to+\infty}b(f,L)
\in(0,+\infty],
\\
c(f,+\infty)
&:=
\lim_{L\to+\infty}c(f,L)
\in(0,+\infty)
\end{aligned}
\]
exist and satisfy
\begin{equation}\label{eqt:b_c_whole_line_representation}
b(f,+\infty)
=
\frac{2}{\pi}\int_0^{+\infty}f(y)\idiff y,
\qquad
c(f,+\infty)
=
\frac{2}{\pi}\int_0^{+\infty}
\frac{1-f(y)}{y^2}\idiff y.
\end{equation}
Here the first identity is understood in the extended sense. In particular,
\(b(f,+\infty)<+\infty\) if and only if
\(f\in L^1(0,+\infty)\), while no integrability assumption on \(f\) is
required for \(c(f,+\infty)\).
\end{lemma}

\begin{proof}
By the evenness of \(f\), the function
\(x\big(f(L)-f(x)\big)\) is odd. Therefore,
\[
\begin{aligned}
b(f,L)
&=
\frac{1}{2L}\,\operatorname{p.v.}\int_{-L}^L
y\big(f(L)-f(y)\big)
\cot\left(-\frac{\pi y}{2L}\right)\idiff y
\\
&=
\frac{1}{L}\int_0^L
y\big(f(y)-f(L)\big)
\cot\left(\frac{\pi y}{2L}\right)\idiff y,
\end{aligned}
\]
which proves \eqref{eqt:b_integral_representation}.

To derive the representation of \(c(f,L)\), we decompose
\[
y\big(f(L)-f(y)\big)
=
-\big(1-f(L)\big)y+y\big(1-f(y)\big).
\]
The periodic Hilbert transform of the first term is
\[
H_L\Big(-\big(1-f(L)\big)y\Big)(x)
=
\big(1-f(L)\big)\frac{2L}{\pi}
\log\left(
2\cos\left(\frac{\pi x}{2L}\right)
\right),
\qquad |x|<L.
\]
Consequently, we have
\[
-\frac12\partial_x^2
H_L\Big(-\big(1-f(L)\big)y\Big)(0)
=
\frac{\pi\big(1-f(L)\big)}{4L}.
\]
On the other hand, the lower bound \(f(y)\geq(1-y^2)_+\) gives
\[
y\big(1-f(y)\big)=O(y^3)
\qquad\text{as }y\to0,
\]
For a smooth profile, a symmetric difference quotient, followed by dominated
convergence after the linear term has been separated, yields
\[
-\frac12
\partial_x^2
H_L\big(y(1-f(y))\big)(0)
=
\frac{\pi^2}{4L^3}
\int_0^L
\frac{
y\big(1-f(y)\big)
\cos\left(\frac{\pi y}{2L}\right)
}{
\sin^3\left(\frac{\pi y}{2L}\right)
}
\idiff y.
\]
Combining the two terms proves the asserted smooth-case derivative identity.
For a general \(f\in\mathbb D\), the positive integral in
\eqref{eqt:c_integral_representation} is finite because
\(1-f(y)\leq\min\{y^2,1\}\). Both \(b(f,L)\) and \(c(f,L)\) are strictly
positive: the slope condition implies \(f(y)<1\) for every \(y>0\), and the
integrands above are positive on a set of positive measure.

We also record continuity at finite \(L\). After the change of variables
\(y=Ls\), the correct formulas are
\[
b(f,L)
=
L\int_0^1s\big(f(Ls)-f(L)\big)
\cot\left(\frac{\pi s}{2}\right)\idiff s
\]
and
\[
c(f,L)
=
\frac{\pi(1-f(L))}{4L}
+
\frac{\pi^2}{4L}
\int_0^1
\frac{s(1-f(Ls))\cos(\pi s/2)}{\sin^3(\pi s/2)}
\idiff s.
\]
On every compact \(L\)-interval, the lower cap supplies a uniform integrable
majorant near \(s=0\); away from zero, continuity and boundedness suffice.
Dominated convergence proves continuity of both functionals.

We next consider the limit \(L\to+\infty\). Since the function
\(t\mapsto t\cot t\) is strictly decreasing on \((0,\pi/2)\), for each fixed
\(y>0\), the quantity
\[
K_{b(\cdot,L)}(y)
:=
\mathbf 1_{\{y<L\}}
\frac{y}{L}
\cot\left(\frac{\pi y}{2L}\right)
\]
is nondecreasing in \(L\), with
\[
K_{b(\cdot,L)}(y)
\longrightarrow
\frac{2}{\pi}
\qquad\text{as }L\to+\infty.
\]
Moreover, since \(f\) is nonincreasing and \(f(L)\to0\), for each fixed
\(y>0\),
\[
K_{b(\cdot,L)}(y)
\mathbf 1_{\{y<L\}}\big(f(y)-f(L)\big)
\]
is nondecreasing in \(L\) and converges pointwise to \(2f(y)/\pi\). The
increase is strict after integration because the slope condition makes
\(f(y)-f(L)>0\) on a nonempty interval below every \(L>0\). Hence,
by the monotone convergence theorem, we obtain
\[
\begin{aligned}
b(f,+\infty)
&=
\lim_{L\to+\infty}
\int_0^{+\infty}
K_{b(\cdot,L)}(y)\big(f(y)-f(L)\big)\idiff y=
\frac{2}{\pi}\int_0^{+\infty}f(y)\idiff y,
\end{aligned}
\]
where the integral is understood in the extended sense and may be infinite.

For \(c(f,+\infty)\), consider the kernel
\[
K_{c(\cdot,L)}(y)
:=
\mathbf 1_{\{y<L\}}
\frac{\pi^2}{4L^3}
\frac{
y\cos\left(\frac{\pi y}{2L}\right)
}{
\sin^3\left(\frac{\pi y}{2L}\right)
}.
\]
For each fixed \(y>0\), the convergence
\[
K_{c(\cdot,L)}(y)
\longrightarrow
\frac{2}{\pi y^2}
\qquad\text{as }L\to+\infty
\]
follows directly from the expansion of the trigonometric functions near the
origin. With \(t=\pi y/(2L)\in(0,\pi/2)\), the identity
\[
K_{c(\cdot,L)}(y)
=
\mathbf 1_{\{y<L\}}
\frac{2}{\pi y^2}
\left(\frac{t}{\sin t}\right)^3\cos t
\]
and the inequality
\(\left(\sin t/t\right)^3\geq\cos t\) for
\(0\leq t\leq\pi/2\) imply the bound
\[
0\leq K_{c(\cdot,L)}(y)\leq\frac{2}{\pi y^2}.
\]
For completeness, the trigonometric inequality follows by setting
\(h(t)=3\log(\sin t/t)-\log\cos t\) on \((0,\pi/2)\). Then \(h(0)=0\) by
continuity and
\[
\begin{aligned}
t\sin t\cos t\,h'(t)
&=
t(1+2\cos^2t)-3\sin t\cos t=:N(t),
\\
N'(t)
&=
4\sin t(\sin t-t\cos t)\geq0,
\end{aligned}
\]
because \((\sin t-t\cos t)'=t\sin t\geq0\). Thus \(h\geq0\), and the
endpoint follows by continuity.
Using the bounds defining \(\mathbb D\), we obtain the integrable domination
\[
0
\leq
K_{c(\cdot,L)}(y)\big(1-f(y)\big)
\leq
\frac{2}{\pi}\mathbf 1_{(0,1)}(y)
+
\frac{2}{\pi y^2}\mathbf 1_{[1,+\infty)}(y).
\]
The dominated convergence theorem, together with
\(\pi/(4L)\to0\), therefore gives
\[
\begin{aligned}
c(f,+\infty)
&=
\lim_{L\to+\infty}
\left(
\frac{\pi\big(1-f(L)\big)}{4L}
+
\int_0^{+\infty}
K_{c(\cdot,L)}(y)\big(1-f(y)\big)\idiff y
\right)
\\
&=
\frac{2}{\pi}\int_0^{+\infty}
\frac{1-f(y)}{y^2}\idiff y.
\end{aligned}
\]
The last integral is finite by the preceding bound. This completes the proof.
\end{proof}

\begin{lemma}\label{lem:convex_cap_mixture}
Let \(f\) be even, continuous, nonnegative, nonincreasing on
\([0,\infty)\), and convex in \(x^2\), with \(f(0)=1\) and
\(f(x)\to0\) as \(x\to\infty\). There is a probability measure \(\nu_f\) on
\((0,+\infty)\) such that
\begin{equation}\label{eqt:f_convex_mixture_representation}
f(x)
=
\int_0^{+\infty}
\left(1-\frac{x^2}{r^2}\right)_+
\idiff\nu_f(r).
\end{equation}
If also \(f(x)\geq(1-x^2)_+\), as for every \(f\in\mathbb D\), then
\begin{equation}\label{eqt:cap_measure_negative_moment}
\int_0^{+\infty}\frac{1}{r^2}\idiff\nu_f(r)\leq1.
\end{equation}
\end{lemma}

\begin{proof}
Set \(F(s):=f(\sqrt{s})\). The function \(F\) is nonnegative,
nonincreasing, and convex, with \(F(0)=1\) and \(F(s)\to0\) as
\(s\to+\infty\). Let \(\sigma_f=F''\) be its nonnegative distributional
second derivative. Convexity and the two endpoint limits give
\[
F(s)=\int_s^{+\infty}(t-s)\idiff\sigma_f(t),
\qquad
\int_0^{+\infty}t\idiff\sigma_f(t)=1.
\]
Define \(\nu_f\) as the pushforward of \(t\idiff\sigma_f(t)\) under
\(t\mapsto r=\sqrt t\). Then
\[
F(s)
=
\int_0^{+\infty}\left(1-\frac{s}{r^2}\right)_+
\idiff\nu_f(r),
\]
which is \eqref{eqt:f_convex_mixture_representation}.
Under the additional lower-cap bound, writing \(H(s):=1-F(s)\) gives
\(H(s)\leq s\) for \(0<s\leq1\), while the mixture representation gives
\[
\frac{H(s)}{s}
=
\int_0^{+\infty}
\min\left\{\frac{1}{r^2},\frac{1}{s}\right\}
\idiff\nu_f(r).
\]
Monotone convergence as \(s\to0^+\) proves
\eqref{eqt:cap_measure_negative_moment}.
\end{proof}

\subsection{Cap mixtures and the period ratio}

We next isolate the cap calculation that makes the period selection unique.
Let
\[
\phi(x):=(1-x^2)_+,
\qquad
B(\lambda):=b(\phi,\lambda),
\qquad
C(\lambda):=c(\phi,\lambda),
\]
and define
\begin{equation}\label{eqt:beta_gamma_definition}
\begin{aligned}
\beta_0
&:=
\int_0^1s(1-s^2)\cot\left(\frac{\pi s}{2}\right)\idiff s
=\frac{9\zeta(3)}{\pi^3},\\
\gamma_0
&:=
\frac{\pi}{4}
+
\frac{\pi^2}{4}
\int_0^1
\frac{s^3\cos(\pi s/2)}{\sin^3(\pi s/2)}
\idiff s
=\frac{6\log2}{\pi}.
\end{aligned}
\end{equation}
Here \(\zeta\) denotes the Riemann zeta function. To verify the second
value, integration by parts gives
\begin{equation}\label{eqt:gamma_integral_evaluation}
\begin{aligned}
 \gamma_0
 &=\frac{3\pi}{4}\int_0^1\frac{s^2}{\sin^2(\pi s/2)}\idiff s
 =3\int_0^1s\cot\left(\frac{\pi s}{2}\right)\idiff s\\
 &=-\frac6\pi\int_0^1\log\sin\left(\frac{\pi s}{2}\right)\idiff s
 =\frac{6\log2}{\pi}.
\end{aligned}
\end{equation}
For the first value, integrate by parts and use the Abel-summed Fourier
series \(\log\sin(\pi s/2)=-\log2-\sum_{j\geq1}\cos(j\pi s)/j\).
The Abel sums are bounded by \(C+|\log s|\), so dominated
convergence permits passage to the Abel limit. Since
\(\int_0^1s^2\cos(j\pi s)\idiff s=2(-1)^j/(j\pi)^2\), the resulting
series is absolutely convergent and yields
\[
\begin{aligned}
 \beta_0
 &=-\frac2\pi\int_0^1(1-3s^2)
       \log\sin\left(\frac{\pi s}{2}\right)\idiff s\\
 &=-\frac{12}{\pi^3}\sum_{j=1}^\infty\frac{(-1)^j}{j^3}
 =\frac{9\zeta(3)}{\pi^3}.
\end{aligned}
\]

\begin{lemma}\label{lem:elementary_b_c_ratio}
The functions \(B,C\) are positive, continuous, and locally Lipschitz on
\((0,+\infty)\). For \(0<\lambda\leq1\),
\begin{equation}\label{eqt:B_C_small_scale}
B(\lambda)=\beta_0\lambda^3,
\qquad
C(\lambda)=\gamma_0\lambda.
\end{equation}
For \(\lambda\geq1\), with \(\alpha=\pi/(2\lambda)\),
\begin{equation}\label{eqt:B_C_large_scale}
\begin{aligned}
B(\lambda)
&=
\frac{1}{\lambda}\int_0^1
y(1-y^2)\cot(\alpha y)\idiff y,\\
C(\lambda)
&=
\frac{1}{\lambda}
\left(
3\int_0^1y\cot(\alpha y)\idiff y-\cot\alpha
\right).
\end{aligned}
\end{equation}
Moreover,
\[
B'(\lambda)>0
\quad\text{for a.e. }\lambda>0,
\]
whereas
\[
C'(\lambda)>0\quad(0<\lambda<1),
\qquad
C'(\lambda)<0\quad(\lambda>1).
\]
For every \(\mu\geq1\),
\begin{equation}\label{eqt:B_C_cross_bounds}
B(\mu)\leq\beta_0\mu,
\qquad
C(\mu)\geq\frac{\gamma_0}{\mu}.
\end{equation}
Finally, for \(\lambda,\mu>0\), define
\begin{equation}\label{eqt:E_pairwise_definition}
\begin{aligned}
E(\lambda,\mu):={}&
\mu B'(\lambda)C(\mu)
+\lambda B'(\mu)C(\lambda)\\
&-\frac{\mu^2}{\lambda}B(\lambda)C'(\mu)
-\frac{\lambda^2}{\mu}B(\mu)C'(\lambda).
\end{aligned}
\end{equation}
At every point where the derivatives exist,
\[
E(\lambda,\mu)>0.
\]
The same conclusion holds with one-sided derivatives when
\(\lambda=1\) or \(\mu=1\).
\end{lemma}

\begin{proof}
For \(0<\lambda\leq1\), we have
\(\phi(y)=1-y^2\) on \([0,\lambda]\). Substitution in
\eqref{eqt:b_integral_representation} and
\eqref{eqt:c_integral_representation}, followed by \(y=\lambda s\),
gives \eqref{eqt:B_C_small_scale}.

Suppose that \(\lambda\geq1\). Then \(\phi(\lambda)=0\), and the
representation of \(B\) reduces to the first formula in
\eqref{eqt:B_C_large_scale}. Integrating the formula for \(C\) by parts,
using that \(\phi\) is supported on \([0,1]\), gives the second formula.

The derivative signs follow from the cotangent expansion. The partial-fraction expansion of the cotangent gives
\[
z\cot z
=
1-\sum_{n\geq1}a_nz^{2n},
\qquad
a_n=\frac{2\zeta(2n)}{\pi^{2n}}>0,
\qquad |z|<\pi.
\]
The series and its derivatives converge uniformly for
\(0\leq z\leq\pi/2\). Put \(A=\pi/2\). The first formula in
\eqref{eqt:B_C_large_scale} becomes
\[
B(\lambda)
=
\frac{1}{A}
\left(
\frac23
-
\sum_{n\geq1}
\frac{2a_n\alpha^{2n}}{(2n+1)(2n+3)}
\right).
\]
Since \(\alpha=A/\lambda\), this proves \(B'(\lambda)>0\) for
\(\lambda>1\). If
\[
\Phi(\alpha):=\int_0^\alpha t\cot t\idiff t,
\]
then the second formula in \eqref{eqt:B_C_large_scale} gives
\[
C(\lambda)
=
\frac{1}{A}
\left(
\frac{3\Phi(\alpha)}{\alpha}-\alpha\cot\alpha
\right)
=
\frac{1}{A}
\left(
2+\sum_{n\geq2}
a_n\alpha^{2n}\frac{2n-2}{2n+1}
\right).
\]
This expression is strictly increasing in \(\alpha\), hence strictly
decreasing in \(\lambda>1\). The signs on \((0,1)\) follow from
\eqref{eqt:B_C_small_scale}. These calculations also show continuity and
local Lipschitz regularity at \(\lambda=1\).

To prove \eqref{eqt:B_C_cross_bounds}, apply concavity to
\(H(s):=1-\phi(\sqrt{s})=\min\{s,1\}\). For \(0\leq y\leq\mu\),
\[
H(\mu^2)-H(y^2)
\leq
1-\frac{y^2}{\mu^2},
\qquad
H(y^2)\geq\frac{y^2}{\mu^2}.
\]
Substitution in the positive integral representations of \(B(\mu)\) and
\(C(\mu)\) gives exactly \eqref{eqt:B_C_cross_bounds}.

It remains to prove the pairwise inequality. If
\(0<\lambda,\mu<1\), then \eqref{eqt:B_C_small_scale} gives
\[
E(\lambda,\mu)
=
4\beta_0\gamma_0\lambda^2\mu^2>0.
\]
If \(\lambda,\mu>1\), all four terms in
\eqref{eqt:E_pairwise_definition} are positive because
\(B,C,B'>0\) and \(C'<0\). Finally, suppose that
\(0<\lambda<1<\mu\). Direct substitution yields
\[
E(\lambda,\mu)
=
\lambda^2
\left(
3\beta_0\mu C(\mu)
-\beta_0\mu^2C'(\mu)
+\gamma_0B'(\mu)
-\gamma_0\frac{B(\mu)}{\mu}
\right).
\]
The two derivative terms have positive sign, and
\eqref{eqt:B_C_cross_bounds} gives
\[
3\beta_0\mu C(\mu)
-\gamma_0\frac{B(\mu)}{\mu}
\geq
2\beta_0\gamma_0>0.
\]
The case \(0<\mu<1<\lambda\) follows by symmetry. The one-sided limits in
the small-scale formulas and the large-scale series supply the conclusions at
the junction.
\end{proof}

\begin{lemma}\label{lem:b_over_c_monotonicity}
For every \(f\in\mathbb D\), the function
\[
L\longmapsto
Q_f(L):=\frac{b(f,L)}{c(f,L)}
\]
is continuous and strictly increasing on \((0,+\infty)\). Moreover,
\begin{equation}\label{eqt:b_over_c_small_L}
0<Q_f(L)
\leq
\frac{\beta_0}{\gamma_0}L^2,
\qquad L>0,
\end{equation}
where
\begin{equation}\label{eqt:beta_gamma_ratio}
0<\frac{\beta_0}{\gamma_0}<\frac13.
\end{equation}
In particular, \(Q_f(L)\to0\) as \(L\to0^+\).
\end{lemma}

\begin{proof}
The scaling identities for the caps
\[
\phi_r(x):=\left(1-\frac{x^2}{r^2}\right)_+
\]
are
\[
b(\phi_r,L)=rB(L/r),
\qquad
c(\phi_r,L)=\frac1rC(L/r).
\]
Lemma~\ref{lem:convex_cap_mixture} and linearity therefore give
\begin{equation}\label{eqt:b_c_mixture_representation}
b(f,L)
=
\int_0^{+\infty}rB(L/r)\idiff\nu_f(r),
\qquad
c(f,L)
=
\int_0^{+\infty}\frac1rC(L/r)\idiff\nu_f(r).
\end{equation}

The bounds furnished by \eqref{eqt:B_C_small_scale} and the convergent
series in the proof of Lemma~\ref{lem:elementary_b_c_ratio} justify
differentiation under these integrals on compact \(L\)-intervals. More
explicitly, as \(\lambda\to0\),
\[
B(\lambda)=O(\lambda^3),\quad B'(\lambda)=O(\lambda^2),
\quad C(\lambda)=O(\lambda),\quad C'(\lambda)=O(1),
\]
whereas as \(\lambda\to+\infty\),
\[
B(\lambda)=O(1),\quad B'(\lambda)=O(\lambda^{-3}),
\quad C(\lambda)=O(1),\quad C'(\lambda)=O(\lambda^{-5}).
\]
Together with \eqref{eqt:cap_measure_negative_moment} and the fact that
\(\nu_f\) is a probability measure, these estimates provide integrable
majorants both as \(r\to0\) and as \(r\to+\infty\). Indeed, for
\(L\in[L_0,L_1]\Subset(0,+\infty)\),
\[
|B'(L/r)|
\lesssim_{L_0,L_1}
\begin{cases}
r^3,&0<r\leq L_0,\\
r^{-2},&r\geq L_1,
\end{cases}
\]
and
\[
\left|r^{-2}C'(L/r)\right|
\lesssim_{L_0,L_1}
\begin{cases}
r^3,&0<r\leq L_0,\\
r^{-2},&r\geq L_1.
\end{cases}
\]
On the intervening compact \(r\)-interval, local Lipschitz bounds apply.
These envelopes justify differentiation and Fubini and show that both mixture
functionals are locally absolutely continuous. For a fixed measure
\(\nu_f\), the values of \(L\) at which an atom lies at the junction
\(L/r=1\) form at most a countable set. At every other differentiation
point, symmetrization yields
\begin{equation}\label{eqt:b_c_ratio_derivative}
2\big(\partial_Lb(f,L)c(f,L)-b(f,L)\partial_Lc(f,L)\big)
=
\frac1L
\iint E(L/r,L/s)\idiff\nu_f(r)\idiff\nu_f(s).
\end{equation}
The integrand is strictly positive by
Lemma~\ref{lem:elementary_b_c_ratio}. Hence \(Q_f'(L)>0\) for almost every
\(L\), and local absolute continuity shows that \(Q_f\) is strictly
increasing.

For the quantitative bound, set \(H(s):=1-f(\sqrt{s})\). This function is
nonnegative and concave, with \(H(0)=0\). Thus, for \(0\leq y\leq L\),
\[
H(y^2)\geq\frac{y^2}{L^2}H(L^2),
\qquad
H(L^2)-H(y^2)
\leq
H(L^2)\left(1-\frac{y^2}{L^2}\right).
\]
The integral representations imply
\[
b(f,L)\leq\beta_0L\big(1-f(L)\big),
\qquad
c(f,L)\geq\frac{\gamma_0}{L}\big(1-f(L)\big).
\]
The slope condition ensures \(f(L)<1\) for every \(L>0\), so division gives
\eqref{eqt:b_over_c_small_L}. Finally,
\eqref{eqt:beta_gamma_definition} and
\eqref{eqt:gamma_integral_evaluation} give
\[
 \frac{\gamma_0}{3}-\beta_0
 =\int_0^1s^3\cot\left(\frac{\pi s}{2}\right)\idiff s>0,
\]
which proves \(\beta_0<\gamma_0/3\).
\end{proof}

\begin{corollary}\label{cor:ell_well_defined}
For every \(f\in\mathbb D\), the quantity \(\ell(f)\) in
\eqref{eqt:ell_a_definition} is well-defined without any selection
ambiguity. More precisely, if
\[
c_l(f,+\infty)<0,
\]
there is a unique \(\ell(f)\in(0,+\infty)\) such that
\[
c_l(f,\ell(f))=0.
\]
If \(c_l(f,+\infty)\geq0\), no finite zero exists and
\(\ell(f)=+\infty\). In either case,
\begin{equation}\label{eqt:ell_uniform_lower_bound}
\ell(f)
\geq
\left(
\frac{(1-a/3)/(2a)}{\beta_0/\gamma_0}
\right)^{1/2}
>1.
\end{equation}
\end{corollary}

\begin{proof}
Recall the constant \(k_a\) from \eqref{eqt:k_a_definition}.
Since \(c(f,L)>0\),
\begin{equation}\label{eqt:cl_ratio_factorization}
c_l(f,L)
=
a c(f,L)\big(k_a-Q_f(L)\big).
\end{equation}
By Lemma~\ref{lem:b_over_c_monotonicity}, \(Q_f\) is strictly increasing
from zero to its extended limit
\[
Q_f(+\infty)
=
\frac{b(f,+\infty)}{c(f,+\infty)}
\in(0,+\infty].
\]
If this limit is larger than \(k_a\), there is a unique finite crossing. If
it is at most \(k_a\), strict increase implies \(Q_f(L)<k_a\) for every
finite \(L\), so no finite crossing exists. These are exactly the two signs
of \(c_l(f,+\infty)\) in the statement.

If \(\ell(f)<+\infty\), then
\(Q_f(\ell(f))=k_a\), and
\eqref{eqt:b_over_c_small_L} gives the first inequality in
\eqref{eqt:ell_uniform_lower_bound}. The same inequality is automatic when
\(\ell(f)=+\infty\). Since \(k_a\geq1/3\) for \(a\in(0,1)\), while
\(\beta_0/\gamma_0<1/3\), the lower bound is strictly larger than one.
\end{proof}

\begin{lemma}\label{lem:ell_c_continuity}
The maps
\[
f\longmapsto\ell(f)
\qquad\text{and}\qquad
f\longmapsto c(f,\ell(f))
\]
are continuous on \(\mathbb D\), where convergence to \(+\infty\) has its
usual one-point meaning. More precisely, if \(f_n\to f\) in
\(L^\infty(\mathbb R)\), then
\[
\ell(f_n)\longrightarrow\ell(f)
\quad\text{in }(0,+\infty],
\]
and
\[
c(f_n,\ell(f_n))
\longrightarrow
c(f,\ell(f)).
\]
\end{lemma}

\begin{proof}
Set \(\delta_n:=\|f_n-f\|_{L^\infty}\). For every fixed finite
\(L>0\), the endpoint-corrected representation gives
\[
|b(f_n,L)-b(f,L)|
\leq
2L\delta_n
\int_0^1s\cot\left(\frac{\pi s}{2}\right)\idiff s.
\]
The same estimate is locally uniform for \(L\) in compact subsets of
\((0,+\infty)\).

For \(c\), the lower cap and the upper bound by one imply
\[
|f_n(y)-f(y)|\leq\min\{y^2,\delta_n\}.
\]
The boundary term in \eqref{eqt:c_integral_representation} satisfies
\[
\frac{\pi}{4L}|f_n(L)-f(L)|
\leq
\frac{\pi}{4}\sqrt{\delta_n}.
\]
Using
\[
0\leq
\mathbf1_{\{y<L\}}
\frac{\pi^2}{4L^3}
\frac{y\cos(\pi y/(2L))}{\sin^3(\pi y/(2L))}
\leq
\frac{2}{\pi y^2}
\]
and splitting the integral at \(y=\sqrt{\delta_n}\), we obtain
\begin{equation}\label{eqt:c_uniform_continuity_in_f}
\sup_{L\in(0,+\infty]}
|c(f_n,L)-c(f,L)|
\leq
\left(\frac{\pi}{4}+\frac4\pi\right)\sqrt{\delta_n}.
\end{equation}
It follows that \(Q_{f_n}\to Q_f\) locally uniformly on
\((0,+\infty)\).

Suppose first that \(\ell(f)<+\infty\). For every sufficiently small
\(\varepsilon>0\), strict monotonicity gives
\[
Q_f(\ell(f)-\varepsilon)<k_a
<
Q_f(\ell(f)+\varepsilon).
\]
Local uniform convergence gives the same inequalities for \(Q_{f_n}\) when
\(n\) is large. The characterization in
Corollary~\ref{cor:ell_well_defined} then yields
\[
\ell(f)-\varepsilon<\ell(f_n)<\ell(f)+\varepsilon.
\]
Thus \(\ell(f_n)\to\ell(f)\).

Suppose instead that \(\ell(f)=+\infty\). For every finite \(M>0\),
strict monotonicity and the no-crossing characterization give
\(Q_f(M)<k_a\). Hence \(Q_{f_n}(M)<k_a\) for all sufficiently large \(n\),
and therefore \(\ell(f_n)>M\). Since \(M\) is arbitrary,
\(\ell(f_n)\to+\infty\).

Finally, \(L\mapsto c(f,L)\) is continuous on \((0,+\infty]\) by
Lemma~\ref{lem:b_c_integral_representation}. Therefore,
\[
\begin{aligned}
|c(f_n,\ell(f_n))-c(f,\ell(f))|
\leq{}&
|c(f_n,\ell(f_n))-c(f,\ell(f_n))|\\
&+
|c(f,\ell(f_n))-c(f,\ell(f))|.
\end{aligned}
\]
The first term tends to zero by
\eqref{eqt:c_uniform_continuity_in_f}, and the second tends to zero by
continuity in \(L\).
\end{proof}

\begin{lemma}\label{lem:c_ell_bounds}
For every \(f\in\mathbb D\) and every \(x>0\),
\begin{equation}\label{eqt:c_ell_pointwise_bound}
\begin{aligned}
\frac{2\eta}{\pi}
&\leq
c(f,\ell(f))
\\
&\leq
\frac{\pi}{4\ell(f)}
+
c(f,+\infty)
\\
&\leq
\frac{\pi}{4\ell(f)}
+
\min\left\{
\frac{4(x+1)}{\pi x}\big(1-f(x)\big)^{1/2},
\frac4\pi
\right\}.
\end{aligned}
\end{equation}
Here \(\pi/(4\ell(f))\) is zero when \(\ell(f)=+\infty\).
\end{lemma}

\begin{proof}
Set \(F(s)=f(\sqrt{s})\). Convexity and
\(F'_-(1/4)=f'_-(1/2)\leq-\eta\) imply
\[
1-f(y)
\geq
\eta\min\{y^2,1/4\}
=
\frac{\eta}{4}\left(1-\phi_{1/2}(y)\right).
\]
The integral representation of \(c\), its positivity, and scaling give, for
every \(L>1\),
\[
c(f,L)
\geq
\frac{\eta}{4}c(\phi_{1/2},L)
=
\frac{\eta}{2}C(2L).
\]
By the large-scale series in
Lemma~\ref{lem:elementary_b_c_ratio}, \(C(\lambda)>4/\pi\) for every
finite \(\lambda>1\), and \(C(\lambda)\to4/\pi\) as
\(\lambda\to+\infty\). The lower bound in
\eqref{eqt:c_ell_pointwise_bound} now follows from
Corollary~\ref{cor:ell_well_defined}.

For the upper bounds, the kernel estimate used above and
\(1-f(L)\leq1\) give
\[
c(f,L)
\leq
\frac{\pi}{4L}
+
\frac2\pi\int_0^{+\infty}
\frac{1-f(y)}{y^2}\idiff y
=
\frac{\pi}{4L}+c(f,+\infty).
\]
The same formula holds at \(L=+\infty\) under the stated convention.

It remains to estimate the whole-line term. Since
\(H(s):=1-f(\sqrt{s})\) is nonnegative and concave with \(H(0)=0\),
\(H(s)/s\) is nonincreasing. For each \(x>0\),
\[
1-f(y)
\leq
\begin{cases}
\min\{y^2,1-f(x)\},&0<y\leq x,\\
\min\{(1-f(x))y^2/x^2,1\},&y\geq x.
\end{cases}
\]
Splitting the integral at \(x\) and at the points where the two minima
change gives
\[
\int_0^{+\infty}\frac{1-f(y)}{y^2}\idiff y
\leq
\frac{2(x+1)}{x}\big(1-f(x)\big)^{1/2}.
\]
The simpler bound \(1-f(y)\leq\min\{y^2,1\}\) gives an integral at most
\(2\). Multiplication by \(2/\pi\) proves the remaining estimates.
\end{proof}

For a fixed finite scale \(L\in(0,+\infty)\), write
\begin{equation}\label{eqt:T_fixed_scale_definition}
\calT_L(f)(x)
:=
\frac1x\int_0^xH_Lq_{f,L}(z)\idiff z-b(f,L),
\qquad 0<|x|\leq L,
\end{equation}
with \(\calT_L(f)(0)=0\). For finite \(L\), continue this function by the
constant value \(\calT_L(f)(L)\) on \(|x|>L\). If
\(b(f,+\infty)<+\infty\), the same formula defines
\(\calT_{+\infty}(f)\) on \(\mathbb R\), with
\(q_{f,+\infty}(y)=-yf(y)\), the whole-line Hilbert transform, and no
constant continuation. This finiteness is automatic whenever
\(\ell(f)=+\infty\), by the selection inequality.

\subsection{Shape preservation by the periodic velocity operator}

The next lemma controls both monotonicity and convexity of the averaged velocity. Its endpoint sign will also determine the vanishing order of a finite-period profile.

\begin{lemma}\label{lem:periodic_T_shape}
Let \(f\in\mathbb D\) and \(L\in(0,+\infty)\). The endpoint-corrected
operator \(\calT_L\) in \eqref{eqt:T_fixed_scale_definition} is strictly
decreasing on \((0,L)\), and
\(s\mapsto\calT_L(f)(\sqrt{s})\) is strictly convex on \((0,L^2)\). More
precisely, its one-sided derivatives satisfy
\begin{equation}\label{eqt:T_shape_derivative_signs}
\partial_x^\pm\calT_L(f)(x)<0,
\qquad 0<x<L.
\end{equation}
The distributional second derivative of
\(s\mapsto\calT_L(f)(\sqrt{s})\) is a nonnegative Radon measure with
positive mass on every nonempty subinterval of \((0,L^2)\). Moreover,
\begin{equation}\label{eqt:T_shape_values}
\begin{gathered}
\calT_L(f)(0)=0,\qquad \calT_L(f)(L)=-b(f,L),\\
-b(f,L)<\calT_L(f)(x)<0\quad(0<x<L),
\end{gathered}
\end{equation}
and
\begin{equation}\label{eqt:T_origin_expansion}
\calT_L(f)(x)
=
-\frac{c(f,L)}{3}x^2+o(x^2)
\qquad\text{as }x\to0.
\end{equation}
The left endpoint derivative is
\begin{equation}\label{eqt:T_terminal_derivative}
\partial_x\calT_L(f)(L-)
=
-\frac1{L^2}
\int_0^L
y\big(f(y)-f(L)\big)
\tan\left(\frac{\pi y}{2L}\right)
\idiff y
<0.
\end{equation}
The constant continuation on \([L,+\infty)\) remains convex as a function of
\(x^2\).
\end{lemma}

\begin{proof}
We first prove the result for the quadratic caps
\[
\phi_r(x):=\left(1-\frac{x^2}{r^2}\right)_+.
\]
Scaling gives
\begin{equation}\label{eqt:T_cap_scaling}
\calT_L(\phi_r)(x)
=
r\calT_{L/r}(\phi)(x/r),
\qquad
\phi:=\phi_1.
\end{equation}
It is therefore enough to analyze the dimensionless ratio
\(\lambda=L/r\).

Suppose first that \(0<\lambda\leq1\). Rescaling the period to
\([-1,1]\) makes the corrected odd input a positive multiple of
\[
q_1(y):=y^3-y.
\]
Its sine coefficients and periodic Hilbert transform are
\[
\int_{-1}^1q_1(y)\sin(n\pi y)\idiff y
=
\frac{12(-1)^n}{\pi^3n^3},
\]
and
\[
H_1q_1(x)
=
\frac{12}{\pi^3}
\sum_{n=1}^{+\infty}
\frac{(-1)^{n+1}\cos(n\pi x)}{n^3}.
\]
Put
\[
S_4(x)
:=
\sum_{n=1}^{+\infty}
\frac{(-1)^{n+1}\sin(n\pi x)}{n^4},
\qquad
d_3:=\sum_{n=1}^{+\infty}\frac{(-1)^{n+1}}{n^3}.
\]
Averaging the preceding series gives
\begin{equation}\label{eqt:T_unit_cap_fourier}
\calT_1(\phi)(x)
=
\frac{12}{\pi^4}
\left(\frac{S_4(x)}{x}-\pi d_3\right).
\end{equation}
For \(|x|<1\), the classical Fourier identity
\[
\sum_{n=1}^{+\infty}
\frac{(-1)^{n+1}\cos(n\pi x)}{n}
=
\log\left(2\cos\frac{\pi x}{2}\right)
\]
implies
\[
S_4''''(x)=\frac{\pi^4}{2}\tan\frac{\pi x}{2},
\qquad
S_4'(0)=\pi d_3,
\qquad
S_4'''(0)=-\pi^3\log2.
\]
Every coefficient in the Taylor series of the tangent is positive. Four
integrations therefore yield
\begin{equation}\label{eqt:T_unit_cap_positive_series}
\calT_1(\phi)(\sqrt{s})
=
-\frac{2\log2}{\pi}s
+
\sum_{j=2}^{+\infty}A_js^j,
\qquad
A_j>0,
\qquad 0\leq s<1.
\end{equation}
Thus the function is strictly convex in \(s\). On the other hand,
\(S_4(1)=0\) and
\[
S_4'(1)=-\pi\zeta(3),
\]
so
\[
\partial_x\calT_1(\phi)(1-)
=
-\frac{12\zeta(3)}{\pi^3}<0.
\]
The \(s\)-derivative in
\eqref{eqt:T_unit_cap_positive_series} is strictly increasing and is still
negative at \(s=1\); hence it is negative throughout \((0,1)\). This proves
both signs for every \(0<\lambda\leq1\).

We next consider \(\lambda>1\), keep the cap radius equal to one, and write
\(L=\lambda\). The corrected input is
\[
q(y)=(y^3-y)\mathbf1_{\{|y|<1\}}.
\]
Its whole-line Hilbert transform is
\begin{equation}\label{eqt:unit_cap_whole_line_transform}
F(x)
:=
H_{\mathbb R}q(x)
=
\frac1\pi
\left(
\frac43-2x^2
+
(x^3-x)\log\left|\frac{x+1}{x-1}\right|
\right).
\end{equation}
Polynomial division proves this identity directly. Its interior and exterior
expansions are
\begin{equation}\label{eqt:unit_cap_F_series}
F(x)
=
\frac4{3\pi}-\frac4\pi x^2
+
\frac4\pi
\sum_{j=2}^{+\infty}
\frac{x^{2j}}{(2j-3)(2j-1)},
\qquad |x|<1,
\end{equation}
and
\begin{equation}\label{eqt:unit_cap_F_exterior_series}
F(x)
=
-\frac4\pi
\sum_{n=1}^{+\infty}
\frac{x^{-2n}}{(2n+1)(2n+3)},
\qquad |x|>1.
\end{equation}
In particular, every coefficient in the exterior expansion is negative.

The partial-fraction expansion of the cotangent gives the absolutely and
locally uniformly convergent image representation
\begin{equation}\label{eqt:periodic_cap_image_representation}
R_L(x):=H_Lq(x)
=
\sum_{k\in\mathbb Z}F(x+2kL).
\end{equation}
For \(|x|<1\), write
\[
R_L(x)=F(x)+P_L(x),
\qquad
P_L(x)
:=
\sum_{k=1}^{+\infty}
\big(F(2kL+x)+F(2kL-x)\big).
\]
By \eqref{eqt:unit_cap_F_exterior_series}, the coefficient of \(x^{2j}\)
in \(P_L\) is
\begin{equation}\label{eqt:periodic_image_coefficients}
-\frac8\pi
\sum_{k=1}^{+\infty}
\sum_{n=1}^{+\infty}
\frac{\binom{2n+2j-1}{2j}}
{(2n+1)(2n+3)(2kL)^{2n+2j}}<0.
\end{equation}
For each fixed \(j\), this coefficient is strictly increasing in \(L\).
If
\[
R_L(x)=b_L+\sum_{j=1}^{+\infty}r_j(L)x^{2j},
\qquad |x|<1,
\]
then the case \(L=1\), already proved above, gives
\(r_j(1)>0\) for \(j\geq2\). Hence
\[
r_j(L)>r_j(1)>0\quad(j\geq2),
\qquad
r_1(L)=-c(\phi,L)<0.
\]
Since averaging divides the coefficient of \(x^{2j}\) by \(2j+1\),
\begin{equation}\label{eqt:T_cap_interior_series}
\calT_L(\phi)(\sqrt{s})
=
\sum_{j=1}^{+\infty}\frac{r_j(L)}{2j+1}s^j,
\qquad 0\leq s<1.
\end{equation}
This proves strict \(s\)-convexity inside the cap support and gives a negative
initial derivative.

For \(1<x<L\), every argument in
\eqref{eqt:periodic_cap_image_representation} is outside the cap support.
Set
\[
S_{n,L}(x):=\sum_{k\in\mathbb Z}(x+2kL)^{-2n}>0.
\]
The exterior expansion gives
\[
R_L(x)
=
-\frac4\pi
\sum_{n=1}^{+\infty}
\frac{S_{n,L}(x)}{(2n+1)(2n+3)}.
\]
Pairing the two sides of the period yields
\[
S_{n,L}'(x)
=
-2n
\sum_{k=0}^{+\infty}
\left(
(2kL+x)^{-(2n+1)}
-(2(k+1)L-x)^{-(2n+1)}
\right)<0.
\]
Thus \(R_L'(x)>0\) on \((1,L)\).

Put
\[
T_{\rm cap}(s):=\calT_L(\phi)(\sqrt{s}),
\qquad
R_{\rm cap}(s):=R_L(\sqrt{s}).
\]
Differentiating \(x\calT_L=\int_0^xR_L-b_Lx\) gives
\begin{equation}\label{eqt:R_T_differential_identity}
R_{\rm cap}(s)-b_L=T_{\rm cap}(s)+2sT_{\rm cap}'(s),
\qquad
R_{\rm cap}'(s)=3T_{\rm cap}'(s)+2sT_{\rm cap}''(s).
\end{equation}
The first identity shows that \(T_{\rm cap}'\) is continuous through \(s=1\).
Moreover, the endpoint calculation below gives
\(T_{\rm cap}'(L^2-)<0\). At any hypothetical zero of \(T_{\rm cap}'\) in \((1,L^2)\), the
second identity and \(R_{\rm cap}'>0\) would give \(T_{\rm cap}''>0\). A positive component of
\(T_{\rm cap}'\) could therefore have no right endpoint, contradicting the negative
terminal value. Hence \(T_{\rm cap}'<0\) throughout \((1,L^2)\), and the second
identity then gives \(T_{\rm cap}''>0\). Continuity at \(s=1\) and strict increase of
\(T_{\rm cap}'\) on \((0,1)\) now also give \(T_{\rm cap}'<0\) throughout the cap interior. This
completes the proof of both cap signs.

For later use, we record the endpoint calculation for a general admissible
profile. If \(R_{f,L}:=H_Lq_{f,L}\), then
\[
\calT_L(f)(x)+x\partial_x\calT_L(f)(x)
=
R_{f,L}(x)-b(f,L).
\]
Since \(\calT_L(f)(L)=-b(f,L)\), pairing \(y\) and \(-y\) in the
Hilbert-transform integral at \(x=L\) gives
\[
L\partial_x\calT_L(f)(L-)
=
R_{f,L}(L)
=
-\frac1L
\int_0^L
y\big(f(y)-f(L)\big)
\tan\left(\frac{\pi y}{2L}\right)
\idiff y.
\]
The integral is finite because the common Lipschitz bound gives
\(f(y)-f(L)=O(L-y)\). This proves
\eqref{eqt:T_terminal_derivative} and supplies the terminal sign used above.

Finally, pass from caps to \(f\). By
Lemma~\ref{lem:convex_cap_mixture},
\[
f(x)=\int_0^{+\infty}\phi_r(x)\idiff\nu_f(r).
\]
The positive integral formulas give the mixture identities for \(b\) and
\(c\) by Tonelli's theorem. For \(0<x<L\), the kernel formula also gives
\[
\calT_L(f)(x)+b(f,L)
=
\frac1{\pi x}\int_0^L
y\big(f(y)-f(L)\big)
\log\left|
\frac{\sin\left(\frac{\pi(x+y)}{2L}\right)}
{\sin\left(\frac{\pi(x-y)}{2L}\right)}
\right|\idiff y.
\]
The logarithm is nonnegative because, with
\(A=\pi x/(2L)\) and \(B=\pi y/(2L)\),
\(\sin^2(A+B)-\sin^2(A-B)=\sin(2A)\sin(2B)>0\).
Applying Tonelli to this quantity and then
subtracting the already justified \(b\)-identity gives
\begin{equation}\label{eqt:T_b_c_cap_mixtures}
\begin{aligned}
\calT_L(f)(x)
&=
\int_0^{+\infty}\calT_L(\phi_r)(x)\idiff\nu_f(r),\\
b(f,L)
&=
\int_0^{+\infty}b(\phi_r,L)\idiff\nu_f(r),\\
c(f,L)
&=
\int_0^{+\infty}c(\phi_r,L)\idiff\nu_f(r).
\end{aligned}
\end{equation}
Testing the distributional second \(s\)-derivative against a nonnegative
test function and applying Tonelli transfers strict convexity; integrating
the strict cap inequalities transfers strict decrease.

For each fixed cap, smoothness near the origin and the definition of \(c\)
give
\[
\calT_L(\phi_r)(x)
=
-\frac{c(\phi_r,L)}3x^2+o(x^2).
\]
Cap convexity and decrease imply
\[
-\frac{c(\phi_r,L)}3s
\leq\calT_L(\phi_r)(\sqrt{s})\leq0.
\]
The last line of \eqref{eqt:T_b_c_cap_mixtures} is integrable, so dominated
convergence proves \eqref{eqt:T_origin_expansion} without imposing a
classical second derivative on \(f\). The remaining value and continuation
claims now follow from strict decrease and the endpoint formula.
\end{proof}

\begin{corollary}\label{cor:selected_G_structure}
For every \(f\in\mathbb D\), the selected function \(\calG(f)\) takes values
in \([0,1]\), is nonincreasing on \([0,+\infty)\), is convex as a function
of \(x^2\), and satisfies
\begin{equation}\label{eqt:G_origin_expansion}
\calG(f)(x)
=
1-\frac{2a}{3-a}x^2+o(x^2)
\qquad\text{as }x\to0.
\end{equation}
If \(\ell(f)<+\infty\), then
\begin{equation}\label{eqt:G_finite_positive_set}
\begin{gathered}
\calG(f)(0)=1,\qquad
\calG(f)(x)=0\quad(|x|\geq\ell(f)),\\
0<\calG(f)(x)<1\quad(0<|x|<\ell(f)),
\end{gathered}
\end{equation}
and the terminal zero is simple. If \(\ell(f)=+\infty\), then
\(\calG(f)(x)>0\) for every finite \(x\), and
\[
\lim_{x\to+\infty}\calG(f)(x)
=
1-\frac{b(f,+\infty)}{k_ac(f,+\infty)}
\geq0.
\]
Consequently,
\begin{equation}\label{eqt:tau_selected_scale}
\tau(f)
=
\begin{cases}
\ell(f),&\ell(f)<+\infty,\\
+\infty,&\ell(f)=+\infty.
\end{cases}
\end{equation}
\end{corollary}

\begin{proof}
Suppose first that \(L:=\ell(f)<+\infty\). Scale selection gives
\(b(f,L)=k_ac(f,L)\), and hence
\[
\calG(f)(x)=1+\frac{\calT_L(f)(x)}{b(f,L)}.
\]
All the finite-branch assertions follow immediately from
Lemma~\ref{lem:periodic_T_shape}. In particular,
\[
\partial_x\calG(f)(L-)
=
-\frac1{b(f,L)L^2}
\int_0^L
y\big(f(y)-f(L)\big)
\tan\left(\frac{\pi y}{2L}\right)
\idiff y<0.
\]
The jump from this negative left derivative to the zero right derivative
preserves convexity in \(x^2\).

If \(\ell(f)=+\infty\), then the decay built into \(\mathbb D\) gives
\(f(+\infty)=0\), while selection gives
\(b(f,+\infty)\leq k_ac(f,+\infty)<+\infty\). For each cap, the periodic
image representation \eqref{eqt:periodic_cap_image_representation} converges
locally to its whole-line term as \(L\to+\infty\). Hence
Lemma~\ref{lem:periodic_T_shape} shows that the whole-line cap operator is
nonincreasing and convex in \(x^2\). The positive logarithmic formula above,
with the sines replaced by their arguments, justifies the \(\calT\)-mixture
identity at \(L=+\infty\); the \(b,c\) identities follow directly from their
positive integral formulas. Therefore the same properties pass to
\(\calT_{+\infty}(f)\). Moreover, the cap bound
\[
-\frac{c(\phi_r,+\infty)}3x^2
\leq
\calT_{+\infty}(\phi_r)(x)
\leq0
\]
gives dominated convergence. Thus
\[
\calT_{+\infty}(f)(x)
=
-\frac{c(f,+\infty)}3x^2+o(x^2)
\qquad\text{as }x\to0.
\]
Moreover, for \(x>0\),
\[
\calT_{+\infty}(f)(x)+b(f,+\infty)
=
\frac1{\pi x}
\int_0^{+\infty}
yf(y)\log\left|\frac{x+y}{x-y}\right|
\idiff y
>0.
\]
For every cap, the exterior series
\eqref{eqt:unit_cap_F_exterior_series} gives
\[
\calT_{+\infty}(\phi_r)(x)+b(\phi_r,+\infty)\longrightarrow0
\qquad\text{as }x\to+\infty.
\]
The cap shape and the positive logarithmic formula give the domination
\[
0
\leq
\calT_{+\infty}(\phi_r)(x)+b(\phi_r,+\infty)
\leq
b(\phi_r,+\infty).
\]
Since the last majorant integrates to \(b(f,+\infty)<+\infty\), dominated
convergence in the cap measure proves
\(\calT_{+\infty}(f)(x)\to-b(f,+\infty)\). This proves strict positivity of
\(\calG\) at every finite point and the asserted limit. Finally,
\eqref{eqt:G_origin_expansion} follows from
\eqref{eqt:T_origin_expansion} in the finite branch, from the expansion just
proved in the infinite branch, and from \(1/(3k_a)=2a/(3-a)\).
\end{proof}

\subsection{Continuity and quantitative origin bounds}

The local kernel estimates in this subsection give continuity as the selected period changes. The midpoint estimates then supply the uniform slope and decay needed for the nonlinear reconstruction.

\begin{lemma}\label{lem:G_local_continuity}
If \(f_n\to f\) in \(L^\infty(\mathbb R)\), then, for every \(R>0\),
\[
\|\calG(f_n)-\calG(f)\|_{L^\infty([-R,R])}
\longrightarrow0.
\]
In particular, \(\calG\) is continuous from \(\mathbb D\) into
\(C_{\mathrm{loc}}(\mathbb R)\).
\end{lemma}

\begin{proof}
For finite \(L\), a direct integration of the periodic Hilbert kernel in
\eqref{eqt:T_fixed_scale_definition} gives
\begin{equation}\label{eqt:T_kernel_representation}
\calT_L(f)(x)
=
\int_0^L K_L(x,y)\big(f(y)-f(L)\big)\idiff y,
\end{equation}
where
\begin{equation}\label{eqt:T_kernel_definition}
K_L(x,y)
=
\frac{1}{\pi}
\left(
\frac{y}{x}
\log\left|
\frac{
\sin\left(\frac{\pi(x+y)}{2L}\right)
}{
\sin\left(\frac{\pi(x-y)}{2L}\right)
}
\right|
-
\frac{\pi y}{L}
\cot\left(\frac{\pi y}{2L}\right)
\right).
\end{equation}
We set \(K_L(0,y)=0\), which is the continuous \(x\to0\) limit for every
fixed \(0<y<L\).
At the endpoint the logarithmic term vanishes, and hence
\begin{equation}\label{eqt:T_L_endpoint}
\calT_L(f)(L)=-b(f,L).
\end{equation}
For \(L=+\infty\), the corresponding kernel is
\begin{equation}\label{eqt:T_whole_line_kernel}
K_\infty(x,y)
=
\frac{1}{\pi}
\left(
\frac{y}{x}
\log\left|\frac{x+y}{x-y}\right|-2
\right).
\end{equation}
Here likewise \(K_\infty(0,y)=0\).

For every fixed \(R>0\), there exists \(C_R>0\) such that, whenever
\(L\geq2R\), \(0\leq x\leq R\), and \(0<y<L\),
\[
|K_L(x,y)|
\leq
C_R
\begin{cases}
1+|\log|x-y||,
& 0<y\leq2R,\\[1mm]
R^2y^{-2},
& 2R<y<L.
\end{cases}
\]
Indeed, the first estimate follows from the logarithmic singularity of the
kernel at \(y=x\), while the second follows by expanding the logarithm and
the cotangent around \(x=0\); their constant and linear terms cancel because
of the subtraction of \(b(f,L)\).

The logarithmic singularity is uniformly integrable under translations:
on \(|x-y|<\varepsilon\) its contribution is
\(O(\varepsilon|\log\varepsilon|)\), uniformly for \(x\) in a fixed compact
interval, while on the complement the kernels converge uniformly. Together
with the far-field bound, this proves the following joint continuity. If
\(f_n\to f\) uniformly and \(L_n\to L\), then
\[
\|\calT_{L_n}(f_n)-\calT_L(f)\|_{L^\infty([-R,R])}
\longrightarrow0
\]
whenever \(R<L\). If \(L=+\infty\), this holds for every fixed \(R>0\): in
that case the common upper barrier gives
\(f_n(L_n)\leq f_{a,\delta_a}^{\mathrm{up}}(L_n)\to0\), while the displayed
kernel bound controls the tail.

It remains to justify the moving endpoint when \(L<+\infty\). Let
\(h_{f,L}:=H_Lq_{f,L}\). If \(L\) ranges in a compact subset of
\((0,+\infty)\), the periodic extension of \(q_{f,L}\) has a uniform
Lipschitz constant; one may take \(M_L=1+2L\). Subtracting
\(q_{f,L}(z)\) inside the principal-value integral and writing
\(d_L(y,z)\) for periodic distance gives
\[
\begin{aligned}
|h_{f,L}(z)|
&\leq
\frac{M_L}{2L}
\int_{-L}^L
d_L(y,z)
\left|\cot\left(\frac{\pi(z-y)}{2L}\right)\right|
\idiff y
\leq
\frac{2M_LL}{\pi}.
\end{aligned}
\]
The function \(h_{f,L}\) is even and has zero mean over a period, so
\(\int_0^Lh_{f,L}=0\). Consequently, for \(L/2\leq x\leq L\),
\begin{equation}\label{eqt:T_endpoint_uniform_modulus}
|\calT_L(f)(x)-\calT_L(f)(L)|
=
\frac1x\left|\int_x^Lh_{f,L}(z)\idiff z\right|
\leq
\frac{4M_L}{\pi}(L-x).
\end{equation}
Thus the endpoint modulus is uniform for \(f\in\mathbb D\) and for \(L\) in
a fixed compact interval. Interior convergence, this modulus,
\(b(f_n,L_n)\to b(f,L)\), and the frozen definition together show that the
extensions converge locally uniformly across \(L_n\to L\). Therefore, by
Lemma~\ref{lem:ell_c_continuity},
\[
\|\calT(f_n)-\calT(f)\|_{L^\infty([-R,R])}
\longrightarrow0.
\]

Moreover, Lemma~\ref{lem:c_ell_bounds} gives
\(c(f,\ell(f))\geq2\eta/\pi\). Since
\[
c(f_n,\ell(f_n))
\longrightarrow c(f,\ell(f)),
\]
we conclude from
\[
\calG(f)
=
1+
\frac{\calT(f)}{k_a c(f,\ell(f))}
\]
that
\[
\calG(f_n)\longrightarrow\calG(f)
\]
uniformly on every compact interval.
\end{proof}

\begin{lemma}\label{lem:periodic_T_midpoint_bounds}
Let \(f\in\mathbb D\), put \(L:=\ell(f)>1\), and write
\[
T:=\calT_L(f),
\qquad
c_L:=c(f,L).
\]
For \(0<x\leq1\), the one-sided derivatives satisfy
\begin{align}
T'_-(x)
&\leq
\frac{x f'_-(x)}{2\pi},
\label{eqt:T_derivative_slope_bound}\\
T'_-(x)
&\leq
-\frac{\big(1-f(x)\big)^{5/2}}
{3\sqrt2\,\pi x^3}.
\label{eqt:T_derivative_height_bound}
\end{align}
In particular,
\begin{equation}\label{eqt:T_midpoint_interpolated_bound}
-T'_-\left(\frac12\right)
\geq
\frac{
\left|f'_-\left(\frac12\right)\right|^{4/5}
\big(1-f(1/2)\big)^{1/2}
}{
\pi(96\sqrt2)^{1/5}
}.
\end{equation}
At the selected scale,
\begin{equation}\label{eqt:selected_c_midpoint_upper_bound}
c_L
\leq
\frac{\pi^2+24}{2\pi}
\big(1-f(1/2)\big)^{1/2}.
\end{equation}
\end{lemma}

\begin{proof}
We first derive a differentiated-kernel comparison. For finite \(L\), let
\[
A_L(x,y):=\int_0^yK_L(x,t)\idiff t,
\qquad
W_L(x,y):=-\partial_xA_L(x,y),
\]
where \(K_L\) is defined in \eqref{eqt:T_kernel_definition}. Integration by
parts in \eqref{eqt:T_kernel_representation} gives
\begin{equation}\label{eqt:T_derivative_W_representation}
T'(x)=\int_0^LW_L(x,y)f'(y)\idiff y
\end{equation}
at differentiability points. Both boundary terms vanish:
\(A_L(x,0)=0\), while \(f(L)-f(L)=0\).
The differentiation in \(x\) is justified by first splitting the
\(y\)-integral into \(\lvert y-x\rvert<\varepsilon\) and its complement.
The difference \(W_L-W_\infty\) is analytic on the domain used below, while
\(W_\infty\) has only an integrable logarithmic singularity at \(y=x\).
Since every \(f\in\mathbb D\) is Lipschitz, the first part is
\(O(\varepsilon\lvert\log\varepsilon\rvert)\), locally uniformly in \(x\),
and ordinary dominated convergence applies to the complement. The same
argument on increasing finite intervals gives the whole-line identity; its
boundary term vanishes because \(f(y)\to0\) and \(A_\infty(x,y)\) has a
finite limit as \(y\to+\infty\).

Let \(K_\infty\) be the kernel in
\eqref{eqt:T_whole_line_kernel}, and define \(W_\infty\) in the same way.
Put \(\alpha=\pi/(2L)\) and
\[
p(z):=\log\left(\frac{\sin z}{z}\right)
=
-\sum_{m=1}^{+\infty}
\frac{\zeta(2m)}{m\pi^{2m}}z^{2m},
\qquad |z|<\pi.
\]
Up to a term independent of \(x\), the difference \(K_L-K_\infty\) is
\[
\frac{y}{\pi x}
\big(p(\alpha(x+y))-p(\alpha(x-y))\big).
\]
For \(0<x\leq1\), \(0<y<L\), and \(L>1\), one has
\(x+y<L+1<2L\), so the series may be differentiated absolutely. Its
\(x\)-dependent factor expands as
\[
\begin{aligned}
&\frac{y}{x}
\big(p(\alpha(x+y))-p(\alpha(x-y))\big)\\
&\quad=
-2\sum_{m=1}^{+\infty}
\frac{\zeta(2m)}{m(2L)^{2m}}
\sum_{j=0}^{m-1}
\binom{2m}{2j+1}
x^{2(m-j-1)}y^{2j+2}.
\end{aligned}
\]
This is a negative series in nonnegative even powers of \(x\), and hence
\(\partial_x(K_L-K_\infty)\leq0\). It follows that
\begin{equation}\label{eqt:W_periodic_whole_line_comparison}
W_L(x,y)\geq W_\infty(x,y),
\qquad
0<x\leq1,\quad 0<y<L.
\end{equation}
The pointwise comparison is understood off the diagonal \(y=x\); both
kernels have an integrable logarithmic singularity there, so it holds almost
everywhere in every integral below.

For completeness, define
\[
\mathscr F(t)
:=
\frac{t^2-1}{2t}
\log\left|\frac{t+1}{t-1}\right|+1.
\]
A direct integration of \(K_\infty\) gives
\[
W_\infty(x,y)=\frac1\pi\mathscr F'(x/y).
\]
For \(0<u<1\),
\[
\begin{aligned}
\mathscr F'(u)
&=
\sum_{n=1}^{+\infty}
\frac{4n}{4n^2-1}u^{2n-1}>0,\\
\mathscr F'(1/u)
&=
\sum_{n=1}^{+\infty}
\frac{4n}{4n^2-1}u^{2n+1}>0.
\end{aligned}
\]
Thus \(W_\infty>0\), and \eqref{eqt:W_periodic_whole_line_comparison}
implies \(W_L>0\). For \(0<y<x\), the second series also yields
\begin{equation}\label{eqt:W_whole_line_local_bounds}
W_\infty(x,y)
\geq
\frac4{3\pi}\left(\frac yx\right)^3,
\qquad
\int_0^x\frac yxW_\infty(x,y)\idiff y
=
\frac{x}{2\pi}.
\end{equation}
The last identity follows by termwise integration and the telescoping sum
\[
\sum_{n=1}^{+\infty}
\frac{4n}{(4n^2-1)(2n+3)}
=
\frac12.
\]

The convexity of \(s\mapsto f(\sqrt{s})\) makes \(f'(y)/y\)
nondecreasing. Hence, for \(0<y<x\),
\[
\frac{f'(y)}{y}\leq\frac{f'_-(x)}{x}.
\]
Using \eqref{eqt:T_derivative_W_representation}, positivity, and
\eqref{eqt:W_whole_line_local_bounds}, we obtain
\[
\begin{aligned}
T'(x)
&\leq
\int_0^xW_L(x,y)f'(y)\idiff y\\
&\leq
\frac{f'_-(x)}x
\int_0^xyW_\infty(x,y)\idiff y
=
\frac{x f'_-(x)}{2\pi}.
\end{aligned}
\]
For the second estimate, the case \(f(x)=1\) is immediate.
Otherwise, the lower cap gives \(1-f(x)\leq x^2\), so
\(z:=\sqrt{(1-f(x))/2}<x\). Again using positivity and
\eqref{eqt:W_whole_line_local_bounds},
\[
\begin{aligned}
T'(x)
&\leq\frac4{3\pi x^3}\int_z^xy^3f'(y)\idiff y
 \leq\frac{4z^3}{3\pi x^3}\big(f(x)-f(z)\big)\\
&\leq-\frac{\bigl(1-f(x)\bigr)^{5/2}}{3\sqrt2\,\pi x^3}.
\end{aligned}
\]
The calculation applies directly at differentiability points. Taking
differentiability points increasingly to \(x\), the one-sided derivatives of
the two convex functions \(s\mapsto f(\sqrt{s})\) and
\(s\mapsto T(\sqrt{s})\) converge to their left derivatives. This proves
\eqref{eqt:T_derivative_slope_bound} and
\eqref{eqt:T_derivative_height_bound}. The whole-line case follows from the
same calculation with \(W_L=W_\infty\).

At \(x=1/2\), the two estimates say
\[
-T'_-(1/2)
\geq
\frac{|f'_-(1/2)|}{4\pi},
\qquad
-T'_-(1/2)
\geq
\frac8{3\sqrt2\pi}\big(1-f(1/2)\big)^{5/2}.
\]
Their weighted geometric mean, with weights \(4/5\) and \(1/5\), is
\eqref{eqt:T_midpoint_interpolated_bound}.

It remains to prove \eqref{eqt:selected_c_midpoint_upper_bound}. If
\(H(s):=1-f(\sqrt{s})\), then \(H\) is nonnegative and concave, and
\(H(s)/s\) is nonincreasing. For \(0<x\leq L\),
\[
\frac{1-f(L)}L
\leq
\min\left\{
\frac{L}{x^2}\big(1-f(x)\big),\frac1L
\right\}
\leq
\frac{\big(1-f(x)\big)^{1/2}}x.
\]
Combining this with the proof of Lemma~\ref{lem:c_ell_bounds} gives
\[
c(f,L)
\leq
\left(
\frac{\pi}{4x}
+
\frac{4(x+1)}{\pi x}
\right)
\big(1-f(x)\big)^{1/2}.
\]
The first term is absent at \(L=+\infty\). Setting \(x=1/2\) proves
\eqref{eqt:selected_c_midpoint_upper_bound}.
\end{proof}

\section{Existence and reconstruction of the profiles}
\label{sec:fixed_point}

We now apply the estimates of Section~\ref{sec:period_selection} to the
nonlinear reconstruction. The proof has three parts: preservation of the
shape class and its decay, continuity through a changing terminal point,
and recovery of the physical stationary equation from a fixed point.

\subsection{Invariance of the admissible set}
\begin{theorem}\label{thm:R_invariant}
For the exponent \(\delta_a\) in \eqref{eqt:delta_a_definition},
\(\calR(\mathbb D)\subset\mathbb D\).
\end{theorem}

\begin{proof}
Fix \(f\in\mathbb D\), and abbreviate
\[
g:=\calG(f),
\qquad
h:=\calR(f),
\qquad
L:=\ell(f),
\qquad
c_L:=c(f,L).
\]
Corollary~\ref{cor:selected_G_structure} gives
\[
g(0)=1,\qquad 0\leq g\leq1,
\]
and shows that \(g\) is nonincreasing and convex as a function of \(x^2\).
It also identifies its positive component as \([0,L)\) when \(L<+\infty\)
and as \([0,+\infty)\) otherwise.

On this positive component, logarithmic differentiation gives
\begin{equation}\label{eqt:R_logarithmic_derivative}
\frac{h'(x)}{h(x)}
=
\frac1a\frac{g'(x)}{g(x)}
+
\frac{1-a}{a}\frac{g(x)-1}{xg(x)}.
\end{equation}
Both terms are nonpositive. Thus \(0<h\leq g^{1/a}\leq1\) and \(h\) is
nonincreasing. At a finite terminal point,
\(h(x)\leq g(x)^{1/a}\to0\), so the zero continuation is continuous.

We next prove convexity in \(x^2\) without a smoothness assumption. On the
positive component put
\[
G(s):=g(\sqrt{s}),
\qquad
H(s):=h(\sqrt{s}).
\]
After changing variables in the defining integral for \(h\),
\[
H(s)
=
G(s)^{1/a}
\exp\left(
\frac{1-a}{2a}
\int_0^s\frac{G(t)-1}{tG(t)}\idiff t
\right).
\]
At every interior point its right derivative is
\[
H'_+(s)=A_+(s)B(s),
\]
where
\[
\begin{aligned}
A_+(s)
&:=
\frac1aG'_+(s)
+
\frac{1-a}{2a}\frac{G(s)-1}{s},\\
B(s)
&:=
\frac{H(s)}{G(s)}.
\end{aligned}
\]
The convexity of \(G\) makes both \(G'_+\) and the secant slope
\((G(s)-G(0))/s\) nondecreasing. Hence \(A_+\) is nondecreasing, and
\(A_+\leq0\) because \(G\) is nonincreasing. Moreover,
\[
B(s)
=
G(s)^{(1-a)/a}
\exp\left(
\frac{1-a}{2a}
\int_0^s\frac{G(t)-1}{tG(t)}\idiff t
\right)
\]
is positive and nonincreasing. Thus, for \(s_1<s_2\),
\[
\begin{aligned}
&A_+(s_2)B(s_2)-A_+(s_1)B(s_1)\\
&\quad=
\big(A_+(s_2)-A_+(s_1)\big)B(s_2)\\
&\qquad+
A_+(s_1)\big(B(s_2)-B(s_1)\big)
\geq0.
\end{aligned}
\]
The right derivative of \(H\) is therefore nondecreasing, which proves
convexity. At a finite terminal point, the left derivative is nonpositive
and the zero continuation has right derivative zero, so the upward derivative
jump preserves convexity.

The origin expansion \eqref{eqt:G_origin_expansion} gives
\[
\log h(x)
=
-\frac1a\frac{2a}{3-a}x^2
-\frac{1-a}{2a}\frac{2a}{3-a}x^2
+o(x^2)
=
-x^2+o(x^2).
\]
Hence \(h(x)=1-x^2+o(x^2)\). The supporting line at the origin of the convex
function \(s\mapsto h(\sqrt{s})\) now gives
\[
h(x)\geq(1-x^2)_+.
\]

It remains to verify the midpoint slope and the upper barrier. Convexity of
\(s\mapsto g(\sqrt{s})\) implies
\[
\frac{g(x)-1}{x}\leq\frac{g'_-(x)}2.
\]
Since \(g'=T'/(k_ac_L)\), where \(T=\calT_L(f)\), equation
\eqref{eqt:R_logarithmic_derivative} yields
\[
h'_-(x)
\leq
\frac{3T'_-(x)}{c_L}\frac{h(x)}{g(x)}.
\]
At \(x=1/2\), the lower cap for \(h\) and \(g\leq1\) give
\[
\frac{h(1/2)}{g(1/2)}\geq\frac34.
\]
Using Lemma~\ref{lem:periodic_T_midpoint_bounds},
\[
-h'_-(1/2)
\geq
K_*\left|f'_-(1/2)\right|^{4/5},
\qquad
K_*:=
\frac9{2(\pi^2+24)(96\sqrt2)^{1/5}}.
\]
The choice of \(\eta\) satisfies
\[
K_*\eta^{-1/5}
=
\frac{36}{\pi^2+24}>1.
\]
Since \(|f'_-(1/2)|\geq\eta\), we conclude that
\[
h'_-(1/2)<-\eta.
\]

The same estimates give a uniform drop of \(g\). Convexity of
\(s\mapsto T(\sqrt{s})\), together with \(T(0)=0\), gives
\[
T(1/2)\leq\frac14T'_-(1/2).
\]
Equations \eqref{eqt:T_midpoint_interpolated_bound} and
\eqref{eqt:selected_c_midpoint_upper_bound}, the input slope condition, and
the exact value of \(\eta\) yield
\[
-\frac{T'_-(1/2)}{c_L}
\geq
\frac{16}{\pi^2+24}\eta.
\]
Consequently,
\[
1-g(1/2)
=
\frac{-T(1/2)}{k_ac_L}
\geq
\frac{24a}{(3-a)(\pi^2+24)}\eta.
\]
Define
\[
p_a^{\mathrm{tail}}
:=
\frac{24(1-a)}{(3-a)(\pi^2+24)}\eta
=
2\delta_a.
\]
Since \(g\) is nonincreasing, for \(x\geq1/2\) before a possible terminal
point,
\[
\frac{1-a}{a}\left(\frac1{g(x)}-1\right)
\geq
\frac{1-a}{a}\big(1-g(1/2)\big)
\geq
p_a^{\mathrm{tail}}.
\]
Dropping the nonpositive \(g'/g\) term in
\eqref{eqt:R_logarithmic_derivative} and integrating gives
\begin{equation}\label{eqt:R_uniform_tail_bound}
h(x)
\leq
h(1/2)(2x)^{-p_a^{\mathrm{tail}}},
\qquad x\geq\frac12.
\end{equation}
The estimate remains true after a finite terminal point.

Convexity and the midpoint slope of \(h\) give
\[
h(x)\leq1-\eta x^2
\quad(0\leq x\leq1/2),
\qquad
h(1/2)\leq1-\frac{\eta}{4}.
\]
Since \(\delta_a\leq\eta\) and
\(s\mapsto(1+s)^{-\delta_a}\) is convex,
\[
h(x)
\leq
1-\eta x^2
\leq
(1+x^2)^{-\delta_a},
\qquad 0\leq x\leq\frac12.
\]
For \(x\geq1/2\), put \(t=2x\geq1\). Then
\[
1+x^2\leq\frac54t^2,
\qquad
\left(\frac54\right)^{-\delta_a}
\geq1-\frac{\delta_a}{4}
\geq1-\frac{\eta}{4}.
\]
Using \eqref{eqt:R_uniform_tail_bound} and
\(p_a^{\mathrm{tail}}=2\delta_a\), we obtain
\[
h(x)
\leq
\left(1-\frac{\eta}{4}\right)t^{-2\delta_a}
\leq
\left(\frac54\right)^{-\delta_a}t^{-2\delta_a}
\leq
(1+x^2)^{-\delta_a}.
\]
Thus \(h\) is even, continuous, normalized, bounded between the two defining
barriers, nonincreasing, convex in \(x^2\), and satisfies the required
midpoint slope. Hence \(h\in\mathbb D\).
\end{proof}

\subsection{Continuity of the nonlinear map}

\begin{proposition}\label{prop:R_continuity}
The map
\[
\calR:\mathbb D\longrightarrow\mathbb D
\]
is continuous with respect to the \(L^\infty(\mathbb R)\)-topology.
More precisely, if \(f_n,f\in\mathbb D\) and
\[
\|f_n-f\|_{L^\infty}\longrightarrow0,
\]
then
\[
\|\calR(f_n)-\calR(f)\|_{L^\infty}
\longrightarrow0.
\]
\end{proposition}

\begin{proof}
Let
\[
g_n:=\calG(f_n),
\qquad
g:=\calG(f),
\qquad
h_n:=\calR(f_n),
\qquad
h:=\calR(f).
\]
By the continuity properties of \(c,\ell,\calT\), and \(\calG\) established
above, we have
\begin{equation}\label{eqt:G_local_convergence_for_R}
g_n\longrightarrow g
\qquad\text{locally uniformly on }\mathbb R.
\end{equation}
Moreover, the structural properties of \(\calG\) give
\begin{equation}\label{eqt:G_basic_properties_for_R}
0\leq g_n(x),g(x)\leq1,
\qquad
g_n,g\text{ are nonincreasing on }[0,+\infty).
\end{equation}

We first prove that \(h_n\to h\) locally uniformly. Since all the functions
under consideration are even, it is enough to work on \([0,+\infty)\).

By the normalization at the origin,
\[
g_n(x)=1-\frac{2a}{3-a}x^2+o(x^2),
\qquad
g(x)=1-\frac{2a}{3-a}x^2+o(x^2).
\]
Since \(g_n(\sqrt{s})\) and \(g(\sqrt{s})\) are convex in \(s\), their
supporting lines at \(s=0\) give the uniform estimate
\begin{equation}\label{eqt:G_uniform_origin_bound}
1-\frac{2a}{3-a}x^2
\leq
g_n(x),g(x)
\leq1,
\qquad x\geq0.
\end{equation}
In particular, if
\[
r_a:=\sqrt{\frac{3-a}{4a}},
\]
then
\[
g_n(x),g(x)\geq\frac12,
\qquad 0\leq x\leq r_a.
\]

For every \(x\) such that \(g_n>0\) on \([0,x]\), define
\[
\Psi_n(x)
:=
\frac{1-a}{a}
\int_0^x
\frac{g_n(y)-1}{y g_n(y)}
\idiff y,
\]
and define \(\Psi\) analogously. Since \(0<g_n,g\leq1\), we have
\begin{equation}\label{eqt:Psi_nonpositive}
\Psi_n(x)\leq0,
\qquad
\Psi(x)\leq0.
\end{equation}
In particular,
\begin{equation}\label{eqt:R_bounded_by_G}
0\leq h_n(x)\leq g_n(x)^{1/a},
\qquad
0\leq h(x)\leq g(x)^{1/a}.
\end{equation}

We next establish the continuity of the exponential factors away from the
first zero of \(g\). Let \(r>0\) be such that
\[
g(x)>0,
\qquad 0\leq x\leq r.
\]
By \eqref{eqt:G_local_convergence_for_R}, there exists \(m_r>0\) such that,
for all sufficiently large \(n\),
\[
g_n(x),g(x)\geq m_r,
\qquad 0\leq x\leq r.
\]
For \(0<\delta<\min\{r,r_a\}\), we split
\[
\begin{aligned}
\sup_{0\leq x\leq r}|\Psi_n(x)-\Psi(x)|
&\leq
\frac{1-a}{a}
\int_0^\delta
\left|
\frac{g_n(y)-1}{y g_n(y)}
-
\frac{g(y)-1}{y g(y)}
\right|
\idiff y\\
&\quad+
\frac{1-a}{a}
\int_\delta^r
\left|
\frac{g_n(y)-1}{y g_n(y)}
-
\frac{g(y)-1}{y g(y)}
\right|
\idiff y.
\end{aligned}
\]
By \eqref{eqt:G_uniform_origin_bound}, for
\(0<y\leq\delta\leq r_a\),
\[
\left|
\frac{g_n(y)-1}{y g_n(y)}
\right|
+
\left|
\frac{g(y)-1}{y g(y)}
\right|
\leq
\frac{8a}{3-a}y.
\]
Therefore,
\begin{equation}\label{eqt:Psi_near_zero_bound}
\int_0^\delta
\left|
\frac{g_n(y)-1}{y g_n(y)}
-
\frac{g(y)-1}{y g(y)}
\right|
\idiff y
\leq
\frac{4a}{3-a}\delta^2.
\end{equation}
On the other hand, for \(y\in[\delta,r]\),
\[
\begin{aligned}
\left|
\frac{g_n(y)-1}{g_n(y)}
-
\frac{g(y)-1}{g(y)}
\right|
&=
\left|
\frac{1}{g(y)}-\frac{1}{g_n(y)}
\right|\\
&\leq
\frac{1}{m_r^2}|g_n(y)-g(y)|.
\end{aligned}
\]
Hence
\begin{equation}\label{eqt:Psi_away_zero_bound}
\int_\delta^r
\left|
\frac{g_n(y)-1}{y g_n(y)}
-
\frac{g(y)-1}{y g(y)}
\right|
\idiff y
\leq
\frac{\log(r/\delta)}{m_r^2}
\|g_n-g\|_{L^\infty([0,r])}.
\end{equation}
Taking \(n\to+\infty\) first and then \(\delta\to0\), we obtain
\begin{equation}\label{eqt:Psi_local_convergence}
\|\Psi_n-\Psi\|_{L^\infty([0,r])}
\longrightarrow0.
\end{equation}

Since \(1/a>1\), the map \(s\mapsto s^{1/a}\) is Lipschitz on
\([0,1]\). Moreover, the exponential map is Lipschitz on
\((-\infty,0]\). It follows from
\eqref{eqt:G_local_convergence_for_R},
\eqref{eqt:Psi_nonpositive}, and
\eqref{eqt:Psi_local_convergence} that
\[
\begin{aligned}
\|h_n-h\|_{L^\infty([0,r])}
&\leq
\left\|g_n^{1/a}-g^{1/a}\right\|_{L^\infty([0,r])}
+
\left\|\econst^{\Psi_n}-\econst^\Psi\right\|_{L^\infty([0,r])}\\
&\lesssim_a
\|g_n-g\|_{L^\infty([0,r])}
+
\|\Psi_n-\Psi\|_{L^\infty([0,r])}
\longrightarrow0.
\end{aligned}
\]

It remains to treat a possible finite zero of \(g\). Recall the positive-set
endpoint \(\tau(f)\) from the definition of \(\calR\). If
\(\tau(f)=+\infty\), the preceding argument already proves that
\(h_n\to h\) uniformly on every compact interval.

Suppose that \(\tau(f)<+\infty\). By continuity and monotonicity,
\[
g(x)>0,\qquad 0\leq x<\tau(f),
\]
and
\[
g(x)=0,\qquad x\geq\tau(f).
\]
Fix \(R>\tau(f)\) and \(\varepsilon>0\). Choose \(r<\tau(f)\) sufficiently
close to \(\tau(f)\) so that
\[
g(r)^{1/a}<\frac{\varepsilon}{4}.
\]
By the local convergence of \(g_n\), for all sufficiently large \(n\),
\[
g_n(r)^{1/a}<\frac{\varepsilon}{2}.
\]
Using the monotonicity of \(g_n\) and \eqref{eqt:R_bounded_by_G}, we obtain
for \(x\in[r,R]\),
\[
0\leq h_n(x)
\leq
g_n(x)^{1/a}
\leq
g_n(r)^{1/a}
<
\frac{\varepsilon}{2}.
\]
Similarly,
\[
0\leq h(x)
\leq
g(x)^{1/a}
\leq
g(r)^{1/a}
<
\frac{\varepsilon}{4}.
\]
Consequently,
\[
\|h_n-h\|_{L^\infty([r,R])}
<
\frac{3\varepsilon}{4}
\]
for all sufficiently large \(n\). On \([0,r]\), the preceding argument gives
uniform convergence. Hence
\[
\|h_n-h\|_{L^\infty([0,R])}\longrightarrow0.
\]
Since \(R>0\) is arbitrary, we have proved
\begin{equation}\label{eqt:R_local_convergence}
\calR(f_n)\longrightarrow\calR(f)
\qquad\text{locally uniformly on }\mathbb R.
\end{equation}

We finally upgrade the local convergence to convergence in
\(L^\infty(\mathbb R)\). Since \(\calR\) maps \(\mathbb D\) into itself,
\[
0\leq h_n(x),h(x)\leq f_{a,\delta_a}^{\mathrm{up}}(x),
\qquad x\in\mathbb R.
\]
By the decay of the upper barrier,
\[
\lim_{R\to+\infty}
\sup_{|x|\geq R}f_{a,\delta_a}^{\mathrm{up}}(x)=0.
\]
Given \(\varepsilon>0\), choose \(R>0\) such that
\[
\sup_{|x|\geq R}f_{a,\delta_a}^{\mathrm{up}}(x)
<
\frac{\varepsilon}{4}.
\]
It follows that
\[
\sup_{|x|\geq R}|h_n(x)-h(x)|
\leq
2\sup_{|x|\geq R}f_{a,\delta_a}^{\mathrm{up}}(x)
<
\frac{\varepsilon}{2}.
\]
By \eqref{eqt:R_local_convergence}, for all sufficiently large \(n\),
\[
\|h_n-h\|_{L^\infty([-R,R])}
<
\frac{\varepsilon}{2}.
\]
Combining the compact and far-field estimates yields
\[
\|\calR(f_n)-\calR(f)\|_{L^\infty}
<
\varepsilon
\]
for all sufficiently large \(n\). This completes the proof.
\end{proof}

\subsection{Schauder existence and physical reconstruction}

\begin{theorem}\label{thm:fixed_point_existence}
For every \(a\in(0,1)\), the map \(\calR\) has a fixed point
\[
f_a\in\mathbb D,
\qquad
\calR(f_a)=f_a.
\]
\end{theorem}

\begin{proof}
By Lemma~\ref{lem:D_properties}, \(\mathbb D\) is a nonempty compact convex
subset of the Banach space \(\mathbb V\). Theorem~\ref{thm:R_invariant} and
Proposition~\ref{prop:R_continuity} show that \(\calR\) is a continuous map
from \(\mathbb D\) into itself. The Schauder fixed-point theorem therefore
gives the asserted fixed point.
\end{proof}

\begin{theorem}\label{thm:fixed_point_profile_reconstruction}
Let \(f\in\mathbb D\) be a fixed point of \(\calR\), and set
\[
L:=\ell(f),
\qquad
b:=b(f,L),
\qquad
c:=c(f,L),
\qquad
\bar c_l:=c_l(f,L).
\]
Define
\begin{equation}\label{eqt:reconstruction_normalization}
\kappa:=\bar c_l+ab=\frac{1-a/3}{2}c>0,
\qquad
\bar c_\omega:=\kappa-b=\bar c_l-(1-a)b.
\end{equation}
Then there are nontrivial odd functions \(\bar\Omega\) and \(\bar U\), with
\(\bar U(0)=0\), such that
\begin{equation}\label{eqt:reconstructed_profile_equation}
(\bar c_lx+a\bar U)\bar\Omega_x
=
(\bar c_\omega+\bar U_x)\bar\Omega,
\qquad
\bar U_x=H_L\bar\Omega.
\end{equation}
If \(L<+\infty\), the system is posed on \(\mathbb T_L\); if
\(L=+\infty\), it is posed on \(\mathbb R\). The identities hold almost
everywhere and hence as distributional identities. Moreover,
\begin{equation}\label{eqt:reconstructed_rate_sign}
\bar c_\omega<0.
\end{equation}

More precisely, exactly one of the following alternatives holds:
\begin{enumerate}
\item \(L<+\infty\), \(\bar c_l=0\), and \(\bar\Omega\) is the odd
\(2L\)-periodic extension of \(-xf(x)\) from \([-L,L]\). In this case
\[
f(L)=0,
\qquad
\bar c_\omega=-(1-a)b<0.
\]
\item \(L=+\infty\) and \(\bar c_l>0\). Then
\(\bar\Omega(x)=-xf(x)\) on \(\mathbb R\), and
\[
\bar c_\omega=\bar c_l-(1-a)b<0.
\]
\item \(L=+\infty\) and \(\bar c_l=0\). Then
\(\bar\Omega(x)=-xf(x)\) on \(\mathbb R\), and
\[
\bar c_\omega=-(1-a)b<0.
\]
\end{enumerate}
In every case,
\begin{equation}\label{eqt:reconstructed_origin_normalization}
f(x)=1-x^2+o(x^2),
\qquad
\bar\Omega(x)=-x+x^3+o(x^3)
\qquad\text{as }x\to0.
\end{equation}
\end{theorem}

\begin{proof}
The selected functionals are positive and finite. This is immediate from
their finite-period integral representations when \(L<+\infty\). When
\(L=+\infty\), scale selection gives
\begin{equation}\label{eqt:infinite_selected_integrability}
b(f,+\infty)\leq k_ac(f,+\infty)<+\infty,
\qquad
b(f,+\infty)=\frac2\pi\int_0^{+\infty}f(y)\idiff y.
\end{equation}
The identity and positivity in \eqref{eqt:reconstruction_normalization} follow
from the definition of \(c_l\) and the selected lower bound for
\(c\).

Suppose first that \(L<+\infty\). Corollary~\ref{cor:selected_G_structure}
gives \(\tau(f)=L\), while \(\calR(f)\) is continued by zero on
\([L,+\infty)\). Since \(f=\calR(f)\) in \(C(\mathbb R)\),
\begin{equation}\label{eqt:finite_fixed_point_terminal_zero}
f(L)=0,
\qquad
f(x)=0\quad(x\geq L).
\end{equation}
Consequently the endpoint-corrected iterate becomes the physical vorticity on
the fundamental cell:
\begin{equation}\label{eqt:endpoint_correction_reconciliation}
q_{f,L}(x)=x\big(f(L)-f(x)\big)=-xf(x),
\qquad -L\leq x\leq L.
\end{equation}
The endpoint values vanish, so this function has a continuous odd
\(2L\)-periodic extension, which we denote by \(\bar\Omega\).

In the infinite case, put \(\bar\Omega(x):=-xf(x)\) on \(\mathbb R\).
Equation \eqref{eqt:infinite_selected_integrability} gives
\[
\int_{\mathbb R}
\frac{|\bar\Omega(y)|}{1+|y|}\idiff y<+\infty.
\]
Thus \(H_{\mathbb R}\bar\Omega\) is well-defined locally as a principal
value almost everywhere and as a tempered distribution. At the origin,
\begin{equation}\label{eqt:reconstructed_H_origin}
H_{\mathbb R}\bar\Omega(0)
=
\frac2\pi\int_0^{+\infty}f(y)\idiff y=b.
\end{equation}

Let
\[
T:=\calT(f),
\qquad
g:=\calG(f),
\]
and define, on the fundamental cell in the finite case and on \(\mathbb R\)
in the infinite case,
\begin{equation}\label{eqt:reconstructed_velocity}
\bar U(x):=x\big(T(x)+b\big).
\end{equation}
The definition of \(\calT\), together with
\eqref{eqt:endpoint_correction_reconciliation} in the finite case, yields
\[
T(x)+b
=
\frac1x\int_0^xH_L\bar\Omega(z)\idiff z.
\]
Therefore \(\bar U(0)=0\) and
\begin{equation}\label{eqt:reconstructed_biot_savart}
\bar U_x=H_L\bar\Omega
\end{equation}
locally in distributions and almost everywhere. Since \(\bar\Omega\) is odd,
its Hilbert transform and \(T\) are even, and hence \(\bar U\) is odd. If
\(L<+\infty\), then
\eqref{eqt:T_L_endpoint} gives \(\bar U(\pm L)=0\), so \(\bar U\) also has
a continuous odd periodic extension.

The normalization \eqref{eqt:reconstruction_normalization}, the definition of
\(\calG\), and \eqref{eqt:reconstructed_velocity} give the exact identity
\begin{equation}\label{eqt:reconstructed_g_identity}
g(x)
=
\frac{\bar c_l+ab+aT(x)}{\kappa}
=
\frac{\bar c_lx+a\bar U(x)}{x\kappa},
\qquad x\neq0.
\end{equation}
In particular,
\[
\bar U_x(0)=b,
\qquad
\kappa=\bar c_l+a\bar U_x(0)
=\bar c_\omega+\bar U_x(0).
\]

On every compact subinterval of the positive component of \(g\), the fixed
point relation reads
\begin{equation}\label{eqt:fixed_point_log_identity}
\log f(x)
=
\frac1a\log g(x)
+
\frac{1-a}{a}
\int_0^x\frac{g(y)-1}{yg(y)}\idiff y.
\end{equation}
Here \(f>0\) because the defining expression for \(\calR\) is positive before
the terminal point. The functions involved are locally absolutely continuous,
so differentiation almost everywhere gives
\begin{equation}\label{eqt:fixed_point_log_derivative}
\frac{f_x}{f}
=
\frac1a\frac{g_x}{g}
+
\frac{1-a}{a}\frac{g-1}{xg}.
\end{equation}
Differentiating \eqref{eqt:reconstructed_g_identity} and using
\eqref{eqt:reconstruction_normalization}, we obtain
\[
\bar c_l+a\bar U_x=\kappa(g+xg_x),
\qquad
a\bar c_\omega-\bar c_l=-(1-a)\kappa.
\]
Hence
\[
\bar c_\omega+\bar U_x
=
\frac{\kappa}{a}\big(g+xg_x-(1-a)\big).
\]
On the other hand, \eqref{eqt:fixed_point_log_derivative} gives
\[
g\left(1+x\frac{f_x}{f}\right)
=
\frac1a\big(g+xg_x-(1-a)\big).
\]
Since \(\bar\Omega=-xf\) and
\(\bar c_lx+a\bar U=x\kappa g\), the last two identities are exactly the first
equation in \eqref{eqt:reconstructed_profile_equation} on the positive
interior. The identity on the negative half follows from the evenness of
\(f,g\) and the oddness of \(\bar\Omega,\bar U\). The origin expansion
\eqref{eqt:G_origin_expansion}, substituted into the fixed-point formula,
gives
\[
\log f(x)
=
-\frac1a\frac{2a}{3-a}x^2
-\frac{1-a}{2a}\frac{2a}{3-a}x^2
+o(x^2)
=-x^2+o(x^2),
\]
which proves \eqref{eqt:reconstructed_origin_normalization} and nontriviality.

It remains to prove the strict sign of \(\bar c_\omega\). If
\(\bar c_l=0\), then \eqref{eqt:reconstruction_normalization} gives
\[
\bar c_\omega=-(1-a)b<0.
\]
This covers both the finite branch and the infinite critical branch. Suppose
now that \(L=+\infty\) and \(\bar c_l>0\). The whole-line terminal value in
Corollary~\ref{cor:selected_G_structure} gives
\[
g_\infty
:=
\lim_{x\to+\infty}g(x)
=
\frac{\bar c_l}{\kappa}>0.
\]
Set
\[
r(x):=\frac{1-a}{a}\left(\frac1{g(x)}-1\right),
\qquad
\rho:=\frac{(1-a)b}{\bar c_l}.
\]
Since \(g\geq g_\infty\), one has \(0\leq r(x)\leq\rho\). The origin
expansion makes \(\int_0^1r(y)\idiff y/y\) finite. Thus the fixed-point
formula implies, for \(x\geq1\),
\[
f(x)
=
g(x)^{1/a}
\exp\left(-\int_0^x\frac{r(y)}y\idiff y\right)
\geq
C g_\infty^{1/a}x^{-\rho}
\]
for a constant \(C>0\). But \eqref{eqt:infinite_selected_integrability} says
that \(f\in L^1(0,+\infty)\), so necessarily \(\rho>1\). Consequently,
\[
\bar c_\omega
=
\bar c_l-(1-a)b
=
\bar c_l(1-\rho)<0.
\]

Finally, we justify the stated almost-everywhere and distributional meanings.
In the finite branch,
\(\bar\Omega\in W^{1,\infty}(\mathbb T_L)\); the matching zero endpoint
values prevent a jump measure in its periodic derivative. Boundedness of the
periodic Hilbert transform on \(W^{1,p}\) gives
\[
\bar U_x\in W^{1,p}(\mathbb T_L),
\qquad
\bar U\in W^{2,p}(\mathbb T_L),
\qquad 1<p<+\infty.
\]
In an infinite branch, monotonicity and integrability of \(f\) imply
\[
xf(x)\leq2\int_{x/2}^xf(y)\idiff y\longrightarrow0,
\]
and hence
\[
\bar\Omega\in W^{1,\infty}_{\mathrm{loc}}(\mathbb R)\cap C_0(\mathbb R),
\qquad
\bar U_x\in L^p_{\mathrm{loc}}(\mathbb R),
\quad 1<p<+\infty.
\]
For the last assertion, on a fixed compact interval split \(\bar\Omega\) into
a compactly supported part and a far part. The Hilbert transform of the first
part is locally in \(L^p\), while the far integral is bounded there by the
weighted \(L^1\) norm established above.
All products in \eqref{eqt:reconstructed_profile_equation} are locally
integrable. The almost-everywhere identity therefore also holds
distributionally; in the finite case, the matching endpoint values create no
periodic Dirac masses. The interior smoothness is proved in Theorem~\ref{thm:fixed_point_profile_regularity}.
\end{proof}

\section{Regularity and endpoint behavior}
\label{sec:regularity}

This section describes the physical profiles produced by the fixed-point
map. We first identify its whole-line members with HQWW profiles and
record the moment information needed below. We then establish the exact
endpoint powers, remove the apparent origin singularity in the
regularity equation, and determine the critical far field.
\subsection{Identification of the whole-line profiles}

The whole-line part of the construction has the same kernel and
normalization as HQWW. We record this correspondence before using their
moment identities; the periodic part of the argument keeps the
cotangent kernel throughout.

Let \(\mathbb D_{\mathrm H}\) denote the parameter-independent HQWW class:
its members are even continuous functions satisfying \(f(0)=1\),
\((1-x^2)_+\leq f(x)\leq1\), \(f'_-(1/2)\leq-\eta\),
monotonicity on \([0,\infty)\), and convexity in \(x^2\).
It is viewed in the weighted Banach space
\begin{equation}\label{eq:hqww_space}
 X=\{f\in C(\mathbb R): f(-x)=f(x),\ 
                \|f\|_X<\infty\},
 \qquad \|f\|_X:=\sup_x(1+|x|)^{-1/2}|f(x)|.
\end{equation}

\begin{proposition}[Correspondence with HQWW]\label{prop:hqww_correspondence}
Let \(0<a<1\). A fixed point of the present map with
\(\ell(f)=\infty\) is a normalized HQWW fixed point.
Conversely, a normalized HQWW fixed point with
\(\bar c_l\geq0\) belongs to the admissible set of the present
construction and has \(\ell(f)=\infty\).
\end{proposition}

\begin{proof}
For a whole-line input the functionals and rates are
\begin{equation}\label{eq:hqww_normalizations}
\begin{aligned}
 b(f)&=\frac2\pi\int_0^\infty f(x)\idiff x,
 &c(f)&=\frac2\pi\int_0^\infty\frac{1-f(x)}{x^2}\idiff x,\\
 \kappa&=\frac{1-a/3}{2}c(f),
 &\bar c_l&=\kappa-ab,\qquad \bar c_\omega=\kappa-b.
\end{aligned}
\end{equation}
The kernel \eqref{eqt:T_whole_line_kernel}, the intermediate function
\(g=1+aT/\kappa\), and the reconstruction
\eqref{eqt:R_map_definition} are exactly those in
\cite[Sections 2--3]{HQWW2024}. The class used here is contained in
\(\mathbb D_{\mathrm H}\), so this proves the first assertion.
For the converse, we first establish membership in \(\mathbb D\).
By \cite[Lemma 4.2]{HQWW2024}, \(b<\infty\). Since
\(\kappa-ab=\bar c_l\geq0\), the whole-line function
\(g=1+aT_\infty f/\kappa\) is positive at every finite \(x\),
nonincreasing, and convex in \(x^2\). Positivity also follows directly
from \(T_\infty f(x)+b>0\), and the shape properties are
\cite[Lemma 3.4]{HQWW2024}. Thus the whole-line reconstruction formula
is defined without making a period selection.

The proof of Lemma~\ref{lem:periodic_T_midpoint_bounds} at
\(L=\infty\) uses only these shape conditions, the lower cap, and the
midpoint slope. It therefore applies to this HQWW input and gives
\[
 -\frac{(T_\infty f)'(1/2)}{c(f)}
 \geq\frac{16\eta}{\pi^2+24},
 \qquad
 1-g(1/2)\geq\frac{24a\eta}{(3-a)(\pi^2+24)}.
\]
Integrating \eqref{eqt:R_logarithmic_derivative}, now with output
\(f\), gives
\(f(x)\leq f(1/2)(2x)^{-2\delta_a}\) for \(x\geq1/2\).
The supporting-line argument in Theorem~\ref{thm:R_invariant}
then yields \(f(x)\leq(1+x^2)^{-\delta_a}\) on the whole line.
All the other defining conditions of \(\mathbb D\) are already HQWW
shape conditions. Now that \(f\in\mathbb D\), the inequality
\(b/c\leq k_a\) and Lemma~\ref{lem:b_over_c_monotonicity}
give \(\ell(f)=\infty\). The two reconstruction formulas therefore
agree at this fixed point.
\end{proof}

We will use the following precise whole-line consequences of HQWW.
They replace the corresponding kernel computations for that already
constructed family.

\begin{proposition}[Whole-line moment information]
\label{thm:whole_line_Q_mu_identity}
For a normalized HQWW fixed point at \(0<a<1\), put
\begin{equation}\label{eqt:whole_line_Q_definition}
 \mathcal Q(f)=\frac1{\pi^2}\int_0^\infty\int_0^\infty f(x)f(y)
 \left(\left(\frac xy+\frac yx\right)
       \log\left|\frac{x+y}{x-y}\right|-2\right)\idiff x\idiff y,
 \qquad \mu(f)=\frac{2\mathcal Q(f)}{b(f)^2}.
\end{equation}
Then \(b<\infty\), \(b_1=(2/\pi)\int_0^\infty xf(x)\idiff x<\infty\),
\(\bar c_\omega<0\), and
\begin{align}
 0<\mu(f)&<\overline\mu:=\frac{9\pi^2}{64}-\frac34<1,
 \label{eqt:whole_line_mu_bounds}\\
 \frac{b(f)}{c(f)}&=\frac{1-a/3}{1+a\mu(f)},
 \label{eqt:whole_line_fixed_point_ratio}\\
 \bar c_l&=\frac b2(1-2a+a\mu(f)),
 \qquad \bar c_\omega=\frac b2(a\mu(f)-1).
 \label{eqt:whole_line_rates_from_mu}
\end{align}
\end{proposition}

\begin{proof}
The functionals \(b,c,\mathcal Q\) have exactly the normalization in
\cite[(3.3)--(3.5), (4.4)--(4.5)]{HQWW2024}. Their
Lemma~4.2 gives \(b<\infty\) and \(\bar c_\omega<0\),
Theorem~4.3 gives \eqref{eqt:whole_line_fixed_point_ratio} and
\(0<\mu\leq\overline\mu\), and Lemma~4.8 gives \(b_1<\infty\).
We supply the strictness argument because it will determine the signs
at the endpoints of the continuation interval.

Use the cap representation of Lemma~\ref{lem:convex_cap_mixture},
and write \(\mathcal Q(f,h)\) for the bilinear form with the kernel
in \eqref{eqt:whole_line_Q_definition}. The cap calculation in the
proof of \cite[Theorem 4.3]{HQWW2024} gives
\begin{equation}\label{eq:hqww_cap_pairing}
 \mathcal Q(\phi_r,\phi_s)=\frac{4rs}{\pi^2}F_4(r/s),
 \qquad \phi_r(x)=(1-x^2/r^2)_+.
\end{equation}
Here the scalar function \(F_4\), defined in
\cite[Appendix A.3]{HQWW2024}, satisfies \(F_4(t)=F_4(1/t)\), and
on \(0\leq t\leq1\) it has the convergent expansion
\begin{equation}\label{eq:hqww_cap_kernel_series}
 F_4(t)=\frac{\pi^2}{32}t+
 \sum_{j=1}^{\infty}
 \frac{8j\,t^{2j}}{(2j-3)(2j-1)^2(2j+1)^2(2j+3)},
 \qquad F_4(1)=\frac{\pi^2}{32}-\frac16.
\end{equation}
The same appendix proves that \(F_4\) is nondecreasing on \([0,1]\).
The series shows that it is real analytic and nonconstant on \((0,1)\).
It cannot therefore be constant on any subinterval, so
\(F_4(t)<F_4(1)\) whenever \(t\neq1\).

The nonnegative kernel permits Tonelli's theorem in the cap mixture.
Since \(b=(4/(3\pi))\int r\idiff\nu_f(r)\),
\eqref{eq:hqww_cap_pairing} yields
\begin{equation}\label{eq:hqww_mu_mixture}
 \mu(f)=\frac92
 \frac{\displaystyle\iint rsF_4(r/s)\idiff\nu_f(r)\idiff\nu_f(s)}
      {\displaystyle\left(\int r\idiff\nu_f(r)\right)^2}
 \leq\frac92F_4(1)=\overline\mu.
\end{equation}
Equality forces \(r=s\) for \(\nu_f\otimes\nu_f\)-almost every
\((r,s)\), and hence \(\nu_f\) is a point mass. Thus equality is
possible only for a single parabolic cap. Such a cap cannot be a fixed
point at \(a<1\): its linear vanishing contradicts the exponent
\(p\geq(3-a)/(2a)>1\) in
\cite[Theorem 4.10]{HQWW2024}. This proves the strict upper bound
in \eqref{eqt:whole_line_mu_bounds}. Substituting the ratio identity
into \eqref{eq:hqww_normalizations} proves
\eqref{eqt:whole_line_rates_from_mu}.
\end{proof}

\begin{corollary}[Whole-line tails and exclusion at large parameters]
\label{cor:whole_line_branch_restrictions}
For a fixed point of the period-selection map with \(L=\infty\),
the focusing exponent satisfies
\begin{equation}\label{eqt:focusing_tail_exponent_improved}
 \rho=\frac{2(1-a)}{1-2a+a\mu(f)}>2
 \qquad\text{if }\bar c_l>0.
\end{equation}
If \(\bar c_l=0\), then
\begin{equation}\label{eqt:critical_mu_constraint}
 \mu(f)=2-\frac1a,
 \qquad a<a_{\mathrm{per}}.
\end{equation}
If \(a_{\mathrm{per}}\leq a<1\), every fixed point of the
period-selection map has finite selected period.
\end{corollary}

\begin{proof}
In the focusing case, subtracting two from the fraction in
\eqref{eqt:focusing_tail_exponent_improved} gives
\(2a(1-\mu)/(1-2a+a\mu)>0\). In the critical case,
\eqref{eqt:whole_line_rates_from_mu} gives \(\mu=2-1/a\),
and \(\mu<\overline\mu\) gives
\(a<1/(2-\overline\mu)=a_{\mathrm{per}}\).
More generally, if \(a\geq a_{\mathrm{per}}\), the same strict
upper bound yields \(1-2a+a\mu<0\), contradicting
\(\bar c_l\geq0\) for a selected whole-line profile.
\end{proof}
\subsection{The finite-period endpoint and the whole-line tails}

\begin{proposition}\label{prop:finite_profile_seam_asymptotic}
Let \(f\in\mathbb D\) be a fixed point for which
\(L:=\ell(f)<+\infty\), and put \(g:=\calG(f)\). Define
\begin{equation}\label{eqt:finite_terminal_slope_constant}
\lambda_L:=-g_x(L-)
=
\frac1{bL^2}
\int_0^L
yf(y)\tan\left(\frac{\pi y}{2L}\right)\idiff y
>0
\end{equation}
and define the velocity slope \(\kappa_L:=-\bar U_x(L-)\).
Since \(g=\bar U/(bx)\) in the finite-period case, these endpoint
slopes satisfy \(\kappa_L=bL\lambda_L\). Set
\begin{equation}\label{eqt:finite_seam_exponent}
p_L
:=
\frac1a+\frac{1-a}{aL\lambda_L}.
\end{equation}
Then there is a constant \(A_{a,f,L}\in(0,+\infty)\) such that
\begin{equation}\label{eqt:finite_seam_asymptotic}
f(x)
=
A_{a,f,L}(L-x)^{p_L}(1+o(1))
\qquad\text{as }x\uparrow L.
\end{equation}
Moreover,
\begin{equation}\label{eqt:finite_seam_exponent_bound}
p_L\geq\frac{3-a}{2a}>1.
\end{equation}
If \(f_{\rm per}\) denotes the even \(2L\)-periodization of the cell trace
of \(f\), then, as \(h\to0\),
\begin{align}
f_{\rm per}(L+h)
&=
A_{a,f,L}|h|^{p_L}(1+o(1)),
\label{eqt:finite_periodic_f_seam}\\
\bar\Omega(L+h)
&=
LA_{a,f,L}\operatorname{sgn}(h)|h|^{p_L}(1+o(1)).
\label{eqt:finite_periodic_omega_seam}
\end{align}
In particular, \(f_{\rm per}\) and \(\bar\Omega\) are \(C^1\) across the
periodic joining point, with zero first derivative there. The complete regularity classification at this point is given in Theorem~\ref{thm:periodic-endpoint-regularity}.
\end{proposition}

\begin{proof}
The formula for \(\lambda_L\) is the terminal derivative formula in
Corollary~\ref{cor:selected_G_structure}, using
\eqref{eqt:finite_fixed_point_terminal_zero}. The integral is finite because
the Lipschitz bound on \(f\) gives \(f(y)=O(L-y)\), and it is strictly
positive because \(f>0\) on \((0,L)\).

The reconstruction proof gives
\(\bar U\in W^{2,q}(\mathbb T_L)\) for every \(1<q<+\infty\). In the finite
branch, \eqref{eqt:reconstructed_g_identity} reduces near \(L\) to
\(g(x)=\bar U(x)/(bx)\). Hence, for every fixed \(\alpha\in(0,1)\),
one may choose \(q\) large enough that
\[
g(L-t)=\lambda_Lt+O(t^{1+\alpha}),
\qquad
g_x(L-t)=-\lambda_L+O(t^\alpha).
\]
Consequently,
\[
\frac{g(y)-1}{yg(y)}
=
-\frac1{L\lambda_L(L-y)}
+O\big((L-y)^{\alpha-1}\big)+O(1).
\]
Thus the finite limit
\[
J_{g,L}
:=
\lim_{x\uparrow L}
\left(
\int_0^x\frac{g(y)-1}{yg(y)}\idiff y
-\frac1{L\lambda_L}\log(L-x)
\right)
\]
exists. Substitution into \eqref{eqt:fixed_point_log_identity} proves
\eqref{eqt:finite_seam_asymptotic}, with
\[
A_{a,f,L}
=
\lambda_L^{1/a}
\exp\left(\frac{1-a}{a}J_{g,L}\right)>0.
\]
The same expansion in \eqref{eqt:fixed_point_log_derivative} gives
\[
f_x(L-t)
=
-p_LA_{a,f,L}t^{p_L-1}(1+o(1)).
\]

Finally, convexity of \(s\mapsto g(\sqrt{s})\), together with
\(g(0)=1\) and \(g(L)=0\), gives
\[
\frac{g_x(L-)}{2L}
\geq
\frac{g(L)-g(0)}{L^2}
=
-\frac1{L^2}.
\]
Hence \(L\lambda_L\leq2\), which proves
\eqref{eqt:finite_seam_exponent_bound}. The even periodization of \(f\) and
the odd periodization of \(-xf\) now give
\eqref{eqt:finite_periodic_f_seam}--\eqref{eqt:finite_periodic_omega_seam};
the derivative asymptotic and \(p_L>1\) give the \(C^1\) matching.
\end{proof}

\begin{proposition}\label{prop:focusing_profile_exact_tail}
Let \(f\in\mathbb D\) be a fixed point for which
\(L=+\infty\) and \(\bar c_l>0\). Set
\begin{equation}\label{eqt:focusing_tail_exponent}
\rho
:=
\frac{(1-a)b}{\bar c_l}
=
1-\frac{\bar c_\omega}{\bar c_l}>1.
\end{equation}
Then there is a constant \(A_{a,f}\in(0,+\infty)\) such that
\begin{equation}\label{eqt:focusing_profile_exact_tail}
\lim_{x\to+\infty}x^\rho f(x)=A_{a,f}.
\end{equation}
Consequently,
\begin{equation}\label{eqt:focusing_vorticity_exact_tail}
\bar\Omega(x)
\sim
-A_{a,f}x^{1-\rho}
=
-A_{a,f}x^{\bar c_\omega/\bar c_l}
\qquad\text{as }x\to+\infty.
\end{equation}
\end{proposition}

\begin{proof}
Put
\[
P(x):=T(x)+b
=
\frac1{\pi x}
\int_0^{+\infty}
yf(y)\log\left|\frac{x+y}{x-y}\right|\idiff y>0,
\]
and define
\[
r(x):=\frac{1-a}{a}\left(\frac1{g(x)}-1\right),
\qquad
g_\infty:=\frac{\bar c_l}{\kappa}>0.
\]
By \eqref{eqt:reconstructed_g_identity},
\(g=(\bar c_l+aP)/\kappa\). Thus \(r\) is nondecreasing and
\(r(x)\to\rho\). The origin expansion makes
\(\int_0^1r(y)\idiff y/y\) finite. For every \(\varepsilon>0\), the fixed-point
formula therefore gives
\[
f(x)\leq C_\varepsilon x^{-(\rho-\varepsilon)}
\qquad(x\geq1).
\]
Since \(\rho>1\), choosing \(\varepsilon\in(0,\rho-1)\) yields
\begin{equation}\label{eqt:focusing_log_moment}
\int_0^{+\infty}f(y)\big(1+\log(1+y)\big)\idiff y<+\infty.
\end{equation}

The kernel defining \(P\) is nonnegative. Tonelli's theorem and the change
\(x=yu\) give
\begin{equation}\label{eqt:focusing_P_log_integral}
\begin{aligned}
\int_1^{+\infty}\frac{P(x)}x\idiff x
&=
\frac1\pi\int_0^{+\infty}f(y)J(y)\idiff y,\\
J(y)
&:=
\int_{1/y}^{+\infty}
u^{-2}\log\left|\frac{u+1}{u-1}\right|\idiff u.
\end{aligned}
\end{equation}
For \(0<y\leq2\), \(J(y)\) is bounded by the integral over
\([1/2,+\infty)\); the singularity at \(u=1\) is logarithmic. For \(y\geq2\),
split the integral at \(1/2\) and \(2\). On \((0,1/2)\),
\(\log((1+u)/(1-u))\leq Cu\), the middle integral is finite, and the
integrand is \(O(u^{-3})\) on \((2,+\infty)\). Hence
\[
J(y)\leq C\big(1+\log(1+y)\big).
\]
Equations \eqref{eqt:focusing_log_moment} and
\eqref{eqt:focusing_P_log_integral} show that \(P(x)/x\) is integrable at
infinity.

The exact branch algebra gives
\[
0\leq\rho-r(x)
=
\frac{(1-a)P(x)}{\kappa g(x)g_\infty}
\leq
\frac{1-a}{\kappa g_\infty^2}P(x).
\]
Therefore
\(\int_1^{+\infty}(\rho-r(x))\idiff x/x<+\infty\), and the fixed-point
identity can be rewritten as
\[
x^\rho f(x)
=
g(x)^{1/a}
\exp\left(-\int_0^1\frac{r(y)}y\idiff y\right)
\exp\left(\int_1^x\frac{\rho-r(y)}y\idiff y\right).
\]
Every factor has a finite strictly positive limit. This proves
\eqref{eqt:focusing_profile_exact_tail};
\eqref{eqt:focusing_vorticity_exact_tail} follows from
\(\bar\Omega=-xf\) and \eqref{eqt:focusing_tail_exponent}.
\end{proof}

\begin{proposition}\label{prop:critical_profile_gaussian_rate}
Let \(f\in\mathbb D\) be a fixed point for which
\(L=+\infty\) and \(\bar c_l=0\). Then
\begin{equation}\label{eqt:critical_superalgebraic_decay}
f(x)=O_N(x^{-N})
\qquad\text{as }x\to+\infty
\quad\text{for every }N>0.
\end{equation}
In particular, the number
\begin{equation}\label{eqt:critical_second_moment}
\mathsf b_2
:=
\frac2\pi\int_0^{+\infty}y^2f(y)\idiff y
\end{equation}
is finite and strictly positive. Moreover,
\begin{align}
\lim_{x\to+\infty}x^2g(x)
&=
\frac{\mathsf b_2}{b},
\label{eqt:critical_g_second_moment}\\
\lim_{x\to+\infty}\frac{\log f(x)}{x^2}
&=
-\frac{(1-a)b}{2a\mathsf b_2}<0.
\label{eqt:critical_gaussian_log_rate}
\end{align}
\end{proposition}

\begin{proof}
In this branch \(g(x)\downarrow0\). Since \(g_x\leq0\) almost everywhere,
\eqref{eqt:fixed_point_log_derivative} gives, for every \(N>0\),
\[
\frac{f_x(x)}{f(x)}\leq-\frac Nx
\]
for almost every sufficiently large \(x\). Integration proves
\eqref{eqt:critical_superalgebraic_decay}, and hence all polynomial moments
are finite.

Again write
\[
P(x):=T(x)+b
=
\frac1{\pi x}
\int_0^{+\infty}
yf(y)\log\left|\frac{x+y}{x-y}\right|\idiff y.
\]
We claim that
\begin{equation}\label{eqt:critical_P_second_moment}
x^2P(x)\longrightarrow\mathsf b_2.
\end{equation}
On \(0<y\leq x/2\),
\[
xy\log\left(\frac{x+y}{x-y}\right)\longrightarrow2y^2,
\qquad
0\leq
xy\log\left(\frac{x+y}{x-y}\right)
\leq Cy^2.
\]
Dominated convergence, with the moving indicator
\(\mathbf 1_{\{y\leq x/2\}}\), shows that this part of \(x^2P(x)\) tends to
\(\mathsf b_2\). On \(y\geq x/2\), put \(y=xt\) and use
\eqref{eqt:critical_superalgebraic_decay}. For every \(N>3\), the remaining
part is bounded by
\[
C_Nx^{3-N}
\int_{1/2}^{+\infty}
t^{1-N}\log\left|\frac{1+t}{1-t}\right|\idiff t=o(1).
\]
The last integral is finite because its singularity at \(t=1\) is logarithmic
and its integrand is \(O(t^{-N})\) at infinity. This proves
\eqref{eqt:critical_P_second_moment}.

Here \eqref{eqt:reconstructed_g_identity} reduces to \(g=P/b\), so
\eqref{eqt:critical_g_second_moment} follows. Set
\(A:=\mathsf b_2/b\) and
\[
h(y):=\frac{g(y)-1}{y^2g(y)}.
\]
Then \(h(y)\to-1/A\). For every fixed sufficiently large \(R\), weighted
Cesaro convergence gives
\[
\frac1{x^2}
\int_0^x\frac{g(y)-1}{yg(y)}\idiff y
=
o(1)+\frac1{x^2}\int_R^xyh(y)\idiff y
\longrightarrow
-\frac1{2A}.
\]
Also \(\log g(x)/x^2\to0\). Dividing
\eqref{eqt:fixed_point_log_identity} by \(x^2\) now proves
\eqref{eqt:critical_gaussian_log_rate}.
\end{proof}

\subsection{Interior smoothness and global Sobolev regularity}

The scalar equation appears singular at the origin because it contains \(T/x\). A commutator identity rewrites this quotient as bounded Hilbert and Hardy operators, so the Sobolev bootstrap includes the origin.

\begin{lemma}\label{lem:fixed_point_T_local_regularity}
Let \(f\in\mathbb D\) be a fixed point, set \(L:=\ell(f)\), and write
\(T:=\calT(f)\). For every integer \(m\geq0\),
\begin{equation}\label{eqt:T_local_regularity_implication}
f\in H^m_{\rm loc}((-L,L))
\quad\Longrightarrow\quad
T_x\in H^m_{\rm loc}((-L,L)),
\end{equation}
where \((-L,L)=\mathbb R\) when \(L=+\infty\).

More precisely, suppose first that \(L<+\infty\), and let \(f_{\rm per}\)
be the even \(2L\)-periodization of the cell trace of \(f\). For
\(-L<x<L\), define
\begin{align}
\mathcal C_Lf(x)
&:=
\frac1{2L}\int_{-L}^L
(y-x)f_{\rm per}(y)
\cot\left(\frac{\pi(x-y)}{2L}\right)\idiff y,
\label{eqt:periodic_sawtooth_commutator}\\
\mathcal E_Lf(x)
&:=
\frac{\mathcal C_Lf(x)-\mathcal C_Lf(0)}x,
\label{eqt:periodic_smoothing_remainder}
\end{align}
with the quotient extended at \(x=0\). Then
\(\mathcal E_Lf\in C^\infty_{\rm loc}((-L,L))\) is odd, and
\begin{align}
H_L\bar\Omega
&=
b-x\big(H_Lf_{\rm per}+\mathcal E_Lf\big),
\label{eqt:periodic_sawtooth_H_identity}\\
T_x(x)
&=
-K_L(x)+\int_0^1tK_L(tx)\idiff t,
\qquad
K_L:=H_Lf_{\rm per}+\mathcal E_Lf.
\label{eqt:periodic_T_Hardy_identity}
\end{align}
If \(L=+\infty\), the corresponding identities are
\begin{align}
H_{\mathbb R}\bar\Omega
&=
b-xH_{\mathbb R}f,
\label{eqt:whole_line_x_commutator}\\
T_x(x)
&=
-H_{\mathbb R}f(x)
+\int_0^1tH_{\mathbb R}f(tx)\idiff t.
\label{eqt:whole_line_T_Hardy_identity}
\end{align}
In this case one also has the global estimate
\begin{equation}\label{eqt:whole_line_T_global_Hm}
f\in H^m(\mathbb R)
\quad\Longrightarrow\quad
T_x\in H^m(\mathbb R),
\qquad
\|T_x\|_{H^m(\mathbb R)}
\leq C_m\|f\|_{H^m(\mathbb R)}.
\end{equation}
\end{lemma}

\begin{proof}
In the finite branch, let \(S_L\) denote the periodic sawtooth whose
representative is \(S_L(x)=x\) for \(-L<x<L\). Direct subtraction of the two
cotangent integrals gives
\[
[H_L,S_L]f_{\rm per}(x)
=
\mathcal C_Lf(x).
\]
There is no principal value in
\eqref{eqt:periodic_sawtooth_commutator}: the factor
\((y-x)\cot(\pi(x-y)/(2L))\) has a removable diagonal singularity. If
\(I\Subset(-L,L)\), then \(x-y\) stays a positive distance from the other
cotangent poles \(\pm2L\) on \(I\times[-L,L]\). All \(x\)-derivatives of the
extended kernel are therefore bounded there, so differentiation under the
\(L^1\) integral gives
\(\mathcal C_Lf\in C^\infty(I)\).

The function \(\mathcal C_Lf\) is even because \(f_{\rm per}\) is even.
Moreover, \(\bar\Omega=-xf_{\rm per}\) on the cell and
\[
\mathcal C_Lf(0)=-H_L\bar\Omega(0)=-b.
\]
It follows that
\(\mathcal C_Lf(x)-\mathcal C_Lf(0)=O(x^2)\) smoothly at the origin.
Thus \eqref{eqt:periodic_smoothing_remainder} extends smoothly and oddly
through zero. The commutator identity now gives
\eqref{eqt:periodic_sawtooth_H_identity}. Since
\[
(xT)_x=H_L\bar\Omega-b,
\]
integration from zero and differentiation of the resulting Hardy average
give \eqref{eqt:periodic_T_Hardy_identity}.

For completeness, the local Sobolev input can be seen directly by a cutoff.
If \(v\in L^2(\mathbb T_L)\cap H^m_{\rm loc}(O)\), choose a smooth periodic
cutoff supported in \(O\) and equal to one near a prescribed compact subset.
The periodic Hilbert transform maps the cutoff part in \(H^m(\mathbb T_L)\)
to \(H^m(\mathbb T_L)\); the kernel for the separated remainder is smooth.
Hence
\[
v\in H^m_{\rm loc}(O)
\quad\Longrightarrow\quad
H_Lv\in H^m_{\rm loc}(O).
\]
The Hardy operator
\[
(\mathcal Ph)(x):=\int_0^1th(tx)\idiff t
\]
is bounded on every \(H^m(-R,R)\). Indeed, for \(0\leq j\leq m\),
\[
\left\|\partial_x^j\mathcal Ph\right\|_{L^2(-R,R)}
\leq
\int_0^1t^{j+1-1/2}\idiff t\,
\left\|h^{(j)}\right\|_{L^2(-R,R)}.
\]
Together with \eqref{eqt:periodic_T_Hardy_identity}, this proves
\eqref{eqt:T_local_regularity_implication} for finite \(L\), including on
intervals containing the origin. Notice that the smoothing assertion for
\(\mathcal E_Lf\) is strictly interior-local and makes no claim at the
periodic joining point.

If \(L=+\infty\), direct subtraction of the whole-line kernels gives
\[
H_{\mathbb R}(xf)
=
xH_{\mathbb R}f-b.
\]
This proves \eqref{eqt:whole_line_x_commutator}; the same Hardy calculation
proves \eqref{eqt:whole_line_T_Hardy_identity}. Finally, selected-scale
integrability gives \(f\in L^1\cap L^2\). A compactly supported cutoff is
handled by the global \(H^m\) boundedness of the Hilbert transform, while
every derivative of the separated far kernel is dominated by the \(L^1\)
tail. This proves the whole-line part of
\eqref{eqt:T_local_regularity_implication}. On the full line, for
\(0\leq j\leq m\), Minkowski's inequality and scaling give
\[
\begin{aligned}
\left\|
\partial_x^j\int_0^1th(t\,\cdot)\idiff t
\right\|_{L^2(\mathbb R)}
&\leq
\int_0^1t^{j+1}
\left\|h^{(j)}(t\,\cdot)\right\|_{L^2(\mathbb R)}
\idiff t\\
&=
\left(\int_0^1t^{j+1-1/2}\idiff t\right)
\left\|h^{(j)}\right\|_{L^2(\mathbb R)}.
\end{aligned}
\]
Apply this with \(h=H_{\mathbb R}f\), using the boundedness of the
whole-line Hilbert transform on \(H^m\), to obtain
\eqref{eqt:whole_line_T_global_Hm}.
\end{proof}

\begin{theorem}\label{thm:fixed_point_profile_regularity}
Let \(f\in\mathbb D\) be a fixed point, with reconstructed profile
\((\bar\Omega,\bar U)\).
\begin{enumerate}
\item If \(L:=\ell(f)<+\infty\), then
\[
f,\bar\Omega,\bar U\in C^\infty((-L,L)).
\]
The conclusion includes the origin but not the periodic joining point. At the joining point,
\(f_{\rm per}\) and \(\bar\Omega\) are \(C^1\) and have the asymptotics in
Proposition~\ref{prop:finite_profile_seam_asymptotic}.
\item If \(L=+\infty\), then
\[
f,\bar\Omega,\bar U\in C^\infty(\mathbb R).
\]
If in addition \(\bar c_l>0\), then
\begin{equation}\label{eqt:focusing_global_Sobolev}
f\in H^m(\mathbb R)
\qquad\text{for every integer }m\geq0.
\end{equation}
If \(\bar c_l=0\), then likewise
\begin{equation}\label{eqt:critical_global_Sobolev}
f\in H^m(\mathbb R)
\qquad\text{for every integer }m\geq0.
\end{equation}
\end{enumerate}
\end{theorem}

\begin{proof}
Write \(D:=\kappa>0\). On the positive interior, the fixed-point logarithmic
derivative can be written as
\begin{equation}\label{eqt:fixed_point_regular_bootstrap_identity}
f_x
=
\frac1D
\left(
T_x+(1-a)\frac{T}{x}
\right)\frac fg.
\end{equation}
In the finite branch,
\eqref{eqt:periodic_sawtooth_H_identity} gives
\[
\frac{T}{x}
=
-H_Lf_{\rm per}-\mathcal E_Lf-T_x.
\]
In an infinite branch the same identity holds with
\(\mathcal E_Lf=0\) and \(H_Lf_{\rm per}\) replaced by
\(H_{\mathbb R}f\). Thus \eqref{eqt:fixed_point_regular_bootstrap_identity}
becomes
\begin{equation}\label{eqt:fixed_point_regular_nonsingular_identity}
f_x
=
\frac1D
\left(
aT_x-(1-a)
\big(H_Lf_{\rm per}+\mathcal E_Lf\big)
\right)\frac fg
\end{equation}
in the finite branch, with the analogous remainder-free formula on the line.
This identity has no division by \(x\).

Fix a compact subset of the spatial interior and enlarge it slightly. The
selected shape theorem gives a positive lower bound for \(g\) on the enlarged
compact set. Admissibility gives \(f\in W^{1,\infty}_{\rm loc}\), hence
\(f\in H^1_{\rm loc}\). Suppose inductively that
\(f\in H^m_{\rm loc}\) for an integer \(m\geq1\).
Lemma~\ref{lem:fixed_point_T_local_regularity} gives
\(T_x\in H^m_{\rm loc}\), and periodic or whole-line pseudolocality gives the
same regularity for the Hilbert term. Since
\[
g_x=\frac aD T_x,
\]
the one-dimensional \(H^m\) algebra and smooth-composition properties imply
\(f/g\in H^m_{\rm loc}\). Equation
\eqref{eqt:fixed_point_regular_nonsingular_identity} then gives
\(f_x\in H^m_{\rm loc}\), and hence \(f\in H^{m+1}_{\rm loc}\).
Induction and Sobolev embedding prove the asserted interior smoothness of
\(f\). The formulas
\(\bar\Omega=-xf\) and \(\bar U=x(T+b)\) give the same conclusion for the
physical profiles. Proposition~\ref{prop:finite_profile_seam_asymptotic}
supplies the stated initial endpoint matching; Theorem~\ref{thm:periodic-endpoint-regularity} gives its sharp refinement.

It remains to verify the two global whole-line statements. Suppose first that
\(\bar c_l>0\). Then \(g\geq g_\infty=\bar c_l/D>0\). Since
\(f\in L^1\cap L^\infty\), one has \(f\in L^2\). The global version of
\eqref{eqt:whole_line_T_Hardy_identity} shows that
\[
B:=aT_x-(1-a)H_{\mathbb R}f\in L^2.
\]
Also \(f/g\in L^\infty\), so the whole-line version of
\eqref{eqt:fixed_point_regular_nonsingular_identity} gives
\(f_x\in L^2\). Thus \(f\in H^1\). If \(f\in H^m\) for some \(m\geq1\),
\eqref{eqt:whole_line_T_global_Hm} and the Hilbert-transform bound give
\(B\in H^m\) and
\(g_x\in H^m\). Since \(m\geq1\), one-dimensional Sobolev embedding gives
\(\partial_x^jg\in L^\infty\) for \(1\leq j\leq m\). The lower bound
\(g\geq g_\infty>0\) and the chain rule therefore give
\(\partial_x^j(1/g)\in L^\infty\) for \(0\leq j\leq m\).
Leibniz' rule, with \(f\in H^m\), now gives \(f/g\in H^m\).
Since \(B,f/g\in H^m\) and the one-dimensional \(H^m\) space is an algebra,
the same identity gives \(f_x\in H^m\), closing the induction and proving
\eqref{eqt:focusing_global_Sobolev}.

Finally, suppose that \(\bar c_l=0\). Since \(0<g\leq1\), the fixed-point
formula and the nonpositivity of its integral give
\[
0<\frac fg
=
g^{1/a-1}
\exp\left(
\frac{1-a}{a}
\int_0^{|x|}
\frac{g(y)-1}{yg(y)}\idiff y
\right)
\leq1.
\]
Again \(B\in L^2\), so the nonsingular identity yields
\(f_x\in L^2\), and hence \(f\in H^1\).

To close the higher-order induction despite \(g(x)\to0\), define, for every
integer \(p\geq0\),
\[
q_p:=\frac{f}{g^p},
\qquad
\beta:=\frac{1-a}{a}.
\]
The critical asymptotics
\eqref{eqt:critical_g_second_moment} and
\eqref{eqt:critical_gaussian_log_rate}, extended to the negative half-line by
the evenness of \(f\) and \(g\), imply that for some \(C_0\geq1\),
\[
f(x)\leq C_0\exp(-|x|^2/C_0),
\qquad
g(x)\geq\frac1{C_0(1+x^2)}
\]
outside a compact interval. Since \(g\) has a positive minimum on every
compact interval, it follows that
\begin{equation}\label{eqt:critical_quotient_Linfinity}
q_p\in L^\infty(\mathbb R)
\qquad\text{for every integer }p\geq0.
\end{equation}

The logarithmic derivative identity
\eqref{eqt:fixed_point_log_derivative} gives, on the positive half-line,
\[
\frac{f_x}{f}
=
\frac1a\frac{g_x}{g}
+
\beta\frac{g-1}{xg}.
\]
Both sides extend to the line by parity. Therefore
\begin{equation}\label{eqt:critical_quotient_derivative}
(q_p)_x
=
q_{p+1}A_p,
\qquad
A_p
:=
\left(\frac1a-p\right)g_x
+
\beta\frac{g-1}{x},
\end{equation}
where the last quotient is extended smoothly at the origin.

Suppose inductively that \(f\in H^m(\mathbb R)\) for an integer \(m\geq1\).
Equation \eqref{eqt:whole_line_T_global_Hm} gives
\[
g_x=\frac aD T_x\in H^m(\mathbb R).
\]
Thus \(\partial_x^jg_x\in L^\infty\) for \(0\leq j\leq m-1\). Moreover,
\begin{equation}\label{eqt:critical_g_difference_quotient}
\frac{g(x)-1}{x}
=
\int_0^1g_x(tx)\idiff t,
\qquad
\partial_x^j\left(\frac{g-1}{x}\right)
=
\int_0^1t^j\partial_x^jg_x(tx)\idiff t.
\end{equation}
Consequently, \(A_p\in W^{m-1,\infty}(\mathbb R)\) for every \(p\).

We claim simultaneously for all \(p\geq0\) that
\begin{equation}\label{eqt:critical_quotient_Wminfinity}
q_p\in W^{m,\infty}(\mathbb R).
\end{equation}
The order-zero statement is
\eqref{eqt:critical_quotient_Linfinity}. If all \(q_p\) belong to
\(W^{r,\infty}\) for some \(0\leq r\leq m-1\), then
\eqref{eqt:critical_quotient_derivative} and the preceding bound on \(A_p\)
give \((q_p)_x\in W^{r,\infty}\) for every \(p\). Hence all \(q_p\) belong
to \(W^{r+1,\infty}\), proving the claim by induction on \(r\).

Finally,
\[
B=aT_x-(1-a)H_{\mathbb R}f\in H^m(\mathbb R).
\]
Since multiplication by a \(W^{m,\infty}\) function is bounded on integer
\(H^m\), the whole-line version of
\eqref{eqt:fixed_point_regular_nonsingular_identity}, together with
\eqref{eqt:critical_quotient_Wminfinity} for \(p=1\), gives
\(f_x=D^{-1}Bq_1\in H^m\). Thus \(f\in H^{m+1}\). Starting from the proved
base case \(f\in H^1\), induction proves
\eqref{eqt:critical_global_Sobolev}.
\end{proof}

\subsection{Sharp regularity at the periodic endpoint}
\begin{theorem}[Regularity across the periodic endpoint]
\label{thm:periodic-endpoint-regularity}
Let \(f\) be a finite-period fixed point and retain the notation of
Proposition~\ref{prop:finite_profile_seam_asymptotic}. Write
\(p=p_L>1\). The periodic physical profiles have the following
regularity.
\begin{enumerate}
\item If \(p=m+\alpha\notin\mathbb N\), with
\(m\geq1\) and \(0<\alpha<1\), then
\begin{equation}\label{eq:finite-noninteger-regularity}
 \bar\Omega\in C^{m,\alpha}(\mathbb T_L),
 \qquad
 \bar U\in C^{m+1,\alpha}(\mathbb T_L).
\end{equation}
The first inclusion is sharp at the periodic endpoint:
\(\bar\Omega\notin C^{m,\beta}\) for
\(\alpha<\beta\leq1\).
\item If \(p=2N\) is even, then
\begin{equation}\label{eq:finite-even-regularity}
 \bar\Omega\in C^{2N-1,1}(\mathbb T_L)
                   \setminus C^{2N}(\mathbb T_L),
 \qquad
 \bar U\in C_*^{2N+1}(\mathbb T_L).
\end{equation}
\item If \(p\) is an odd integer, then
\begin{equation}\label{eq:finite-odd-regularity}
 \bar\Omega,\bar U\in C^\infty(\mathbb T_L).
\end{equation}
\end{enumerate}
\end{theorem}

\begin{proof}
Only the endpoint requires consideration, by
Theorem~\ref{thm:fixed_point_profile_regularity}. Use the local
coordinate \(s=x-L\), with periodic extension understood. Oddness and
periodicity give
\(\bar\Omega(L+s)=-\bar\Omega(L-s)\) and
\(\bar U(L+s)=-\bar U(L-s)\). Moreover,
\eqref{eqt:finite_terminal_slope_constant} gives
\begin{equation}\label{eq:finite-velocity-endpoint-slope}
 \bar U'(L)=-bL\lambda_L<0.
\end{equation}
For sufficiently small \(|s|\), define the even functions
\begin{equation}\label{eq:finite-regular-singular-coefficient}
 R(s)=\frac{\bar U(L+s)}s
      =\int_0^1\bar U'(L+ts)\,\mathrm{d}t,
 \qquad
 C(s)=\frac{\bar U'(L+s)+\bar c_\omega}{aR(s)}.
\end{equation}
The quotient \(R\) has the nonzero limit \(\bar U'(L)\), and
\eqref{eqt:finite_seam_exponent} gives \(C(0)=p\).
The profile equation therefore becomes
\begin{equation}\label{eq:finite-regular-singular-ode}
 s\,\partial_s\bar\Omega(L+s)=C(s)\bar\Omega(L+s).
\end{equation}
Its positive amplitude is determined by
\eqref{eqt:finite_periodic_omega_seam}. Integrating
\eqref{eq:finite-regular-singular-ode} gives the exact local factorization
\begin{equation}\label{eq:finite-amplitude-factorization}
 \bar\Omega(L+s)=\operatorname{sgn}(s)|s|^p A(s),
 \qquad
 A(s)=LA_{a,f,L}
       \exp\left(\int_0^1\frac{C(ts)-p}{t}\,\mathrm{d}t\right).
\end{equation}
The function \(A\) is even and positive.

We record the regularity of this amplitude explicitly. If
\(C\in C^{k,\alpha}\), with \(k\in\mathbb N_0\) and
\(0<\alpha<1\), the integral in
\eqref{eq:finite-amplitude-factorization} also belongs to
\(C^{k,\alpha}\). For \(k=0\), this follows from
\begin{equation}\label{eq:finite-amplitude-holder}
 \left|\int_0^1\frac{C(ts)-C(tr)}t\,\mathrm{d}t\right|
 \leq\frac{[C]_{C^{0,\alpha}}}{\alpha}|s-r|^\alpha.
\end{equation}
For \(1\leq j\leq k\), its derivative of order \(j\) is
\begin{equation}\label{eq:finite-amplitude-derivatives}
 \int_0^1t^{j-1}C^{(j)}(ts)\,\mathrm{d}t,
\end{equation}
whose H\"older bound follows by integrating \(t^{j-1+\alpha}\).
Consequently \(A\in C^{k,\alpha}\).

Suppose now that \(\bar\Omega\in C^{k,\alpha}(\mathbb T_L)\).
The periodic Hilbert transform gives
\(\bar U'\in C^{k,\alpha}\). The first formula in
\eqref{eq:finite-regular-singular-coefficient} and its nonzero
endpoint value then imply \(R,C\in C^{k,\alpha}\) locally.
Thus \(A\in C^{k,\alpha}\) as well. Differentiating
\eqref{eq:finite-amplitude-factorization} through
\eqref{eq:finite-regular-singular-ode} yields
\begin{equation}\label{eq:finite-endpoint-bootstrap}
 \partial_s\bar\Omega(L+s)=|s|^{p-1}C(s)A(s).
\end{equation}
This identity provides the extra derivative.

For noninteger \(p=m+\alpha\), the function
\(|s|^{p-1}\) belongs to \(C^{m-1,\alpha}\).
The initial \(C^1\) bound makes \(\bar\Omega\) Lipschitz,
and therefore \(C^{0,\alpha}\); no integer-endpoint Hilbert
transform bound is needed. Apply
\eqref{eq:finite-endpoint-bootstrap} successively for
\(k=0,\ldots,m-1\). This proves
\(\bar\Omega\in C^{m,\alpha}\), and the Hilbert relation gives
\(\bar U\in C^{m+1,\alpha}\).
If \(p=n\) is even, the same induction first gives
\(\bar\Omega\in C^{n-1,\alpha}\) for every \(\alpha<1\).
The coefficient \(CA\) in \eqref{eq:finite-endpoint-bootstrap}
therefore belongs to \(C^{n-1,\alpha}\subset C^{n-2,1}\).
Since \(|s|^{n-1}\in C^{n-2,1}\), the identity yields
\(\bar\Omega\in C^{n-1,1}\). The relation
\(-\Lambda\bar U=\bar\Omega\) gives the asserted Zygmund
regularity of \(\bar U\).
If \(p\) is odd, then \(|s|^{p-1}=s^{p-1}\) is a polynomial,
so the same induction continues to every derivative order.

Finally, the positive value \(A(0)=LA_{a,f,L}\) in
\eqref{eq:finite-amplitude-factorization} gives the sharp exclusions
by the same Taylor argument as
\eqref{eq:steady-origin-cusp}. For even \(p\), the two limits of
\(\bar\Omega(L+s)/s^p\) are opposite; for noninteger \(p\),
a higher H\"older exponent would contradict the nonzero leading
power. This proves all assertions.
\end{proof}

\begin{remark}[The auxiliary periodic profile]
The periodized scalar function \(f_{\rm per}\) has different parity
at the endpoint from the physical vorticity. Formula
\eqref{eq:finite-amplitude-factorization} gives
\begin{equation}\label{eq:finite-auxiliary-amplitude}
 f_{\rm per}(L+s)=\frac{|s|^p A(s)}{L-|s|}
 \qquad(0<|s|\ll1).
\end{equation}
For noninteger \(p=m+\alpha\), this function belongs to
\(C^{m,\alpha}\), with sharp exponent at the endpoint.
When \(p\) is an odd integer, \(A\) is smooth and even, but
\eqref{eq:finite-auxiliary-amplitude} gives
\(f_{\rm per}\in C^{p-1,1}\setminus C^p\).
Thus smoothness of the physical periodic vorticity in
\eqref{eq:finite-odd-regularity} does not imply smoothness of this
auxiliary even periodization.
\end{remark}

\subsection{Weighted derivatives and the critical far field}
\begin{lemma}[Whole-line finite-moment potential tail]
\label{lem:whole_line_finite_moment_potential_tail}
Let \(f\in\mathbb D\) be a fixed point with
\(L:=\ell(f)=+\infty\), and suppose that
\[
\mathsf b_2
:=
\frac2\pi\int_0^{+\infty}y^2f(y)\idiff y
<+\infty.
\]
Put \(b:=b(f,+\infty)\), \(T:=\calT(f)\), and \(P:=T+b\). Then
\begin{equation}\label{eqt:whole_line_finite_moment_P_limit}
\lim_{x\to+\infty}x^2P(x)=\mathsf b_2,
\end{equation}
and, for every \(x>0\),
\begin{equation}\label{eqt:whole_line_finite_moment_P_bounds}
0<P(x)\leq\frac{3\mathsf b_2}{x^2},
\qquad
|P_x(x)|=|T_x(x)|\leq\frac{24\mathsf b_2}{x^3}.
\end{equation}
Consequently, the reconstructed velocity satisfies
\begin{equation}\label{eqt:whole_line_finite_moment_velocity_tail}
\lim_{x\to+\infty}x\bar U(x)=\mathsf b_2,
\qquad
0<\bar U(x)\leq\frac{3\mathsf b_2}{x}.
\end{equation}
If, in addition, \(\bar c_l>0\), then
\begin{equation}\label{eqt:focusing_first_derivative_asymptotic}
\lim_{x\to+\infty}x^{\rho+1}f_x(x)
=
-\rho A_{a,f},
\end{equation}
where \(\rho\) and \(A_{a,f}\) are defined in
Proposition~\ref{prop:focusing_profile_exact_tail}.
\end{lemma}

\begin{proof}
Monotonicity of \(f\) gives
\[
\int_{x/2}^xy^2f(y)\idiff y
\geq
\frac7{24}x^3f(x).
\]
The tail of the finite second-moment integral therefore shows that
\(x^3f(x)\to0\). In particular, \(x^3f(x)\) is bounded.

The positive whole-line kernel gives
\[
x^2P(x)
=
\frac1\pi\int_0^{+\infty}
f(y)xy\log\left|\frac{x+y}{x-y}\right|\idiff y.
\]
After adding and subtracting \(2y^2\) on \(0<y\leq x\) and setting
\(y=tx\) in the remainder, this becomes
\[
\begin{aligned}
x^2P(x)
={}&
\frac2\pi\int_0^xy^2f(y)\idiff y
\\
&+
\frac1\pi\int_0^{+\infty}
f(tx)(tx)^3\mathcal K_2(t)\idiff t,
\end{aligned}
\]
where
\[
\mathcal K_2(t)
:=
\frac1{t^2}\log\left|\frac{1+t}{1-t}\right|
-\frac2t\mathbf 1_{\{t\leq1\}}.
\]
The function \(\mathcal K_2\) is absolutely integrable: it is \(O(t)\)
at zero, has an integrable logarithmic singularity at one, and is
\(O(t^{-3})\) at infinity. For each fixed \(t>0\), the boundedness and
pointwise convergence of \(f(tx)(tx)^3\) allow dominated convergence in the remainder, while
monotone convergence applies to the first term. This proves
\eqref{eqt:whole_line_finite_moment_P_limit}.

For the integration by parts, define
\[
 \mathsf F_0(t):=\frac{t^2-1}{2t}
             \log\left|\frac{t+1}{t-1}\right|+1.
\]
The value at \(t=1\) is its continuous limit. Splitting at its logarithmic
singularity and integrating by parts gives
\[
P(x)
=
-\frac1\pi\int_0^{+\infty}
f_y(y)y\mathsf F_0(y/x)\idiff y.
\]
The artificial boundary terms at \(y=x\) cancel because
\(\mathsf F_0\) is continuous there.
The endpoint terms vanish because \(yf(y)\to0\) and \(y^3f(y)\to0\).
The power series on \((0,1]\) and the inversion identity give
\[
0\leq\mathsf F_0(t)\leq2t^2
\qquad(t>0).
\]
Indeed, \(\mathsf F_0(t)\leq t^2\) for \(t\leq1\), while
\(\mathsf F_0(t)=2-\mathsf F_0(t^{-1})\leq2\leq2t^2\) for \(t\geq1\).
Since \(f_y\leq0\), another integration by parts yields
\[
0<P(x)
\leq
-\frac2{\pi x^2}\int_0^{+\infty}y^3f_y(y)\idiff y
=
\frac{3\mathsf b_2}{x^2}.
\]

It remains to bound the derivative without differentiating this inequality.
Let \(V(s):=P(\sqrt{s})\). Corollary~\ref{cor:selected_G_structure}
shows in the selected whole-line case that \(V\) is nonincreasing and
convex. Since the
fixed point is smooth, for \(s>0\),
\[
0\leq-V'(s)
\leq
\frac{2}{s}\big(V(s/2)-V(s)\big)
\leq
\frac{2V(s/2)}s
\leq
\frac{12\mathsf b_2}{s^2}.
\]
Taking \(s=x^2\) proves the derivative bound in
\eqref{eqt:whole_line_finite_moment_P_bounds}. Since \(\bar U=xP\),
\eqref{eqt:whole_line_finite_moment_velocity_tail} follows immediately.

Finally, suppose that \(\bar c_l>0\), and put
\(g_\infty:=\bar c_l/\kappa>0\). The reconstruction identity and the bounds
above give
\[
g=g_\infty+\frac a{\kappa}P
=g_\infty+O(x^{-2}),
\qquad
g_x=O(x^{-3}).
\]
The logarithmic derivative identity therefore gives
\[
\frac{f_x}{f}
=
\frac1a\frac{g_x}{g}
+
\frac{1-a}{a}\frac{g-1}{xg}
=
-\frac\rho x+O(x^{-3}).
\]
Combining this with
\eqref{eqt:focusing_profile_exact_tail} proves
\eqref{eqt:focusing_first_derivative_asymptotic}.
\end{proof}

\begin{corollary}[Weighted whole-line regularity]
\label{cor:whole_line_weighted_regularity}
Let \(f\in\mathbb D\) be a fixed point with
\(\ell(f)=+\infty\), and retain the reconstructed profiles
\(\bar\Omega=-xf\) and \(\bar U\).
\begin{enumerate}
\item If \(\bar c_l>0\), then for every integer \(k\geq0\),
\begin{equation}\label{eqt:focusing_all_derivative_decay}
\left|\partial_x^kf(x)\right|
\leq
C_k(1+|x|)^{-\rho},
\end{equation}
where \(\rho>2\) is given by
\eqref{eqt:focusing_tail_exponent_improved}. Consequently,
\begin{equation}\label{eqt:focusing_vorticity_weighted_regularity}
\bar\Omega\in W^{m,p}(\mathbb R)
\quad
\text{for every integer }m\geq0
\text{ and every }1\leq p\leq+\infty,
\end{equation}
and \(\bar U_x\in H^m(\mathbb R)\) for every integer \(m\geq0\).
If \(\mathsf b_2<+\infty\), then
\begin{equation}\label{eqt:focusing_velocity_finite_moment_regularity}
\bar U\in H^m(\mathbb R)
\quad\text{for every integer }m\geq0,
\qquad
\bar U\notin L^1(\mathbb R).
\end{equation}
The velocity conclusion in \eqref{eqt:focusing_velocity_finite_moment_regularity} uses the finite second moment.
\item If \(\bar c_l=0\), then
\begin{equation}\label{eqt:critical_Schwartz_regularity}
f,\bar\Omega\in\mathcal S(\mathbb R).
\end{equation}
Moreover,
\begin{equation}\label{eqt:critical_velocity_regular_but_not_Schwartz}
\bar U\in H^m(\mathbb R)
\quad\text{for every integer }m\geq0,
\qquad
\bar U\notin L^1(\mathbb R),
\qquad
\bar U\notin\mathcal S(\mathbb R).
\end{equation}
\end{enumerate}
\end{corollary}

\begin{proof}
The quotient recurrence used in the critical part of
Theorem~\ref{thm:fixed_point_profile_regularity} is valid in either
whole-line branch. For every integer \(n\geq0\), put
\[
q_n:=\frac{f}{g^n},
\qquad
A_n
:=
\left(\frac1a-n\right)g_x
+
\frac{1-a}{a}\frac{g-1}{x}.
\]
Then, on the full line by parity,
\begin{equation}\label{eqt:whole_line_weighted_quotient_recurrence}
(q_n)_x=q_{n+1}A_n.
\end{equation}
The global Sobolev conclusions in
Theorem~\ref{thm:fixed_point_profile_regularity}, together with
\eqref{eqt:whole_line_T_global_Hm}, imply that
\(g_x\in H^m(\mathbb R)\) for every integer \(m\). Hence every derivative
of \(g_x\) is bounded. The identity
\[
\partial_x^j\left(\frac{g-1}{x}\right)
=
\int_0^1t^j\partial_x^jg_x(tx)\idiff t
\]
then shows that every derivative of every \(A_n\) is bounded.

In the positive-\(\bar c_l\) branch, \(g\) has a positive floor, and
Proposition~\ref{prop:focusing_profile_exact_tail} gives
\[
|q_n(x)|\leq C_n(1+|x|)^{-\rho}
\qquad(n\geq0).
\]
Differentiate \eqref{eqt:whole_line_weighted_quotient_recurrence} repeatedly.
If the derivatives through order \(k\) of all \(q_n\) have this decay, then
\[
\partial_x^{k+1}q_n
=
\sum_{j=0}^k
\binom{k}{j}
(\partial_x^jq_{n+1})
(\partial_x^{k-j}A_n)
\]
has the same decay. Simultaneous induction in \(n\) proves
\eqref{eqt:focusing_all_derivative_decay}. Since
\[
\partial_x^m(xf)
=
x\partial_x^mf+m\partial_x^{m-1}f
\qquad(m\geq1),
\]
each derivative of \(xf\) is bounded by a constant times
\((1+|x|)^{1-\rho}\). This weight is integrable because \(\rho>2\), which
proves \eqref{eqt:focusing_vorticity_weighted_regularity} for every
\(1\leq p\leq+\infty\). The whole-line Hilbert-transform bound applied to
\(\bar U_x=H_{\mathbb R}\bar\Omega\) gives the stated conclusion for
\(\bar U_x\). If \(\mathsf b_2<+\infty\),
Lemma~\ref{lem:whole_line_finite_moment_potential_tail} gives
\(\bar U(x)\sim\mathsf b_2/x\). Thus \(\bar U\in L^2\setminus L^1\),
and combining \(\bar U\in L^2\) with
\(\bar U_x\in H^m\) for every \(m\) proves
\eqref{eqt:focusing_velocity_finite_moment_regularity}.

In the critical branch, the Gaussian estimate and the positive quadratic
asymptotic of \(g\) give, for every fixed \(n\),
\[
|q_n(x)|
\leq
C_n(1+x^2)^n\exp(-\kappa_0|x|^2)
\]
for some \(\kappa_0>0\).
The same simultaneous derivative induction, now in the class of rapidly
decreasing functions, proves that every derivative of every \(q_n\) decays
faster than any power. Taking \(n=0\) gives
\(f\in\mathcal S(\mathbb R)\), and multiplication by \(x\) gives
\eqref{eqt:critical_Schwartz_regularity}.

Finally, \eqref{eqt:critical_P_second_moment} and
\(\bar U=xP\) yield
\[
\bar U(x)\sim\frac{\mathsf b_2}{x}
\qquad(x\to+\infty),
\]
where \(\mathsf b_2>0\). Thus \(\bar U\in L^2\setminus L^1\) and is not
Schwartz. On the other hand,
\(\bar U_x=H_{\mathbb R}\bar\Omega\in H^m\) for every \(m\). Combining
this with \(\bar U\in L^2\) proves
\eqref{eqt:critical_velocity_regular_but_not_Schwartz}.
\end{proof}

\begin{proposition}[The complete leading critical tail]
\label{prop:critical-full-tail}
Let \(f\) be a critical whole-line fixed point, so that
\(\bar c_l=0\), and retain \(b\) and \(\mathsf b_2\) from
\eqref{eqt:critical_second_moment}. Define the finite fourth moment by
\begin{equation}\label{eq:critical-fourth-moment}
 \mathsf b_4=\frac2\pi\int_0^\infty y^4f(y)\,\mathrm{d}y.
\end{equation}
Set
\begin{equation}\label{eq:critical-tail-parameters}
 \delta=\frac{(1-a)b}{2a\mathsf b_2},
 \qquad
 \nu=-\frac{1+a}{a}
       +\frac{(1-a)b\mathsf b_4}{3a\mathsf b_2^2}.
\end{equation}
Then there is a constant \(A_\infty>0\) such that
\begin{equation}\label{eq:critical-full-tail}
 f(x)=A_\infty x^\nu\exp(-\delta x^2)
                       \bigl(1+O(x^{-2})\bigr)
 \qquad(x\to+\infty).
\end{equation}
Consequently, the vorticity satisfies
\begin{equation}\label{eq:critical-vorticity-full-tail}
 \bar\Omega(x)=-A_\infty x^{\nu+1}\exp(-\delta x^2)
                       \bigl(1+O(x^{-2})\bigr).
\end{equation}
The expansion on the negative half-line follows by parity.
\end{proposition}

\begin{proof}
Corollary~\ref{cor:whole_line_weighted_regularity} gives
\(f\in\mathcal S(\mathbb R)\), so all polynomial moments are finite.
The positive potential
\(P=T+b=\bar U/x\) has the expansion
\begin{equation}\label{eq:critical-potential-fourth-order}
 P(x)=\frac{\mathsf b_2}{x^2}
       +\frac{\mathsf b_4}{3x^4}+O(x^{-6}).
\end{equation}
To prove this formula, split its defining integral at \(y=x/2\).
On \(0\leq y\leq x/2\), use
\begin{equation}\label{eq:critical-kernel-fourth-order}
 \log\frac{x+y}{x-y}
 =2\frac yx+\frac23\frac{y^3}{x^3}
                       +O\left(\frac{y^5}{x^5}\right).
\end{equation}
Multiplication by \(yf(y)/(\pi x)\) and integration give the two
moments in \eqref{eq:critical-potential-fourth-order}, with an
\(O(x^{-6})\) error controlled by the sixth moment.
For the complementary region, use \(y=xt\) and
\(f(y)\leq C_Ny^{-N}\). Its contribution to \(P\) is at most
\begin{equation}\label{eq:critical-potential-tail-error}
 C_Nx^{1-N}\int_{1/2}^\infty
       t^{1-N}\log\left|\frac{1+t}{1-t}\right|\,\mathrm{d}t.
\end{equation}
The integral is finite for \(N>1\). Choosing \(N>7\) also makes
the omitted tails of the second and fourth moments negligible at the
required order. This proves \eqref{eq:critical-potential-fourth-order}.

In the critical branch, \(g=P/b\). Inverting
\eqref{eq:critical-potential-fourth-order} gives
\begin{equation}\label{eq:critical-inverse-potential-expansion}
 \frac1{g(x)}
 =\frac b{\mathsf b_2}x^2
       -\frac{b\mathsf b_4}{3\mathsf b_2^2}+O(x^{-2}).
\end{equation}
Thus the integrand in the reconstruction formula satisfies
\begin{equation}\label{eq:critical-log-integrand-expansion}
 \frac{g(x)-1}{xg(x)}
 =-\frac b{\mathsf b_2}x
   +\left(1+\frac{b\mathsf b_4}{3\mathsf b_2^2}\right)\frac1x
   +O(x^{-3}).
\end{equation}
The remainder is integrable at infinity, with tail \(O(x^{-2})\).
Furthermore,
\begin{equation}\label{eq:critical-log-potential-expansion}
 \log g(x)=\log(\mathsf b_2/b)-2\log x+O(x^{-2}).
\end{equation}
Substitution of \eqref{eq:critical-log-integrand-expansion} and
\eqref{eq:critical-log-potential-expansion} into
\eqref{eqt:fixed_point_log_identity} shows that
\begin{equation}\label{eq:critical-renormalized-log-limit}
 \log f(x)+\delta x^2-\nu\log x
       =\log A_\infty+O(x^{-2})
\end{equation}
for a finite real constant \(\log A_\infty\).
Exponentiation proves \eqref{eq:critical-full-tail}, with
\(A_\infty>0\); multiplication by \(-x\) gives
\eqref{eq:critical-vorticity-full-tail}.
\end{proof}
\section{The focusing range and existence of a critical profile}
\label{sec:critical_parameter}

The purpose of this section is to improve the explicit focusing range
and prove that a critical whole-line profile occurs. The analytic input
is an optimal inequality for the first-moment interaction. The
topological input is continuation of the HQWW map on a common compact
convex domain.
\subsection{A sharp first-moment inequality and the focusing range}
\label{subsec:sharp_first_moment}

The first moment of a whole-line profile controls the sign of its spatial
scaling rate.  We first prove an inequality for the underlying shape class;
the profile equation then converts it into an explicit parameter bound.
Throughout this subsection, set
\begin{equation}\label{eq:sharp_moment_constants}
 d:=\frac{18-16\log 2}{35},
 \qquad
 a_{\mathrm{foc}}:=\frac{1-d}{1+d}
 =\frac{17+16\log 2}{53-16\log 2}.
\end{equation}

\begin{lemma}[Sharp first-moment inequality]\label{lem:sharp_first_moment}
Let \(f:[0,\infty)\to[0,\infty)\) be continuous and nonincreasing, with
\(f(0)=1\), \(f(x)\to0\) as \(x\to\infty\), and
\(s\mapsto f(\sqrt{s})\) convex.  Suppose that
\[
 b:=\frac2\pi\int_0^\infty f(x)\,\mathrm{d}x<\infty,
 \qquad
 b_1:=\frac2\pi\int_0^\infty xf(x)\,\mathrm{d}x<\infty.
\]
Define the positive quadratic form
\begin{equation}\label{eq:positive_first_moment_form}
 \mathcal U_1(f):=
 \frac2{\pi^2}\int_0^\infty\int_0^\infty
 \frac{xy}{x+y}f(x)f(y)\,\mathrm{d}x\,\mathrm{d}y.
\end{equation}
Then
\begin{equation}\label{eq:sharp_first_moment_bound}
 \mathcal U_1(f)\leq d\,b b_1.
\end{equation}
Equality holds if and only if
\(f(x)=(1-x^2/L^2)_+\) for some \(L>0\).
\end{lemma}

\begin{proof}
Put \(G(s):=f(\sqrt{s})\), and let \(\sigma_f\) be the distributional
second derivative of \(G\) on \((0,\infty)\).  Since \(G\) is convex,
nonincreasing, and tends to zero, \(\sigma_f\) is a positive measure and
\[
 G(s)=\int_{(s,\infty)}(t-s)\,\sigma_f(\mathrm{d}t),
 \qquad
 \int_{(0,\infty)}t\,\sigma_f(\mathrm{d}t)=1.
\]
Indeed, \(G'_+(s)\uparrow0\) as \(s\to\infty\), and integration of
\(-G'_+(s)=\sigma_f((s,\infty))\) proves the first identity by Tonelli's
theorem.  The second follows by taking \(s\downarrow0\).
The pushforward of \(t\sigma_f(\mathrm{d}t)\) under \(t\mapsto\sqrt t\)
is therefore a probability measure \(\nu_f\) satisfying
\begin{equation}\label{eq:moment_cap_mixture}
 f(x)=\int_{(0,\infty)}\phi_L(x)\,\nu_f(\mathrm{d}L),
 \qquad
 \phi_L(x):=(1-x^2/L^2)_+.
\end{equation}
In particular,
\begin{equation}\label{eq:cap_first_moments}
 b=\frac4{3\pi}\int L\,\nu_f(\mathrm{d}L),
 \qquad
 b_1=\frac1{2\pi}\int L^2\,\nu_f(\mathrm{d}L).
\end{equation}

The kernel in \eqref{eq:positive_first_moment_form} has the factorization
\[
 \frac{xy}{x+y}=\int_0^\infty
 (x\mathrm{e}^{-tx})(y\mathrm{e}^{-ty})\,\mathrm{d}t.
\]
Consequently, its bilinear form satisfies the Cauchy--Schwarz inequality.
Writing \(\mathcal U_1(f,g)\) for the corresponding bilinear form,
homogeneity gives
\[
 \mathcal U_1(\phi_L,\phi_M)
 \leq
 \bigl(\mathcal U_1(\phi_L)\mathcal U_1(\phi_M)\bigr)^{1/2}
 =\mathcal U_1(\phi_1)(LM)^{3/2}.
\]
Using \eqref{eq:moment_cap_mixture}, we integrate this bound with respect
to \(\nu_f\otimes\nu_f\). Cauchy--Schwarz then gives
\begin{equation}\label{eq:mixture_moment_bound}
 \mathcal U_1(f)
 \leq\mathcal U_1(\phi_1)
 \left(\int L^{3/2}\,\nu_f(\mathrm{d}L)\right)^2
 \leq\mathcal U_1(\phi_1)
 \left(\int L\,\nu_f(\mathrm{d}L)\right)
 \left(\int L^2\,\nu_f(\mathrm{d}L)\right).
\end{equation}
All these integrals are finite by \eqref{eq:cap_first_moments}.

To evaluate the remaining constant, symmetry and the change of variables
\(y=tx\) on \(0<y<x<1\) give
\begin{align*}
 K&:=\int_0^1\int_0^1
 \frac{xy}{x+y}(1-x^2)(1-y^2)\,\mathrm{d}y\,\mathrm{d}x\\
 &=2\int_0^1\frac{t}{1+t}
 \left(\frac13-\frac{1+t^2}{5}+\frac{t^2}{7}\right)\,\mathrm{d}t\\
 &=\frac4{15}\int_0^1
 \left(\frac{t}{1+t}-\frac37\frac{t^3}{1+t}\right)\,\mathrm{d}t
 =\frac{18-16\log2}{105}.
\end{align*}
Thus \(\mathcal U_1(\phi_1)=2K/\pi^2\), and
\eqref{eq:cap_first_moments}--\eqref{eq:mixture_moment_bound} yield
\eqref{eq:sharp_first_moment_bound} with \(d=3K\).
If equality holds, equality in the last Cauchy--Schwarz inequality in
\eqref{eq:mixture_moment_bound} requires \(L^{1/2}\) and \(L\) to be
proportional \(\nu_f\)-almost everywhere.  Since \(L>0\), the measure
\(\nu_f\) is concentrated at one point.  Conversely, the preceding exact
calculation shows that every single cap gives equality.
\end{proof}

For a whole-line fixed point with \(\bar\Omega=-xf\), the form
\(\mathcal U_1\) is the Hilbert-transform pairing
\begin{equation}\label{eq:U1_Hilbert_pairing}
 \mathcal U_1(f)=-\frac2\pi\int_0^\infty
 (H_{\mathbb R}\bar\Omega)(x)\bar\Omega(x)\,\mathrm{d}x.
\end{equation}
To justify this identity, first truncate \(\bar\Omega\) symmetrically at a
finite radius and use the principal-value cutoff \(|x-y|>\varepsilon\).
The antisymmetric singular terms cancel after interchanging \(x\) and
\(y\); the remaining numerator satisfies
\[
 \frac{xy^2}{y^2-x^2}+\frac{yx^2}{x^2-y^2}=\frac{xy}{x+y}.
\]
Letting \(\varepsilon\downarrow0\), and then letting the truncation
radius tend to infinity, proves \eqref{eq:U1_Hilbert_pairing}.
The pairing converges by the \(L^2\) boundedness of the Hilbert transform,
whereas the nonnegative double integral converges by monotone convergence.
The required condition \(xf\in L^2\) follows from monotonicity and
\(b_1<\infty\): for \(x>0\),
\[
 \frac{3x^2}{8}f(x)
 \leq\int_{x/2}^x yf(y)\,\mathrm{d}y
 \leq\frac\pi2 b_1,
\]
so \(x^2f(x)^2\) is integrable at infinity; continuity handles the origin.

\begin{proposition}[Improved focusing range]\label{prop:improved_focusing_range}
Let \(0<a<1\), and let \(f\in\mathbb D_{\mathrm H}\) be a fixed point of the
whole-line map of \cite{HQWW2024}.  Set
\[
 \kappa:=\frac{1-a/3}{2}c(f),
 \qquad \bar c_l:=\kappa-a b(f).
\]
Then
\begin{equation}\label{eq:improved_mass_ratio}
 \frac{b(f)}{\kappa}
 <\frac2{(1-d)(1+a)}.
\end{equation}
In particular, every such fixed point is focusing when
\(0<a\leq a_{\mathrm{foc}}\).  Every critical fixed point in this class
satisfies \(a>a_{\mathrm{foc}}\).
\end{proposition}

\begin{proof}
By \cite[Lemma 4.8]{HQWW2024}, \(b_1(f)<\infty\).  The identity
\cite[(4.12), with \(p=1\)]{HQWW2024} reads
\begin{equation}\label{eq:first_moment_mass_identity}
 \bigl((1+a)b(f)-2\kappa\bigr)b_1(f)
 =(1+a)\mathcal U_1(f).
\end{equation}
A parabolic cap is not a fixed point at \(0<a<1\).  Indeed, if a fixed
point has compact support with positive endpoint \(L\), then
\cite[Theorem 4.10]{HQWW2024} gives a vanishing exponent at least
\((3-a)/(2a)>1\), whereas a parabolic cap vanishes linearly.
Lemma~\ref{lem:sharp_first_moment} therefore gives
\(\mathcal U_1(f)<d\,b(f)b_1(f)\).  Since \(b_1(f)>0\),
\eqref{eq:first_moment_mass_identity} proves
\eqref{eq:improved_mass_ratio}.  Consequently,
\[
 \frac{\bar c_l}{\kappa}
 >1-\frac{2a}{(1-d)(1+a)}\geq0
 \qquad\text{for }0<a\leq a_{\mathrm{foc}},
\]
as claimed.
\end{proof}

\begin{remark}\label{rem:focusing_range_scope}
Proposition~\ref{prop:improved_focusing_range} applies to normalized
whole-line fixed points.  Such a focusing fixed point also satisfies the
period-selection rule of the present construction with \(L=\infty\).
It therefore yields a whole-line profile in that construction for every
\(0<a\leq a_{\mathrm{foc}}\).  Finite-period fixed points are governed by the separate period-selection rule of Section~\ref{sec:period_selection}.
\end{remark}
\subsection{Existence of a critical profile}
\label{subsec:critical_continuation}

The parameter estimate supplies opposite spatial rates at the two ends
of an interval. To turn these signs into a critical profile, one must
continue a set of fixed points and control the mass along that set.
The argument does not require a choice of one fixed point at every
parameter.

\begin{theorem}[Critical profile by continuation]
\label{thm:critical_profile_continuation}
There exist
\begin{equation}\label{eq:critical_existence_range}
 a_*\in(a_{\mathrm{foc}},a_{\mathrm{per}}),
 \qquad f_*\in\mathbb D_{\mathrm H},
\end{equation}
such that \(f_*\) is an HQWW fixed point and
\begin{equation}\label{eq:critical_mass_balance}
 \frac{1-a_*/3}{2}c(f_*)=a_*b(f_*),
 \qquad
 a_* = \frac{3c(f_*)}{c(f_*)+6b(f_*)}.
\end{equation}
It is also a fixed point of the period-selection map with
\(\ell(f_*)=\infty\) and \(\bar c_l=0\). Its physical vorticity
\(\bar\Omega_*=-xf_*\) is smooth, belongs to
\(\mathcal S(\mathbb R)\), and has the critical asymptotics of
Section~\ref{sec:regularity}.
\end{theorem}

\begin{proof}
Set \(I=[a_{\mathrm{foc}},a_{\mathrm{per}}]\), and denote the
whole-line HQWW map by \(\mathcal R_a^{\mathbb R}\). This notation
distinguishes its compact-support outputs, which may have negative
spatial rate, from the finite-period outputs of the present map.
On \(\mathbb D_{\mathrm H}\), define
\begin{equation}\label{eq:wholeline_continuation_map}
 G_{a,f}(x)=1+\frac{2aT_\infty f(x)}{(1-a/3)c(f)}.
\end{equation}
The output \(\mathcal R_a^{\mathbb R}f\) is given by
\eqref{eqt:R_map_definition} while \(G_{a,f}>0\), and by zero
after its first zero. The invariance and compactness statements of
\cite[Lemmas 3.6 and 3.10]{HQWW2024} give a common compact convex
domain \(\mathbb D_{\mathrm H}\subset X\) for all \(a\in I\).

\medskip
\noindent\textit{Joint continuity of the map.}
We include the parameter-continuity argument to specify the input to
the continuation theorem. Suppose that \(a_n\to a\) in \(I\) and
\(f_n\to f\) in \(X\). Weighted convergence gives uniform
convergence on bounded intervals. In the integral for \(c\), the
common cap bound controls the origin, while \(y^{-2}\) controls
the tail. Thus \(c(f_n)\to c(f)\), with a positive lower bound
provided by the common slope condition. For fixed bounded \(x\),
the kernel \eqref{eqt:T_whole_line_kernel} is
\(O(x^2/y^2)\) as \(y\to\infty\); its diagonal logarithm is
uniformly integrable. It follows that
\(T_\infty f_n\to T_\infty f\) locally uniformly and hence
\(G_{a_n,f_n}\to G_{a,f}\) locally uniformly.

On every compact subinterval where \(G_{a,f}>0\), the integral
formula converges uniformly. At the origin the supporting-line bound
\(G_{a_n,f_n}(x)\geq1-C_Ix^2\) supplies the integrable majorant
\(C_Ix\) for its logarithmic integrand. At and beyond a terminal
zero one instead uses
\begin{equation}\label{eq:continuation_terminal_control}
 0\leq\mathcal R_{a_n}^{\mathbb R}f_n(x)
 \leq\bigl(G_{a_n,f_n}(x)_+\bigr)^{1/a_n}.
\end{equation}
These observations give pointwise convergence of the outputs.
The common Lipschitz bound of the invariant class gives local uniform
convergence, and the bound \(0\leq\mathcal R_a^{\mathbb R}f\leq1\)
controls the weighted tails in \(X\). Thus
\((a,f)\mapsto\mathcal R_a^{\mathbb R}f\) is jointly continuous.

The continuation theorem in \cite[Theorem 1.2]{Amster2023}, a
normed-space form of Browder's theorem \cite{Browder1960}, now gives
a connected subset whose closure \(\mathscr C\) is a compact connected
subset of
\begin{equation}\label{eq:continuation_fixed_point_set}
 \Sigma:=\{(a,f)\in I\times\mathbb D_{\mathrm H}:
                     f=\mathcal R_a^{\mathbb R}f\}
\end{equation}
whose projection onto the parameter axis is all of \(I\).
The set \(\mathscr C\) may contain several fixed points at the same
parameter; only its connectedness and full projection are used below.

\medskip
\noindent\textit{Continuity of the mass on the fixed-point set.}
Let \((a_n,f_n)\to(a,f)\) in \(\Sigma\). Write
\(b_n=b(f_n)\) and \(\kappa_n=(1-a_n/3)c(f_n)/2\), and use the
corresponding symbols without subscripts for the limit.
These quantities are finite, with \(b_n>\kappa_n\) and \(b>\kappa\),
by Proposition~\ref{thm:whole_line_Q_mu_identity}.
The mass identity \eqref{eqt:whole_line_fixed_point_ratio} is equivalent to
\begin{equation}\label{eq:continuation_mass_quadratic_identity}
 a_n\mathcal Q(f_n)
 =\kappa_n b_n-\frac{b_n^2}{2}
 =\frac{\kappa_n^2-(b_n-\kappa_n)^2}{2}.
\end{equation}
The kernel defining \(\mathcal Q\) in
\eqref{eqt:whole_line_Q_definition} is nonnegative. Fatou's lemma
therefore yields
\begin{equation}\label{eq:continuation_fatou}
 b(f)\leq\liminf_n b_n,
 \qquad
 \mathcal Q(f)\leq\liminf_n\mathcal Q(f_n).
\end{equation}
Since \(\kappa_n\to\kappa\) and \(a_n\to a>0\),
\eqref{eq:continuation_mass_quadratic_identity} gives
\begin{equation}\label{eq:continuation_mass_upper}
 \limsup_n(b_n-\kappa_n)^2
 \leq\kappa^2-2a\mathcal Q(f)=(b-\kappa)^2.
\end{equation}
Since \(b_n-\kappa_n>0\), this gives
\(\limsup_n b_n\leq b(f)\).
Together with \eqref{eq:continuation_fatou}, this proves
\(b_n\to b(f)\). In particular, the spatial rate
\begin{equation}\label{eq:continuation_spatial_rate}
 (a,f)\longmapsto
 \bar c_l(a,f)=\frac{1-a/3}{2}c(f)-ab(f)
\end{equation}
is continuous on \(\Sigma\).
This continuity is a property of the fixed-point set, where the moment
identity \eqref{eq:continuation_mass_quadratic_identity} is available.

\medskip
\noindent\textit{The change of sign.}
At every fixed point with \(a=a_{\mathrm{foc}}\),
Proposition~\ref{prop:improved_focusing_range} gives
\(\bar c_l>0\). At every fixed point with
\(a=a_{\mathrm{per}}=1/(2-\overline\mu)\),
\eqref{eqt:whole_line_rates_from_mu} gives
\begin{equation}\label{eq:continuation_negative_endpoint}
 \bar c_l=\frac b2(1-2a_{\mathrm{per}}+a_{\mathrm{per}}\mu(f))<0,
\end{equation}
because \(\mu(f)<\overline\mu\).
The continuous image of the connected set \(\mathscr C\) under
\eqref{eq:continuation_spatial_rate} therefore contains zero.
Choose \((a_*,f_*)\in\mathscr C\) with \(\bar c_l(a_*,f_*)=0\).
The strict endpoint signs give \eqref{eq:critical_existence_range},
and solving the rate identity gives \eqref{eq:critical_mass_balance}.
Proposition~\ref{prop:hqww_correspondence} identifies this fixed point
with a critical fixed point of the period-selection map.
The regularity and asymptotics follow from
Proposition~\ref{prop:critical_profile_gaussian_rate} and
Corollary~\ref{cor:whole_line_weighted_regularity}.
\end{proof}

\begin{remark}[The critical equation and normalization]
At a critical profile, the equality
\(\kappa=a_*b\) makes the intermediate function
\(g=1+T_\infty f_*/b\). Thus its fixed-point equation is the
critical specialization of the HQWW map, with
\(a_*=3c/(c+6b)\); it is also the critical specialization of our
period-selection map. This is the same normalized critical stationary
problem. The continuation argument proves the existence of a solution
to that problem and the strict interval
\eqref{eq:critical_existence_range}.
\end{remark}

\vspace{0.2in}
\noindent{\bf Acknowledgments.}
This work was supported in part by NSF Grant DMS-2508463,
the Choi Family Gift Fund, and the Mike Yan Gift Fund.

\appendix
\section{Fourier identities for the Hilbert transform}\label{app:hilbert}

This appendix proves the sign convention and inversion identities used in the nonlinear eigenvalue reduction in Section~\ref{sec:steady}.
We work on
\[
    \mathbb T=\R/(2\pi\Z),
\]
and use the Fourier convention
\[
    f(x)=\sum_{k\in\Z}\widehat f(k)e^{ikx}.
\]
The periodic Hilbert transform is normalized by
\[
    \widehat{H f}(k)=-i\,\sgn(k)\widehat f(k),
    \qquad \sgn(0)=0.
\]
The positive square-root of the Laplacian is
\[
    \Lambda=(-\partial^2_{x})^{1/2},
    \qquad \widehat{\Lambda f}(k)=|k|\widehat f(k).
\]
\begin{lemma}[Fourier identities for the periodic Hilbert transform]\label{Appendix: lemma fourier identity for half Laplacian}
Let \(s\in\R\). For every \(f\in H^s(\mathbb T)\),
\[
    H\partial_x f=\Lambda f
\]
as an identity in \(H^{s-1}(\mathbb T)\). Moreover,
\[
    H^2 f=-f+\widehat f(0).
\]
Consequently, on the mean-zero subspace,
\[
    H^2=-I,
    \qquad H^{-1}=-H.
\]
In particular, if
\[
    u_x=H\omega,
\]
then
\[
    \omega-\widehat\omega(0)=-\Lambda u.
\]
If, in addition, \(\widehat\omega(0)=0\), then
\[
    \omega=-\Lambda u.
\]
\end{lemma}

\begin{proof}
Write
\[
    f(x)=\sum_{k\in\Z}\widehat f(k)e^{ikx}.
\]
Then
\[
    \widehat{\partial_x f}(k)=ik\widehat f(k).
\]
Therefore
\[
    \widehat{H\partial_x f}(k)
    =-i\sgn(k)\,ik\widehat f(k)
    =k\sgn(k)\widehat f(k)
    =|k|\widehat f(k),
\]
which is precisely the Fourier multiplier defining \(\Lambda f\). Hence
\[
    H\partial_x f=\Lambda f.
\]

Similarly,
\[
    \widehat{H^2 f}(k)
    =\bigl(-i\sgn(k)\bigr)^2\widehat f(k)
    =-\sgn(k)^2\widehat f(k).
\]
For \(k\ne0\), this equals \(-\widehat f(k)\), while for \(k=0\), it is zero. Therefore
\[
    H^2 f=-f+\widehat f(0).
\]
In particular, if \(f\) has zero mean, then
\[
    H^2 f=-f.
\]
Thus \(H^{-1}=-H\) on the mean-zero subspace.

Now suppose
\[
    u_x=H\omega.
\]
Applying \(H\) to both sides gives
\[
    H u_x=H^2\omega=-\omega+\widehat\omega(0).
\]
On the other hand, the first part of the lemma gives
\[
    H u_x=\Lambda u.
\]
Hence
\[
    \Lambda u=-\omega+\widehat\omega(0),
\]
or equivalently
\[
    \omega-\widehat\omega(0)=-\Lambda u.
\]
If \(\widehat\omega(0)=0\), this reduces to
\[
    \omega=-\Lambda u.
\]
The proof is complete.
\end{proof}

\begin{remark}
The identity \(H\partial_x=\Lambda\) does not require a mean-zero assumption on \(f\). The mean-zero assumption is needed only when inverting \(H\), because constants belong to the kernel of \(H\).
\end{remark}

\bibliographystyle{myalpha}
\bibliography{references}
\end{document}